\documentclass[10pt]{amsart}
\usepackage[cp1251]{inputenc}
\usepackage[english]{babel}
\usepackage{amsmath}
\usepackage{amssymb}
\usepackage{amsfonts}

\begin{document}
\renewcommand{\refname}{References}

\thispagestyle{empty}

\title[Optimization of the Mean-Square Approximation Procedures]
{Optimization of the Mean-Square Approximation Procedures 
for Iterated Ito Stochastic Integrals
of Multiplicities 1 to 5 
from the Unified Taylor--Ito Expansion Based on Multiple 
Fourier--Legendre Series}
\author[M.D. Kuznetsov]{Mikhail D. Kuznetsov}
\address{Mikhail Dmitrievich Kuznetsov
\newline\hphantom{iii} Saint-Petersburg Electrotechnical University,
\newline\hphantom{iii} ul. Professora Popova, 5,
\newline\hphantom{iii} 197376, Saint-Petersburg, Russia}%
\email{sde\_kuznetsov@inbox.ru}
\author[D.F. Kuznetsov]{Dmitriy F. Kuznetsov}
\address{Dmitriy Feliksovich Kuznetsov
\newline\hphantom{iii} Peter the Great Saint-Petersburg Polytechnic University,
\newline\hphantom{iii} Polytechnicheskaya ul., 29,
\newline\hphantom{iii} 195251, Saint-Petersburg, Russia}%
\email{sde\_kuznetsov@inbox.ru}

\vspace{10mm}

\maketitle {\small
\begin{quote}
\noindent{\sc Abstract.} 
The article is devoted to 
optimization of the mean-square approximation procedures 
for iterated Ito stochastic integrals
of multiplicities 1 to 5.
The mentioned stochastic integrals are part 
of strong numerical methods with convergence orders
$1.0,$ $1.5,$ $2.0,$ and $2.5$ 
for Ito stochastic differential equations 
with multidimensional non-commutative noise 
based on the unified Taylor--Ito 
expansion and 
multiple Fourier--Legendre series
converging in the sense of norm in Hilbert space $L_2([t, T]^k)$ 
$(k=1,\ldots,5).$
In this article we use multiple Fourier--Legendre series
within the framework of the method of expansion and mean-square approximation 
of iterated Ito stochastic integrals based on 
generalized multiple Fourier series.
We show that the lengths of sequences of independent 
standard Gaussian random variables required for
the mean-square approximation of iterated Ito stochastic 
integrals of multiplicities 1 to 5 
can be significantly reduced without the loss
of the mean-square accuracy of approximation for 
these stochastic integrals.
\medskip
\end{quote}
}

\vspace{15mm}

\setlength{\baselineskip}{2.0em}

\tableofcontents

\setlength{\baselineskip}{1.2em}

\vspace{7mm}

\section{Explicit One-Step Strong Numerical Methods with Convergence Orders
$1.0,$ $1.5,$ $2.0,$ and $2.5$ 
for Ito SDEs Based on the Unified Taylor--Ito Expansion}

\vspace{5mm}

Let $(\Omega,$ ${\rm F},$ ${\sf P})$ be a complete probability space, let 
$\{{\rm F}_t, t\in[0,T]\}$ be a nondecreasing right-continuous 
family of $\sigma$-algebras of ${\rm F},$
and let ${\bf f}_t$ be a standard $m$-dimensional Wiener stochastic 
process, which is
${\rm F}_t$-measurable for any $t\in[0, T].$ We assume that the components
${\bf f}_{t}^{(i)}$ $(i=1,\ldots,m)$ of this process are independent. 
Consider
an Ito stochastic differential equation (SDE) 
in the integral form

\vspace{-2mm}
\begin{equation}
\label{1.5.2}
{\bf x}_t={\bf x}_0+\int\limits_0^t {\bf a}({\bf x}_{\tau},\tau)d\tau+
\sum\limits_{i=1}^m
\int\limits_0^t B_i({\bf x}_{\tau},\tau)d{\bf f}_{\tau}^{(i)},\ \ \
{\bf x}_0={\bf x}(0,\omega).
\end{equation}

\vspace{3mm}
\noindent
Here ${\bf x}_t$ is some $n$-dimensional stochastic process 
satisfying the Ito SDE (\ref{1.5.2}). 
The nonrandom functions ${\bf a}({\bf x},t): 
\mathbb{R}^n\times[0, T]\to\mathbb{R}^n$,
$B({\bf x},t): \mathbb{R}^n\times[0, T]\to\mathbb{R}^{n\times m}$
guarantee the existence and uniqueness up to stochastic equivalence 
of a solution
of the Ito SDE (\ref{1.5.2}) \cite{1982}. The second integral on 
the right-hand side of (\ref{1.5.2}) is 
interpreted as an Ito stochastic integral and $B_i({\bf x},t)$
is the $i$th colomn of the matrix function $B({\bf x},t).$
Let ${\bf x}_0$ be an $n$-dimensional random variable, which is 
${\rm F}_0$-measurable and 
${\sf M}\bigl\{\left|{\bf x}_0\right|^2\bigr\}<\infty$ 
(${\sf M}$ is an expectation operator).
We assume that
${\bf x}_0$ and ${\bf f}_t-{\bf f}_0$ are independent when $t>0.$

It is well known that one of the promising approaches 
to the numerical integration of 
Ito SDEs is an approach based on the stochastic Taylor expansions
\cite{1988}-\cite{2010}. 
The essential feature of such 
stochastic expansions is a presence in them of the so-called iterated
Ito and Stratonovich stochastic integrals, which play the key 
role for solving the 
problem of numerical integration of Ito SDEs and have the 
following form

\vspace{-1mm}
\begin{equation}
\label{ito}
J[\psi^{(k)}]_{T,t}=\int\limits_t^T\psi_k(t_k) \ldots \int\limits_t^{t_{2}}
\psi_1(t_1) d{\bf w}_{t_1}^{(i_1)}\ldots
d{\bf w}_{t_k}^{(i_k)},
\end{equation}

\begin{equation}
\label{str}
J^{*}[\psi^{(k)}]_{T,t}=
\int\limits_t^{*T}\psi_k(t_k) \ldots \int\limits_t^{*t_{2}}
\psi_1(t_1) d{\bf w}_{t_1}^{(i_1)}\ldots
d{\bf w}_{t_k}^{(i_k)},
\end{equation}

\vspace{2mm}
\noindent
where every $\psi_l(\tau)$ $(l=1,\ldots,k)$ is a 
nonrandom function 
on $[t,T],$ ${\bf w}_{\tau}^{(i)}={\bf f}_{\tau}^{(i)}$
for $i=1,\ldots,m$ and
${\bf w}_{\tau}^{(0)}=\tau,$

\vspace{-3mm}
$$
\int\limits\ \hbox{and}\ \int\limits^{*}
$$ 

\vspace{3mm}
\noindent
denote Ito and 
Stratonovich stochastic integrals,
respectively; $i_1,\ldots,i_k = 0, 1,\ldots,m$
(in this paper, 
we use the definition of the Stratonovich stochastic integral from \cite{1992}).

Note that $\psi_l(\tau)\equiv 1$ $(l=1,\ldots,k)$ and
$i_1,\ldots,i_k = 0, 1,\ldots,m$ in  
the classical Taylor--Ito and Taylor--Stratonovich
expansions
\cite{1992}, \cite{1994}. At the same time 
$\psi_l(\tau)\equiv (t-\tau)^{q_l}$ ($l=1,\ldots,k$; 
$q_1,\ldots,q_k=0, 1, 2,\ldots $) and $i_1,\ldots,i_k = 1,\ldots,m$ in
the unified Taylor--Ito and Taylor--Stratonovich
expansions
\cite{1998a}-\cite{2020xx1}.

Let $C^{2,1}({\mathbb{R}}^n\times [0,T])$ 
be the space of functions $R({\bf x}, t): 
{\mathbb{R}}^n\times [0, T] \to {\mathbb{R}}^1$ with the 
following property: these functions are twice
continuously differentiable in ${\bf x}$ and have one continuous 
derivative in $t$. Let us consider the following differential operators on
the space $C^{2,1}({\mathbb{R}}^n\times [0,T])$

\vspace{1mm}
\begin{equation}
\label{2.3}
L= {\partial \over \partial t}
+ \sum^ {n} _ {i=1} a^{(i)} ({\bf x},  t) 
{\partial  \over  \partial  {\bf  x}^{(i)}} 
+ {1\over 2} \sum^ {m} _ {j=1} \sum^ {n} _ {l,i=1}
B^ {(lj)} ({\bf x}, t) B^ {(ij)} ({\bf x}, t) {\partial
^{2} \over\partial{\bf x}^{(l)}\partial{\bf x}^{(i)}},
\end{equation}

\vspace{2mm}
\begin{equation}
\label{2.4}
G^ {(i)} _ {0} = \sum^ {n} _ {j=1} B^ {(ji)} ({\bf x}, t)
{\partial  \over \partial {\bf x} ^ {(j)}},\ \ \
i=1,\ldots,m,
\end{equation}

\vspace{4mm}
\noindent
where $a^{(i)} ({\bf x},  t)$ is the $i$th component of  
the vector function $a({\bf x},  t)$ and $B^ {(ij)} ({\bf x}, t)$
is the $ij$th element of the matrix function $B({\bf x}, t)$.

Consider the following sequence of differential 
operators

$$
G_p^{(i)}=\frac{1}{p}\left(
G_{p-1}^{(i)}L-LG_{p-1}^{(i)}\right),\ \ \
p=1, 2,\ldots,\ \ \ i=1,\ldots,m,
$$

\vspace{3mm}
\noindent
where $L$ and $G_0^{(i)},$ $i=1,\ldots,m$
are defined by the equalities
(\ref{2.3}), (\ref{2.4}).

For the further consideration, we need to introduce
the following set of iterated Ito stochastic 
integrals
\begin{equation}
\label{qqq1x}
I_{(l_1\ldots l_k)s,t}^{(i_1\ldots i_k)}=
\int\limits_t^s
(t-t_{k})^{l_{k}}\ldots 
\int\limits_t^{t_2}
(t-t _ {1}) ^ {l_ {1}} d
{\bf f} ^ {(i_ {1})} _ {t_ {1}} \ldots 
d {\bf f} _ {t_ {k}} ^ {(i_ {k})},
\end{equation}

\vspace{2mm}
\noindent
where $l_1,\ldots,l_k=0,1,\ldots$ and $i_1,\ldots,l_k=1,\ldots,m.$

Assume that 
$R({\bf x},t)$, 
${\bf a}({\bf x},t),$ and $B_i({\bf x},t),$ $i=1,\ldots,m$
are enough smooth functions with respect to the 
variables ${\bf x}$ and $t$. Then 
for all $s, t\in[0,T]$ such that $s>t$
we can write the following
unified Taylor--Ito expansion
\cite{1998a}-\cite{2020xx1}

\vspace{1mm}
$$
R({\bf x}_s,s)=R({\bf x}_t,t)+\sum_{q=1}^r\sum_{(k,j,l_1,\ldots,l_k)\in
{\rm D}_q}
\frac{(s-t)^j}{j!} \sum\limits_{i_1,\ldots,i_k=1}^m
G_{l_1}^{(i_1)}\ldots G_{l_k}^{(i_k)}
L^j R({\bf x}_t,t)\
I^{(i_1\ldots i_k)}_{(l_1\ldots l_k)s,t}+
$$

\vspace{2mm}
\begin{equation}
\label{15.001}
+\left(H_{r+1}\right)_{s,t}\ \ \ \hbox{w.\ p.\ 1},
\end{equation}

\vspace{5mm}
\noindent
where 
$$
L^j R({\bf x},t)\stackrel{\sf def}{=}
\begin{cases}\underbrace{L\ldots L}_j
R({\bf x},t)\ &\hbox{for}\ j\ge 1\cr\cr
R({\bf x},t)\ &\hbox{for}\ j=0
\end{cases},
$$

\vspace{5mm}
\begin{equation}
\label{asas1}
{\rm D}_{q}=\left\{
(k, j, l_ {1},\ldots, l_ {k}):\  k + 2\left(j +
\sum\limits_{p=1}^k l_p\right)= q;\ k, j, l_{1},\ldots, 
l_{k} =0,1,\ldots\right\},
\end{equation}

\vspace{5mm}
\noindent
and $\left(H_{r+1}\right)_{s,t}$ is the 
remainder term in integral form
\cite{1998a}-\cite{2020xx1}.

Consider the partition 
$\{\tau_p\}_{p=0}^N$ of the interval
$[0,T]$ such that

\vspace{-1mm}
\begin{equation}
\label{asas2}
0=\tau_0<\tau_1<\ldots<\tau_N=T,\ \ \
\Delta_N=
\max\limits_{0\le j\le N-1}\left|\tau_{j+1}-\tau_j\right|.
\end{equation}

\vspace{2mm}

Let ${\bf y}_{\tau_j}\stackrel{\sf def}{=}
{\bf y}_{j},$\ $j=0, 1,\ldots,N$ be a time discrete approximation
of the process ${\bf x}_t,$ $t\in[0,T],$ which is a solution of the Ito
SDE (\ref{1.5.2}). 

\vspace{2mm}

{\bf Definiton 1}\ \cite{1992}.\
{\it We will say that a time discrete approximation 
${\bf y}_{j},$\ $j=0, 1,\ldots,N,$
corresponding to the maximal step of discretization $\Delta_N,$
converges strongly with order
$\gamma>0$ at time moment 
$T$ to the process ${\bf x}_t,$ $t\in[0,T]$,
if there exists a constant $C>0,$ which does not depend on 
$\Delta_N,$ and a $\delta>0$ such that 

\vspace{-2mm}
$$
{\sf M}\{|{\bf x}_T-{\bf y}_T|\}\le
C(\Delta_N)^{\gamma}
$$

\vspace{2mm}
\noindent
for each $\Delta_N\in(0, \delta).$}

\vspace{2mm}

From (\ref{15.001}) for $s=\tau_{p+1}$ and
$t=\tau_p$ we obtain the following representation for family of
explicit one-step strong numerical schemes 
for the Ito SDE (\ref{1.5.2})

\vspace{1mm}

$$
{\bf y}_{p+1}={\bf y}_{p}+
\sum_{q=1}^r\sum_{(k,j,l_1,\ldots,l_k)\in{\rm D}_q}
\frac{(\tau_{p+1}-\tau_p)^j}{j!} \sum\limits_{i_1,\ldots,i_k=1}^m
G_{l_1}^{(i_1)}\ldots G_{l_k}^{(i_k)}
L^j\hspace{0.4mm}{\bf y}_{p}\
{\hat I}^{(i_1\ldots i_k)}_{(l_1\ldots l_k)\tau_{p+1},\tau_p}+
$$

\vspace{2mm}
\begin{equation}
\label{15.002}
+{\bf 1}_{\{r=2d-1,
d\in \mathbb{N}\}}\frac{(\tau_{p+1}-\tau_p)^{(r+1)/2}}{\left(
(r+1)/2\right)!}L^{(r+1)/2}\hspace{0.4mm}{\bf y}_{p},
\end{equation}

\vspace{5mm}
\noindent
where
$\hat I^{(i_1\ldots i_k)}_{(l_1\ldots l_k)\tau_{p+1},\tau_p}$ 
is an approximation of the iterated Ito
stochastic integral 
(\ref{qqq1x}).
The equality (\ref{15.002}) should be understood componentwise
with respect to the components ${\bf y}_p^{(i)}$ of the column
${\bf y}_p.$
Let for simplicity 
$\tau_p=p\Delta$,
$\Delta=T/N,$ $T=\tau_N,$ $p=0,1,\ldots,N.$

It is known \cite{1992} that under the standard conditions
the numerical scheme (\ref{15.002}) has strong order of convergence $r/2$
($r\in\mathbb{N}$).

Further, we consider particular cases of the numerical scheme
(\ref{15.002}) for $r=2,3,4,$ and $5,$ i.e. 
explicit one-step strong numerical schemes 
with convergence orders $1.0, 1.5, 2.0,$ and $2.5$ 
for the Ito SDE (\ref{1.5.2}).
At that for simplicity 
we will write ${\bf a},$ $L {\bf a},$ 
$B_i,$ $G_0^{(i)}B_{j}$ etc.
instead of ${\bf a}({\bf y}_p,\tau_p),$ 
$L {\bf a}({\bf y}_p,\tau_p),$ 
$B_i({\bf y}_p,\tau_p),$ 
$G_0^{(i)}B_{j}({\bf y}_p,\tau_p)$ etc. correspondingly.
Here $L$ and $G_0^{(i)},$ $i=1,\ldots,m$
are defined by the equalities
(\ref{2.3}), (\ref{2.4}). Thus, we obtain
the following numerical schemes.

\vspace{8mm}

\centerline{\bf Milstein scheme (Scheme with strong order 1.0)}

\vspace{3mm}

\begin{equation}
\label{al1}
{\bf y}_{p+1}={\bf y}_{p}+\sum_{i_{1}=1}^{m}B_{i_{1}}
\hat I_{(0)\tau_{p+1},\tau_p}^{(i_{1})}+\Delta{\bf a}
+\sum_{i_{1},i_{2}=1}^{m}G_0^{(i_{1})}
B_{i_{2}}\hat I_{(00)\tau_{p+1},\tau_p}^{(i_{1}i_{2})}.
\end{equation}

\vspace{7mm}

\centerline{\bf Scheme with strong order 1.5}

\vspace{6mm}
$$
{\bf y}_{p+1}={\bf y}_{p}+\sum_{i_{1}=1}^{m}B_{i_{1}}
\hat I_{(0)\tau_{p+1},\tau_p}^{(i_{1})}+\Delta{\bf a}
+\sum_{i_{1},i_{2}=1}^{m}G_0^{(i_{1})}
B_{i_{2}}\hat I_{(00)\tau_{p+1},\tau_p}^{(i_{1}i_{2})}+
$$

\vspace{2mm}
$$
+
\sum_{i_{1}=1}^{m}\left[G_0^{(i_{1})}{\bf a}\left(
\Delta \hat I_{(0)\tau_{p+1},\tau_p}^{(i_{1})}+
\hat I_{(1)\tau_{p+1},\tau_p}^{(i_{1})}\right)
- LB_{i_{1}}\hat I_{(1)\tau_{p+1},\tau_p}^{(i_{1})}\right]+
$$

\vspace{2mm}
\begin{equation}
\label{al2}
+\sum_{i_{1},i_{2},i_{3}=1}^{m} G_0^{(i_{1})}G_0^{(i_{2})}
B_{i_{3}}\hat I_{(000)\tau_{p+1},\tau_p}^{(i_{1}i_{2}i_{3})}
+
\frac{\Delta^2}{2}L{\bf a}.
\end{equation}

\vspace{12mm}

\centerline{\bf Scheme with strong order 2.0}

\vspace{6mm}
$$
{\bf y}_{p+1}={\bf y}_{p}+\sum_{i_{1}=1}^{m}B_{i_{1}}
\hat I_{(0)\tau_{p+1},\tau_p}^{(i_{1})}+\Delta{\bf a}
+\sum_{i_{1},i_{2}=1}^{m}G_0^{(i_{1})}
B_{i_{2}}\hat I_{(00)\tau_{p+1},\tau_p}^{(i_{1}i_{2})}+
$$

\vspace{2mm}
$$
+
\sum_{i_{1}=1}^{m}\left[G_0^{(i_{1})}{\bf a}\left(
\Delta \hat I_{(0)\tau_{p+1},\tau_p}^{(i_{1})}+
\hat I_{(1)\tau_{p+1},\tau_p}^{(i_{1})}\right)
- LB_{i_{1}}\hat I_{(1)\tau_{p+1},\tau_p}^{(i_{1})}\right]+
$$

\vspace{2mm}
$$
+\sum_{i_{1},i_{2},i_{3}=1}^{m} G_0^{(i_{1})}G_0^{(i_{2})}
B_{i_{3}}\hat I_{(000)\tau_{p+1},\tau_p}^{(i_{1}i_{2}i_{3})}+
\frac{\Delta^2}{2} L{\bf a}+
$$

\vspace{2mm}
$$
+\sum_{i_{1},i_{2}=1}^{m}
\left[G_0^{(i_{1})} LB_{i_{2}}\left(
\hat I_{(10)\tau_{p+1},\tau_p}^{(i_{1}i_{2})}-
\hat I_{(01)\tau_{p+1},\tau_p}^{(i_{1}i_{2})}
\right)
-LG_0^{(i_{1})}
B_{i_{2}}\hat I_{(10)\tau_{p+1},\tau_p}^{(i_{1}i_{2})}
+\right.
$$

\vspace{2mm}
$$
\left.+G_0^{(i_{1})}G_0^{(i_{2})}{\bf a}\left(
\hat I_{(01)\tau_{p+1},\tau_p}
^{(i_{1}i_{2})}+\Delta \hat I_{(00)\tau_{p+1},\tau_p}^{(i_{1}i_{2})}
\right)\right]+
$$

\vspace{2mm}
\begin{equation}
\label{al3}
+
\sum_{i_{1},i_{2},i_{3},i_{4}=1}^{m}G_0^{(i_{1})}G_0^{(i_{2})}G_0^{(i_{3})}
B_{i_{4}}\hat I_{(0000)\tau_{p+1},\tau_p}^{(i_{1}i_{2}i_{3}i_{4})}.
\end{equation}

\vspace{12mm}

\centerline{\bf Scheme with strong order 2.5}

\vspace{6mm}
$$
{\bf y}_{p+1}={\bf y}_{p}+\sum_{i_{1}=1}^{m}B_{i_{1}}
\hat I_{(0)\tau_{p+1},\tau_p}^{(i_{1})}+\Delta{\bf a}
+\sum_{i_{1},i_{2}=1}^{m}G_0^{(i_{1})}
B_{i_{2}}\hat I_{(00)\tau_{p+1},\tau_p}^{(i_{1}i_{2})}+
$$

\vspace{2mm}
$$
+
\sum_{i_{1}=1}^{m}\left[G_0^{(i_{1})}{\bf a}\left(
\Delta \hat I_{(0)\tau_{p+1},\tau_p}^{(i_{1})}+
\hat I_{(1)\tau_{p+1},\tau_p}^{(i_{1})}\right)
- LB_{i_{1}}\hat I_{(1)\tau_{p+1},\tau_p}^{(i_{1})}\right]+
$$

\vspace{2mm}
$$
+\sum_{i_{1},i_{2},i_{3}=1}^{m} G_0^{(i_{1})}G_0^{(i_{2})}
B_{i_{3}}\hat I_{(000)\tau_{p+1},\tau_p}^{(i_{1}i_{2}i_{3})}+
\frac{\Delta^2}{2} L{\bf a}+
$$

\vspace{2mm}
$$
+\sum_{i_{1},i_{2}=1}^{m}
\left[G_0^{(i_{1})} LB_{i_{2}}\left(
\hat I_{(10)\tau_{p+1},\tau_p}^{(i_{1}i_{2})}-
\hat I_{(01)\tau_{p+1},\tau_p}^{(i_{1}i_{2})}
\right)
- LG_0^{(i_{1})}
B_{i_{2}}\hat I_{(10)\tau_{p+1},\tau_p}^{(i_{1}i_{2})}
+\right.
$$

\vspace{2mm}
$$
\left.+G_0^{(i_{1})}G_0^{(i_{2})}{\bf a}\left(
\hat I_{(01)\tau_{p+1},\tau_p}
^{(i_{1}i_{2})}+\Delta \hat I_{(00)\tau_{p+1},\tau_p}^{(i_{1}i_{2})}
\right)\right]+
$$

$$
+
\sum_{i_{1},i_{2},i_{3},i_{4}=1}^{m}G_0^{(i_{1})}G_0^{(i_{2})}G_0^{(i_{3})}
B_{i_{4}}\hat I_{(0000)\tau_{p+1},\tau_p}^{(i_{1}i_{2}i_{3}i_{4})}+
$$

\vspace{2mm}
$$
+\sum_{i_{1}=1}^{m}\Biggl[G_0^{(i_{1})} L{\bf a}\left(\frac{1}{2}
\hat I_{(2)\tau_{p+1},\tau_p}
^{(i_{1})}+\Delta \hat I_{(1)\tau_{p+1},\tau_p}^{(i_{1})}+
\frac{\Delta^2}{2}\hat I_{(0)\tau_{p+1},\tau_p}^{(i_{1})}\right)\Biggr.+
$$

\vspace{2mm}
$$
+\frac{1}{2} L LB_{i_{1}}\hat I_{(2)\tau_{p+1},\tau_p}^{(i_{1})}-
LG_0^{(i_{1})}{\bf a}\Biggl.
\left(\hat I_{(2)\tau_{p+1},\tau_p}^{(i_{1})}+
\Delta \hat I_{(1)\tau_{p+1},\tau_p}^{(i_{1})}\right)\Biggr]+
$$

\vspace{2mm}
$$
+
\sum_{i_{1},i_{2},i_{3}=1}^m\left[
G_0^{(i_{1})} LG_0^{(i_{2})}B_{i_{3}}
\left(\hat I_{(100)\tau_{p+1},\tau_p}
^{(i_{1}i_{2}i_{3})}-\hat I_{(010)\tau_{p+1},\tau_p}
^{(i_{1}i_{2}i_{3})}\right)
\right.+
$$

\vspace{2mm}
$$
+G_0^{(i_{1})}G_0^{(i_{2})} LB_{i_{3}}\left(
\hat I_{(010)\tau_{p+1},\tau_p}^{(i_{1}i_{2}i_{3})}-
\hat I_{(001)\tau_{p+1},\tau_p}^{(i_{1}i_{2}i_{3})}\right)+
$$

\vspace{2mm}
$$
+G_0^{(i_{1})}G_0^{(i_{2})}G_0^{(i_{3})} {\bf a}
\left(\Delta \hat I_{(000)\tau_{p+1},\tau_p}^{(i_{1}i_{2}i_{3})}+
\hat I_{(001)\tau_{p+1},\tau_p}^{(i_{1}i_{2}i_{3})}\right)-
$$

\vspace{2mm}

$$
\left.- LG_0^{(i_{1})}G_0^{(i_{2})}B_{i_{3}}
\hat I_{(100)\tau_{p+1},\tau_p}^{(i_{1}i_{2}i_{3})}\right]+
$$

\vspace{2mm}
$$
+\sum_{i_{1},i_{2},i_{3},i_{4},i_{5}=1}^m
G_0^{(i_{1})}G_0^{(i_{2})}G_0^{(i_{3})}G_0^{(i_{4})}B_{i_{5}}
\hat I_{(00000)\tau_{p+1},\tau_p}^{(i_{1}i_{2}i_{3}i_{4}i_{5})}+
$$

\vspace{1mm}
\begin{equation}
\label{al4}
+
\frac{\Delta^3}{6}LL{\bf a}.
\end{equation}

\vspace{6mm}

It is well known \cite{1992} that under the standard conditions
the numerical schemes (\ref{al2})--(\ref{al4}) 
have strong orders of convergence 1.0, 1.5, 2.0, and 2.5
correspondingly.
Among these conditions we consider only the condition
for approximations of iterated Ito stochastic 
integrals from the numerical
schemes (\ref{al2})--(\ref{al4}) \cite{1992} (also see \cite{2006})

\begin{equation}
\label{uslov}
{\sf M}\left\{\Biggl(I_{(l_{1}\ldots l_{k})\tau_{p+1},\tau_p}
^{(i_{1}\ldots i_{k})} 
-\hat I_{(l_{1}\ldots l_{k})\tau_{p+1},\tau_p}^{(i_{1}\ldots i_{k})}
\Biggr)^2\right\}\le C\Delta^{r+1},
\end{equation}

\vspace{4mm}
\noindent
where 
constant $C$ is independent of $\Delta$ and
$r/2$ is the strong order of convergence for the numerical schemes
(\ref{al2})--(\ref{al4}), i.e. $r/2=1.0, 1.5, 2.0,$ and $2.5.$

Note that
the numerical schemes (\ref{al1})--(\ref{al4})
are unrealizable in practice without 
effective procedures for the numerical simulation 
of iterated Ito stochastic integrals
from (\ref{15.002}).
That is why in the next section, 
we consider the effective method
of the mean-square approximation of
iterated Ito and Stratonovich stochastic integrals
of arbitrary multiplicity $k$ ($k\in\mathbb{N}$).

\vspace{5mm}

\section{Method of Expansion and Mean-Square Approximation of
Iterated Ito and Stratonovich Stochastic Integrals
Based on Generalized Multiple Fourier Series}

\vspace{5mm}

Let us consider the effective approach to expansion 
and mean-square approximation of iterated Ito
stochastic integrals \cite{2006} (2006), \cite{2007a}-\cite{2020xx1}, 
\cite{5a}-\cite{Mikh-1} 
(the so-called
method of generalized
multiple Fourier series).

The idea of this method is as follows: 
the iterated Ito stochastic 
integral (\ref{ito}) of the multiplicity $k$ ($k\in\mathbb{N}$)
is represented as the multiple stochastic 
integral from the certain discontinuous nonrandom function of $k$ variables 
defined on the hypercube $[t, T]^k.$ Here $[t, T]$ is the interval of 
integration of the iterated Ito stochastic integral (\ref{ito}). 
Then, the mentioned
nonrandom function of $k$ variables
is expanded in the hypercube $[t, T]^k$ into the generalized 
multiple Fourier series converging 
in the mean-square sense
in the space 
$L_2([t,T]^k)$. After a number of nontrivial transformations we come 
to the 
mean-square converging expansion of the iterated Ito stochastic 
integral (\ref{ito}) into the multiple 
series of products
of standard  Gaussian random 
variables. The coefficients of this 
series are the coefficients of 
generalized multiple Fourier series for the mentioned nonrandom function 
of $k$ variables, which can be calculated using the explicit formula 
regardless 
of the multiplicity $k$ of the iterated Ito stochastic integral (\ref{ito}).

Suppose that $\psi_1(\tau),\ldots,\psi_k(\tau)\in L_2([t, T])$.
Define the following function on the hypercube $[t, T]^k$

\vspace{-1mm}
\begin{equation}
\label{ppp}
K(t_1,\ldots,t_k)=
\begin{cases}
\psi_1(t_1)\ldots \psi_k(t_k)\ \ \hbox{for}\ \ t_1<\ldots<t_k\cr\cr\cr
0\ \ \hbox{otherwise}
\end{cases},\ \ t_1,\ldots,t_k\in[t, T],\ \ k\ge 2,
\end{equation}

\vspace{4mm}
\noindent
and 
$K(t_1)\equiv\psi_1(t_1)$ for $t_1\in[t, T].$

Suppose that $\{\phi_j(x)\}_{j=0}^{\infty}$
is a complete orthonormal system of functions in the space
$L_2([t, T])$. 
The function $K(t_1,\ldots,t_k)$ belongs to the space $L_2([t, T]^k).$
At this situation it is well known that the generalized 
multiple Fourier series 
of $K(t_1,\ldots,t_k)\in L_2([t, T]^k)$ is converging 
to $K(t_1,\ldots,t_k)$ in the hypercube $[t, T]^k$ in 
the mean-square sense, i.e.

\vspace{1mm}
$$
\hbox{\vtop{\offinterlineskip\halign{
\hfil#\hfil\cr
{\rm lim}\cr
$\stackrel{}{{}_{p_1,\ldots,p_k\to \infty}}$\cr
}} }\Biggl\Vert
K(t_1,\ldots,t_k)-
\sum_{j_1=0}^{p_1}\ldots \sum_{j_k=0}^{p_k}
C_{j_k\ldots j_1}\prod_{l=1}^{k} \phi_{j_l}(t_l)\Biggr\Vert_{L_2([t,T]^k)}=0,
$$

\vspace{3mm}
\noindent
where

\vspace{-1mm}
\begin{equation}
\label{ppppa}
C_{j_k\ldots j_1}=\int\limits_{[t,T]^k}
K(t_1,\ldots,t_k)\prod_{l=1}^{k}\phi_{j_l}(t_l)dt_1\ldots dt_k
\end{equation}

\vspace{4mm}
\noindent
is the Fourier coefficient,

$$
\left\Vert f\right\Vert_{L_2([t,T]^k)}=\left(\int\limits_{[t,T]^k}
f^2(t_1,\ldots,t_k)dt_1\ldots dt_k\right)^{1/2}.
$$

\vspace{4mm}

Consider the partition $\{\tau_j\}_{j=0}^N$ of the interval
$[t,T]$ such that

\begin{equation}
\label{1111}
t=\tau_0<\ldots <\tau_N=T,\ \ \
\Delta_N=
\hbox{\vtop{\offinterlineskip\halign{
\hfil#\hfil\cr
{\rm max}\cr
$\stackrel{}{{}_{0\le j\le N-1}}$\cr
}} }\Delta\tau_j\to 0\ \ \hbox{if}\ \ N\to \infty,\ \ \
\Delta\tau_j=\tau_{j+1}-\tau_j.
\end{equation}

\vspace{5mm}

{\bf Theorem 1} \cite{2006} (2006), \cite{2007a}-\cite{2020xx1}, 
\cite{5a}-\cite{Mikh-1}.
{\it Suppose that
every $\psi_l(\tau)$ $(l=1,\ldots, k)$ is a continuous nonrandom function on 
$[t, T]$ and
$\{\phi_j(x)\}_{j=0}^{\infty}$ is a complete orthonormal system  
of continuous functions in the space $L_2([t,T]).$ Then

\vspace{1mm}
$$
J[\psi^{(k)}]_{T,t}\  =\ 
\hbox{\vtop{\offinterlineskip\halign{
\hfil#\hfil\cr
{\rm l.i.m.}\cr
$\stackrel{}{{}_{p_1,\ldots,p_k\to \infty}}$\cr
}} }\sum_{j_1=0}^{p_1}\ldots\sum_{j_k=0}^{p_k}
C_{j_k\ldots j_1}\Biggl(
\prod_{l=1}^k\zeta_{j_l}^{(i_l)}\ -
\Biggr.
$$

\vspace{4mm}
\begin{equation}
\label{tyyy}
-\ \Biggl.
\hbox{\vtop{\offinterlineskip\halign{
\hfil#\hfil\cr
{\rm l.i.m.}\cr
$\stackrel{}{{}_{N\to \infty}}$\cr
}} }\sum_{(l_1,\ldots,l_k)\in {\rm G}_k}
\phi_{j_{1}}(\tau_{l_1})
\Delta{\bf w}_{\tau_{l_1}}^{(i_1)}\ldots
\phi_{j_{k}}(\tau_{l_k})
\Delta{\bf w}_{\tau_{l_k}}^{(i_k)}\Biggr),
\end{equation}

\vspace{6mm}
\noindent
where $J[\psi^{(k)}]_{T,t}$ is defined by {\rm (\ref{ito}),}

\vspace{1mm}
$$
{\rm G}_k={\rm H}_k\backslash{\rm L}_k,\ \ \
{\rm H}_k=\{(l_1,\ldots,l_k):\ l_1,\ldots,l_k=0,\ 1,\ldots,N-1\},
$$

\vspace{1mm}
$$
{\rm L}_k=\{(l_1,\ldots,l_k):\ l_1,\ldots,l_k=0,\ 1,\ldots,N-1;\
l_g\ne l_r\ (g\ne r);\ g, r=1,\ldots,k\},
$$

\vspace{5mm}
\noindent
${\rm l.i.m.}$ is a limit in the mean-square sense$,$
$i_1,\ldots,i_k=0,1,\ldots,m,$

\vspace{-1mm}
\begin{equation}
\label{rr23}
\zeta_{j}^{(i)}=
\int\limits_t^T \phi_{j}(s) d{\bf w}_s^{(i)}
\end{equation} 

\vspace{3mm}
\noindent
are independent standard Gaussian random variables
for various
$i$ or $j$ {\rm(}if $i\ne 0${\rm),}
$C_{j_k\ldots j_1}$ is the Fourier coefficient {\rm(\ref{ppppa}),}
$\Delta{\bf w}_{\tau_{j}}^{(i)}=
{\bf w}_{\tau_{j+1}}^{(i)}-{\bf w}_{\tau_{j}}^{(i)}$
$(i=0, 1,\ldots,m),$
$\left\{\tau_{j}\right\}_{j=0}^{N}$ is a partition of
the interval $[t, T],$ which satisfies the condition {\rm (\ref{1111})}.
}

\vspace{4mm}

Note that the adaptation of Theorem 1 
for complete orthonormal systems of Haar and 
Rademacher--Walsh functions in the space $L_2([t,T])$ can be found in
\cite{2006}-\cite{2020xx1},
\cite{5a}-\cite{7}, \cite{9a}-\cite{11}.

In order to evaluate the significance of Theorem 1 for practice we will
demonstrate its transformed particular cases for 
$k=1,\ldots,6$ 
\cite{2006}-\cite{2020xx1}, \cite{5a}-\cite{Mikh-1}

\vspace{1mm}
\begin{equation}
\label{a1}
J[\psi^{(1)}]_{T,t}
=\hbox{\vtop{\offinterlineskip\halign{
\hfil#\hfil\cr
{\rm l.i.m.}\cr
$\stackrel{}{{}_{p_1\to \infty}}$\cr
}} }\sum_{j_1=0}^{p_1}
C_{j_1}\zeta_{j_1}^{(i_1)},
\end{equation}

\vspace{4mm}
\begin{equation}
\label{leto5001}
J[\psi^{(2)}]_{T,t}
=\hbox{\vtop{\offinterlineskip\halign{
\hfil#\hfil\cr
{\rm l.i.m.}\cr
$\stackrel{}{{}_{p_1,p_2\to \infty}}$\cr
}} }\sum_{j_1=0}^{p_1}\sum_{j_2=0}^{p_2}
C_{j_2j_1}\Biggl(\zeta_{j_1}^{(i_1)}\zeta_{j_2}^{(i_2)}
-{\bf 1}_{\{i_1=i_2\ne 0\}}
{\bf 1}_{\{j_1=j_2\}}\Biggr),
\end{equation}

\vspace{6mm}
$$
J[\psi^{(3)}]_{T,t}=
\hbox{\vtop{\offinterlineskip\halign{
\hfil#\hfil\cr
{\rm l.i.m.}\cr
$\stackrel{}{{}_{p_1,\ldots,p_3\to \infty}}$\cr
}} }\sum_{j_1=0}^{p_1}\sum_{j_2=0}^{p_2}\sum_{j_3=0}^{p_3}
C_{j_3j_2j_1}\Biggl(
\zeta_{j_1}^{(i_1)}\zeta_{j_2}^{(i_2)}\zeta_{j_3}^{(i_3)}
-\Biggr.
$$

\begin{equation}
\label{leto5002}
\Biggl.-{\bf 1}_{\{i_1=i_2\ne 0\}}
{\bf 1}_{\{j_1=j_2\}}
\zeta_{j_3}^{(i_3)}
-{\bf 1}_{\{i_2=i_3\ne 0\}}
{\bf 1}_{\{j_2=j_3\}}
\zeta_{j_1}^{(i_1)}-
{\bf 1}_{\{i_1=i_3\ne 0\}}
{\bf 1}_{\{j_1=j_3\}}
\zeta_{j_2}^{(i_2)}\Biggr),
\end{equation}

\vspace{7mm}
$$
J[\psi^{(4)}]_{T,t}
=
\hbox{\vtop{\offinterlineskip\halign{
\hfil#\hfil\cr
{\rm l.i.m.}\cr
$\stackrel{}{{}_{p_1,\ldots,p_4\to \infty}}$\cr
}} }\sum_{j_1=0}^{p_1}\ldots\sum_{j_4=0}^{p_4}
C_{j_4\ldots j_1}\Biggl(
\prod_{l=1}^4\zeta_{j_l}^{(i_l)}
\Biggr.
-
$$
$$
-
{\bf 1}_{\{i_1=i_2\ne 0\}}
{\bf 1}_{\{j_1=j_2\}}
\zeta_{j_3}^{(i_3)}
\zeta_{j_4}^{(i_4)}
-
{\bf 1}_{\{i_1=i_3\ne 0\}}
{\bf 1}_{\{j_1=j_3\}}
\zeta_{j_2}^{(i_2)}
\zeta_{j_4}^{(i_4)}-
$$
$$
-
{\bf 1}_{\{i_1=i_4\ne 0\}}
{\bf 1}_{\{j_1=j_4\}}
\zeta_{j_2}^{(i_2)}
\zeta_{j_3}^{(i_3)}
-
{\bf 1}_{\{i_2=i_3\ne 0\}}
{\bf 1}_{\{j_2=j_3\}}
\zeta_{j_1}^{(i_1)}
\zeta_{j_4}^{(i_4)}-
$$
$$
-
{\bf 1}_{\{i_2=i_4\ne 0\}}
{\bf 1}_{\{j_2=j_4\}}
\zeta_{j_1}^{(i_1)}
\zeta_{j_3}^{(i_3)}
-
{\bf 1}_{\{i_3=i_4\ne 0\}}
{\bf 1}_{\{j_3=j_4\}}
\zeta_{j_1}^{(i_1)}
\zeta_{j_2}^{(i_2)}+
$$
$$
+
{\bf 1}_{\{i_1=i_2\ne 0\}}
{\bf 1}_{\{j_1=j_2\}}
{\bf 1}_{\{i_3=i_4\ne 0\}}
{\bf 1}_{\{j_3=j_4\}}
+
$$
$$
+
{\bf 1}_{\{i_1=i_3\ne 0\}}
{\bf 1}_{\{j_1=j_3\}}
{\bf 1}_{\{i_2=i_4\ne 0\}}
{\bf 1}_{\{j_2=j_4\}}+
$$
\begin{equation}
\label{leto5003}
+\Biggl.
{\bf 1}_{\{i_1=i_4\ne 0\}}
{\bf 1}_{\{j_1=j_4\}}
{\bf 1}_{\{i_2=i_3\ne 0\}}
{\bf 1}_{\{j_2=j_3\}}\Biggr),
\end{equation}

\vspace{10mm}

$$
J[\psi^{(5)}]_{T,t}
=\hbox{\vtop{\offinterlineskip\halign{
\hfil#\hfil\cr
{\rm l.i.m.}\cr
$\stackrel{}{{}_{p_1,\ldots,p_5\to \infty}}$\cr
}} }\sum_{j_1=0}^{p_1}\ldots\sum_{j_5=0}^{p_5}
C_{j_5\ldots j_1}\Biggl(
\prod_{l=1}^5\zeta_{j_l}^{(i_l)}
-\Biggr.
$$
$$
-
{\bf 1}_{\{i_1=i_2\ne 0\}}
{\bf 1}_{\{j_1=j_2\}}
\zeta_{j_3}^{(i_3)}
\zeta_{j_4}^{(i_4)}
\zeta_{j_5}^{(i_5)}-
{\bf 1}_{\{i_1=i_3\ne 0\}}
{\bf 1}_{\{j_1=j_3\}}
\zeta_{j_2}^{(i_2)}
\zeta_{j_4}^{(i_4)}
\zeta_{j_5}^{(i_5)}-
$$
$$
-
{\bf 1}_{\{i_1=i_4\ne 0\}}
{\bf 1}_{\{j_1=j_4\}}
\zeta_{j_2}^{(i_2)}
\zeta_{j_3}^{(i_3)}
\zeta_{j_5}^{(i_5)}-
{\bf 1}_{\{i_1=i_5\ne 0\}}
{\bf 1}_{\{j_1=j_5\}}
\zeta_{j_2}^{(i_2)}
\zeta_{j_3}^{(i_3)}
\zeta_{j_4}^{(i_4)}-
$$
$$
-
{\bf 1}_{\{i_2=i_3\ne 0\}}
{\bf 1}_{\{j_2=j_3\}}
\zeta_{j_1}^{(i_1)}
\zeta_{j_4}^{(i_4)}
\zeta_{j_5}^{(i_5)}-
{\bf 1}_{\{i_2=i_4\ne 0\}}
{\bf 1}_{\{j_2=j_4\}}
\zeta_{j_1}^{(i_1)}
\zeta_{j_3}^{(i_3)}
\zeta_{j_5}^{(i_5)}-
$$
$$
-
{\bf 1}_{\{i_2=i_5\ne 0\}}
{\bf 1}_{\{j_2=j_5\}}
\zeta_{j_1}^{(i_1)}
\zeta_{j_3}^{(i_3)}
\zeta_{j_4}^{(i_4)}
-{\bf 1}_{\{i_3=i_4\ne 0\}}
{\bf 1}_{\{j_3=j_4\}}
\zeta_{j_1}^{(i_1)}
\zeta_{j_2}^{(i_2)}
\zeta_{j_5}^{(i_5)}-
$$
$$
-
{\bf 1}_{\{i_3=i_5\ne 0\}}
{\bf 1}_{\{j_3=j_5\}}
\zeta_{j_1}^{(i_1)}
\zeta_{j_2}^{(i_2)}
\zeta_{j_4}^{(i_4)}
-{\bf 1}_{\{i_4=i_5\ne 0\}}
{\bf 1}_{\{j_4=j_5\}}
\zeta_{j_1}^{(i_1)}
\zeta_{j_2}^{(i_2)}
\zeta_{j_3}^{(i_3)}+
$$
$$
+
{\bf 1}_{\{i_1=i_2\ne 0\}}
{\bf 1}_{\{j_1=j_2\}}
{\bf 1}_{\{i_3=i_4\ne 0\}}
{\bf 1}_{\{j_3=j_4\}}\zeta_{j_5}^{(i_5)}+
{\bf 1}_{\{i_1=i_2\ne 0\}}
{\bf 1}_{\{j_1=j_2\}}
{\bf 1}_{\{i_3=i_5\ne 0\}}
{\bf 1}_{\{j_3=j_5\}}\zeta_{j_4}^{(i_4)}+
$$
$$
+
{\bf 1}_{\{i_1=i_2\ne 0\}}
{\bf 1}_{\{j_1=j_2\}}
{\bf 1}_{\{i_4=i_5\ne 0\}}
{\bf 1}_{\{j_4=j_5\}}\zeta_{j_3}^{(i_3)}+
{\bf 1}_{\{i_1=i_3\ne 0\}}
{\bf 1}_{\{j_1=j_3\}}
{\bf 1}_{\{i_2=i_4\ne 0\}}
{\bf 1}_{\{j_2=j_4\}}\zeta_{j_5}^{(i_5)}+
$$
$$
+
{\bf 1}_{\{i_1=i_3\ne 0\}}
{\bf 1}_{\{j_1=j_3\}}
{\bf 1}_{\{i_2=i_5\ne 0\}}
{\bf 1}_{\{j_2=j_5\}}\zeta_{j_4}^{(i_4)}+
{\bf 1}_{\{i_1=i_3\ne 0\}}
{\bf 1}_{\{j_1=j_3\}}
{\bf 1}_{\{i_4=i_5\ne 0\}}
{\bf 1}_{\{j_4=j_5\}}\zeta_{j_2}^{(i_2)}+
$$
$$
+
{\bf 1}_{\{i_1=i_4\ne 0\}}
{\bf 1}_{\{j_1=j_4\}}
{\bf 1}_{\{i_2=i_3\ne 0\}}
{\bf 1}_{\{j_2=j_3\}}\zeta_{j_5}^{(i_5)}+
{\bf 1}_{\{i_1=i_4\ne 0\}}
{\bf 1}_{\{j_1=j_4\}}
{\bf 1}_{\{i_2=i_5\ne 0\}}
{\bf 1}_{\{j_2=j_5\}}\zeta_{j_3}^{(i_3)}+
$$
$$
+
{\bf 1}_{\{i_1=i_4\ne 0\}}
{\bf 1}_{\{j_1=j_4\}}
{\bf 1}_{\{i_3=i_5\ne 0\}}
{\bf 1}_{\{j_3=j_5\}}\zeta_{j_2}^{(i_2)}+
{\bf 1}_{\{i_1=i_5\ne 0\}}
{\bf 1}_{\{j_1=j_5\}}
{\bf 1}_{\{i_2=i_3\ne 0\}}
{\bf 1}_{\{j_2=j_3\}}\zeta_{j_4}^{(i_4)}+
$$
$$
+
{\bf 1}_{\{i_1=i_5\ne 0\}}
{\bf 1}_{\{j_1=j_5\}}
{\bf 1}_{\{i_2=i_4\ne 0\}}
{\bf 1}_{\{j_2=j_4\}}\zeta_{j_3}^{(i_3)}+
{\bf 1}_{\{i_1=i_5\ne 0\}}
{\bf 1}_{\{j_1=j_5\}}
{\bf 1}_{\{i_3=i_4\ne 0\}}
{\bf 1}_{\{j_3=j_4\}}\zeta_{j_2}^{(i_2)}+
$$
$$
+
{\bf 1}_{\{i_2=i_3\ne 0\}}
{\bf 1}_{\{j_2=j_3\}}
{\bf 1}_{\{i_4=i_5\ne 0\}}
{\bf 1}_{\{j_4=j_5\}}\zeta_{j_1}^{(i_1)}+
{\bf 1}_{\{i_2=i_4\ne 0\}}
{\bf 1}_{\{j_2=j_4\}}
{\bf 1}_{\{i_3=i_5\ne 0\}}
{\bf 1}_{\{j_3=j_5\}}\zeta_{j_1}^{(i_1)}+
$$
\begin{equation}
\label{a5}
+\Biggl.
{\bf 1}_{\{i_2=i_5\ne 0\}}
{\bf 1}_{\{j_2=j_5\}}
{\bf 1}_{\{i_3=i_4\ne 0\}}
{\bf 1}_{\{j_3=j_4\}}\zeta_{j_1}^{(i_1)}\Biggr),
\end{equation}

\vspace{10mm}

$$
J[\psi^{(6)}]_{T,t}
=\hbox{\vtop{\offinterlineskip\halign{
\hfil#\hfil\cr
{\rm l.i.m.}\cr
$\stackrel{}{{}_{p_1,\ldots,p_6\to \infty}}$\cr
}} }\sum_{j_1=0}^{p_1}\ldots\sum_{j_6=0}^{p_6}
C_{j_6\ldots j_1}\Biggl(
\prod_{l=1}^6
\zeta_{j_l}^{(i_l)}
-\Biggr.
$$
$$
-
{\bf 1}_{\{i_1=i_6\ne 0\}}
{\bf 1}_{\{j_1=j_6\}}
\zeta_{j_2}^{(i_2)}
\zeta_{j_3}^{(i_3)}
\zeta_{j_4}^{(i_4)}
\zeta_{j_5}^{(i_5)}-
{\bf 1}_{\{i_2=i_6\ne 0\}}
{\bf 1}_{\{j_2=j_6\}}
\zeta_{j_1}^{(i_1)}
\zeta_{j_3}^{(i_3)}
\zeta_{j_4}^{(i_4)}
\zeta_{j_5}^{(i_5)}-
$$
$$
-
{\bf 1}_{\{i_3=i_6\ne 0\}}
{\bf 1}_{\{j_3=j_6\}}
\zeta_{j_1}^{(i_1)}
\zeta_{j_2}^{(i_2)}
\zeta_{j_4}^{(i_4)}
\zeta_{j_5}^{(i_5)}-
{\bf 1}_{\{i_4=i_6\ne 0\}}
{\bf 1}_{\{j_4=j_6\}}
\zeta_{j_1}^{(i_1)}
\zeta_{j_2}^{(i_2)}
\zeta_{j_3}^{(i_3)}
\zeta_{j_5}^{(i_5)}-
$$
$$
-
{\bf 1}_{\{i_5=i_6\ne 0\}}
{\bf 1}_{\{j_5=j_6\}}
\zeta_{j_1}^{(i_1)}
\zeta_{j_2}^{(i_2)}
\zeta_{j_3}^{(i_3)}
\zeta_{j_4}^{(i_4)}-
{\bf 1}_{\{i_1=i_2\ne 0\}}
{\bf 1}_{\{j_1=j_2\}}
\zeta_{j_3}^{(i_3)}
\zeta_{j_4}^{(i_4)}
\zeta_{j_5}^{(i_5)}
\zeta_{j_6}^{(i_6)}-
$$
$$
-
{\bf 1}_{\{i_1=i_3\ne 0\}}
{\bf 1}_{\{j_1=j_3\}}
\zeta_{j_2}^{(i_2)}
\zeta_{j_4}^{(i_4)}
\zeta_{j_5}^{(i_5)}
\zeta_{j_6}^{(i_6)}-
{\bf 1}_{\{i_1=i_4\ne 0\}}
{\bf 1}_{\{j_1=j_4\}}
\zeta_{j_2}^{(i_2)}
\zeta_{j_3}^{(i_3)}
\zeta_{j_5}^{(i_5)}
\zeta_{j_6}^{(i_6)}-
$$
$$
-
{\bf 1}_{\{i_1=i_5\ne 0\}}
{\bf 1}_{\{j_1=j_5\}}
\zeta_{j_2}^{(i_2)}
\zeta_{j_3}^{(i_3)}
\zeta_{j_4}^{(i_4)}
\zeta_{j_6}^{(i_6)}-
{\bf 1}_{\{i_2=i_3\ne 0\}}
{\bf 1}_{\{j_2=j_3\}}
\zeta_{j_1}^{(i_1)}
\zeta_{j_4}^{(i_4)}
\zeta_{j_5}^{(i_5)}
\zeta_{j_6}^{(i_6)}-
$$
$$
-
{\bf 1}_{\{i_2=i_4\ne 0\}}
{\bf 1}_{\{j_2=j_4\}}
\zeta_{j_1}^{(i_1)}
\zeta_{j_3}^{(i_3)}
\zeta_{j_5}^{(i_5)}
\zeta_{j_6}^{(i_6)}-
{\bf 1}_{\{i_2=i_5\ne 0\}}
{\bf 1}_{\{j_2=j_5\}}
\zeta_{j_1}^{(i_1)}
\zeta_{j_3}^{(i_3)}
\zeta_{j_4}^{(i_4)}
\zeta_{j_6}^{(i_6)}-
$$
$$
-
{\bf 1}_{\{i_3=i_4\ne 0\}}
{\bf 1}_{\{j_3=j_4\}}
\zeta_{j_1}^{(i_1)}
\zeta_{j_2}^{(i_2)}
\zeta_{j_5}^{(i_5)}
\zeta_{j_6}^{(i_6)}-
{\bf 1}_{\{i_3=i_5\ne 0\}}
{\bf 1}_{\{j_3=j_5\}}
\zeta_{j_1}^{(i_1)}
\zeta_{j_2}^{(i_2)}
\zeta_{j_4}^{(i_4)}
\zeta_{j_6}^{(i_6)}-
$$
$$
-
{\bf 1}_{\{i_4=i_5\ne 0\}}
{\bf 1}_{\{j_4=j_5\}}
\zeta_{j_1}^{(i_1)}
\zeta_{j_2}^{(i_2)}
\zeta_{j_3}^{(i_3)}
\zeta_{j_6}^{(i_6)}+
$$
$$
+
{\bf 1}_{\{i_1=i_2\ne 0\}}
{\bf 1}_{\{j_1=j_2\}}
{\bf 1}_{\{i_3=i_4\ne 0\}}
{\bf 1}_{\{j_3=j_4\}}
\zeta_{j_5}^{(i_5)}
\zeta_{j_6}^{(i_6)}+
{\bf 1}_{\{i_1=i_2\ne 0\}}
{\bf 1}_{\{j_1=j_2\}}
{\bf 1}_{\{i_3=i_5\ne 0\}}
{\bf 1}_{\{j_3=j_5\}}
\zeta_{j_4}^{(i_4)}
\zeta_{j_6}^{(i_6)}+
$$
$$
+
{\bf 1}_{\{i_1=i_2\ne 0\}}
{\bf 1}_{\{j_1=j_2\}}
{\bf 1}_{\{i_4=i_5\ne 0\}}
{\bf 1}_{\{j_4=j_5\}}
\zeta_{j_3}^{(i_3)}
\zeta_{j_6}^{(i_6)}
+
{\bf 1}_{\{i_1=i_3\ne 0\}}
{\bf 1}_{\{j_1=j_3\}}
{\bf 1}_{\{i_2=i_4\ne 0\}}
{\bf 1}_{\{j_2=j_4\}}
\zeta_{j_5}^{(i_5)}
\zeta_{j_6}^{(i_6)}+
$$
$$
+
{\bf 1}_{\{i_1=i_3\ne 0\}}
{\bf 1}_{\{j_1=j_3\}}
{\bf 1}_{\{i_2=i_5\ne 0\}}
{\bf 1}_{\{j_2=j_5\}}
\zeta_{j_4}^{(i_4)}
\zeta_{j_6}^{(i_6)}
+{\bf 1}_{\{i_1=i_3\ne 0\}}
{\bf 1}_{\{j_1=j_3\}}
{\bf 1}_{\{i_4=i_5\ne 0\}}
{\bf 1}_{\{j_4=j_5\}}
\zeta_{j_2}^{(i_2)}
\zeta_{j_6}^{(i_6)}+
$$
$$
+
{\bf 1}_{\{i_1=i_4\ne 0\}}
{\bf 1}_{\{j_1=j_4\}}
{\bf 1}_{\{i_2=i_3\ne 0\}}
{\bf 1}_{\{j_2=j_3\}}
\zeta_{j_5}^{(i_5)}
\zeta_{j_6}^{(i_6)}
+
{\bf 1}_{\{i_1=i_4\ne 0\}}
{\bf 1}_{\{j_1=j_4\}}
{\bf 1}_{\{i_2=i_5\ne 0\}}
{\bf 1}_{\{j_2=j_5\}}
\zeta_{j_3}^{(i_3)}
\zeta_{j_6}^{(i_6)}+
$$
$$
+
{\bf 1}_{\{i_1=i_4\ne 0\}}
{\bf 1}_{\{j_1=j_4\}}
{\bf 1}_{\{i_3=i_5\ne 0\}}
{\bf 1}_{\{j_3=j_5\}}
\zeta_{j_2}^{(i_2)}
\zeta_{j_6}^{(i_6)}
+
{\bf 1}_{\{i_1=i_5\ne 0\}}
{\bf 1}_{\{j_1=j_5\}}
{\bf 1}_{\{i_2=i_3\ne 0\}}
{\bf 1}_{\{j_2=j_3\}}
\zeta_{j_4}^{(i_4)}
\zeta_{j_6}^{(i_6)}+
$$
$$
+
{\bf 1}_{\{i_1=i_5\ne 0\}}
{\bf 1}_{\{j_1=j_5\}}
{\bf 1}_{\{i_2=i_4\ne 0\}}
{\bf 1}_{\{j_2=j_4\}}
\zeta_{j_3}^{(i_3)}
\zeta_{j_6}^{(i_6)}
+
{\bf 1}_{\{i_1=i_5\ne 0\}}
{\bf 1}_{\{j_1=j_5\}}
{\bf 1}_{\{i_3=i_4\ne 0\}}
{\bf 1}_{\{j_3=j_4\}}
\zeta_{j_2}^{(i_2)}
\zeta_{j_6}^{(i_6)}+
$$
$$
+
{\bf 1}_{\{i_2=i_3\ne 0\}}
{\bf 1}_{\{j_2=j_3\}}
{\bf 1}_{\{i_4=i_5\ne 0\}}
{\bf 1}_{\{j_4=j_5\}}
\zeta_{j_1}^{(i_1)}
\zeta_{j_6}^{(i_6)}
+
{\bf 1}_{\{i_2=i_4\ne 0\}}
{\bf 1}_{\{j_2=j_4\}}
{\bf 1}_{\{i_3=i_5\ne 0\}}
{\bf 1}_{\{j_3=j_5\}}
\zeta_{j_1}^{(i_1)}
\zeta_{j_6}^{(i_6)}+
$$
$$
+
{\bf 1}_{\{i_2=i_5\ne 0\}}
{\bf 1}_{\{j_2=j_5\}}
{\bf 1}_{\{i_3=i_4\ne 0\}}
{\bf 1}_{\{j_3=j_4\}}
\zeta_{j_1}^{(i_1)}
\zeta_{j_6}^{(i_6)}
+
{\bf 1}_{\{i_6=i_1\ne 0\}}
{\bf 1}_{\{j_6=j_1\}}
{\bf 1}_{\{i_3=i_4\ne 0\}}
{\bf 1}_{\{j_3=j_4\}}
\zeta_{j_2}^{(i_2)}
\zeta_{j_5}^{(i_5)}+
$$
$$
+
{\bf 1}_{\{i_6=i_1\ne 0\}}
{\bf 1}_{\{j_6=j_1\}}
{\bf 1}_{\{i_3=i_5\ne 0\}}
{\bf 1}_{\{j_3=j_5\}}
\zeta_{j_2}^{(i_2)}
\zeta_{j_4}^{(i_4)}
+
{\bf 1}_{\{i_6=i_1\ne 0\}}
{\bf 1}_{\{j_6=j_1\}}
{\bf 1}_{\{i_2=i_5\ne 0\}}
{\bf 1}_{\{j_2=j_5\}}
\zeta_{j_3}^{(i_3)}
\zeta_{j_4}^{(i_4)}+
$$
$$
+
{\bf 1}_{\{i_6=i_1\ne 0\}}
{\bf 1}_{\{j_6=j_1\}}
{\bf 1}_{\{i_2=i_4\ne 0\}}
{\bf 1}_{\{j_2=j_4\}}
\zeta_{j_3}^{(i_3)}
\zeta_{j_5}^{(i_5)}
+
{\bf 1}_{\{i_6=i_1\ne 0\}}
{\bf 1}_{\{j_6=j_1\}}
{\bf 1}_{\{i_4=i_5\ne 0\}}
{\bf 1}_{\{j_4=j_5\}}
\zeta_{j_2}^{(i_2)}
\zeta_{j_3}^{(i_3)}+
$$
$$
+
{\bf 1}_{\{i_6=i_1\ne 0\}}
{\bf 1}_{\{j_6=j_1\}}
{\bf 1}_{\{i_2=i_3\ne 0\}}
{\bf 1}_{\{j_2=j_3\}}
\zeta_{j_4}^{(i_4)}
\zeta_{j_5}^{(i_5)}
+
{\bf 1}_{\{i_6=i_2\ne 0\}}
{\bf 1}_{\{j_6=j_2\}}
{\bf 1}_{\{i_3=i_5\ne 0\}}
{\bf 1}_{\{j_3=j_5\}}
\zeta_{j_1}^{(i_1)}
\zeta_{j_4}^{(i_4)}+
$$
$$
+
{\bf 1}_{\{i_6=i_2\ne 0\}}
{\bf 1}_{\{j_6=j_2\}}
{\bf 1}_{\{i_4=i_5\ne 0\}}
{\bf 1}_{\{j_4=j_5\}}
\zeta_{j_1}^{(i_1)}
\zeta_{j_3}^{(i_3)}
+
{\bf 1}_{\{i_6=i_2\ne 0\}}
{\bf 1}_{\{j_6=j_2\}}
{\bf 1}_{\{i_3=i_4\ne 0\}}
{\bf 1}_{\{j_3=j_4\}}
\zeta_{j_1}^{(i_1)}
\zeta_{j_5}^{(i_5)}+
$$
$$
+
{\bf 1}_{\{i_6=i_2\ne 0\}}
{\bf 1}_{\{j_6=j_2\}}
{\bf 1}_{\{i_1=i_5\ne 0\}}
{\bf 1}_{\{j_1=j_5\}}
\zeta_{j_3}^{(i_3)}
\zeta_{j_4}^{(i_4)}
+
{\bf 1}_{\{i_6=i_2\ne 0\}}
{\bf 1}_{\{j_6=j_2\}}
{\bf 1}_{\{i_1=i_4\ne 0\}}
{\bf 1}_{\{j_1=j_4\}}
\zeta_{j_3}^{(i_3)}
\zeta_{j_5}^{(i_5)}+
$$
$$
+
{\bf 1}_{\{i_6=i_2\ne 0\}}
{\bf 1}_{\{j_6=j_2\}}
{\bf 1}_{\{i_1=i_3\ne 0\}}
{\bf 1}_{\{j_1=j_3\}}
\zeta_{j_4}^{(i_4)}
\zeta_{j_5}^{(i_5)}
+
{\bf 1}_{\{i_6=i_3\ne 0\}}
{\bf 1}_{\{j_6=j_3\}}
{\bf 1}_{\{i_2=i_5\ne 0\}}
{\bf 1}_{\{j_2=j_5\}}
\zeta_{j_1}^{(i_1)}
\zeta_{j_4}^{(i_4)}+
$$
$$
+
{\bf 1}_{\{i_6=i_3\ne 0\}}
{\bf 1}_{\{j_6=j_3\}}
{\bf 1}_{\{i_4=i_5\ne 0\}}
{\bf 1}_{\{j_4=j_5\}}
\zeta_{j_1}^{(i_1)}
\zeta_{j_2}^{(i_2)}
+
{\bf 1}_{\{i_6=i_3\ne 0\}}
{\bf 1}_{\{j_6=j_3\}}
{\bf 1}_{\{i_2=i_4\ne 0\}}
{\bf 1}_{\{j_2=j_4\}}
\zeta_{j_1}^{(i_1)}
\zeta_{j_5}^{(i_5)}+
$$
$$
+
{\bf 1}_{\{i_6=i_3\ne 0\}}
{\bf 1}_{\{j_6=j_3\}}
{\bf 1}_{\{i_1=i_5\ne 0\}}
{\bf 1}_{\{j_1=j_5\}}
\zeta_{j_2}^{(i_2)}
\zeta_{j_4}^{(i_4)}
+
{\bf 1}_{\{i_6=i_3\ne 0\}}
{\bf 1}_{\{j_6=j_3\}}
{\bf 1}_{\{i_1=i_4\ne 0\}}
{\bf 1}_{\{j_1=j_4\}}
\zeta_{j_2}^{(i_2)}
\zeta_{j_5}^{(i_5)}+
$$
$$
+
{\bf 1}_{\{i_6=i_3\ne 0\}}
{\bf 1}_{\{j_6=j_3\}}
{\bf 1}_{\{i_1=i_2\ne 0\}}
{\bf 1}_{\{j_1=j_2\}}
\zeta_{j_4}^{(i_4)}
\zeta_{j_5}^{(i_5)}
+
{\bf 1}_{\{i_6=i_4\ne 0\}}
{\bf 1}_{\{j_6=j_4\}}
{\bf 1}_{\{i_3=i_5\ne 0\}}
{\bf 1}_{\{j_3=j_5\}}
\zeta_{j_1}^{(i_1)}
\zeta_{j_2}^{(i_2)}+
$$
$$
+
{\bf 1}_{\{i_6=i_4\ne 0\}}
{\bf 1}_{\{j_6=j_4\}}
{\bf 1}_{\{i_2=i_5\ne 0\}}
{\bf 1}_{\{j_2=j_5\}}
\zeta_{j_1}^{(i_1)}
\zeta_{j_3}^{(i_3)}
+
{\bf 1}_{\{i_6=i_4\ne 0\}}
{\bf 1}_{\{j_6=j_4\}}
{\bf 1}_{\{i_2=i_3\ne 0\}}
{\bf 1}_{\{j_2=j_3\}}
\zeta_{j_1}^{(i_1)}
\zeta_{j_5}^{(i_5)}+
$$
$$
+
{\bf 1}_{\{i_6=i_4\ne 0\}}
{\bf 1}_{\{j_6=j_4\}}
{\bf 1}_{\{i_1=i_5\ne 0\}}
{\bf 1}_{\{j_1=j_5\}}
\zeta_{j_2}^{(i_2)}
\zeta_{j_3}^{(i_3)}
+
{\bf 1}_{\{i_6=i_4\ne 0\}}
{\bf 1}_{\{j_6=j_4\}}
{\bf 1}_{\{i_1=i_3\ne 0\}}
{\bf 1}_{\{j_1=j_3\}}
\zeta_{j_2}^{(i_2)}
\zeta_{j_5}^{(i_5)}+
$$
$$
+
{\bf 1}_{\{i_6=i_4\ne 0\}}
{\bf 1}_{\{j_6=j_4\}}
{\bf 1}_{\{i_1=i_2\ne 0\}}
{\bf 1}_{\{j_1=j_2\}}
\zeta_{j_3}^{(i_3)}
\zeta_{j_5}^{(i_5)}
+
{\bf 1}_{\{i_6=i_5\ne 0\}}
{\bf 1}_{\{j_6=j_5\}}
{\bf 1}_{\{i_3=i_4\ne 0\}}
{\bf 1}_{\{j_3=j_4\}}
\zeta_{j_1}^{(i_1)}
\zeta_{j_2}^{(i_2)}+
$$
$$
+
{\bf 1}_{\{i_6=i_5\ne 0\}}
{\bf 1}_{\{j_6=j_5\}}
{\bf 1}_{\{i_2=i_4\ne 0\}}
{\bf 1}_{\{j_2=j_4\}}
\zeta_{j_1}^{(i_1)}
\zeta_{j_3}^{(i_3)}
+
{\bf 1}_{\{i_6=i_5\ne 0\}}
{\bf 1}_{\{j_6=j_5\}}
{\bf 1}_{\{i_2=i_3\ne 0\}}
{\bf 1}_{\{j_2=j_3\}}
\zeta_{j_1}^{(i_1)}
\zeta_{j_4}^{(i_4)}+
$$
$$
+
{\bf 1}_{\{i_6=i_5\ne 0\}}
{\bf 1}_{\{j_6=j_5\}}
{\bf 1}_{\{i_1=i_4\ne 0\}}
{\bf 1}_{\{j_1=j_4\}}
\zeta_{j_2}^{(i_2)}
\zeta_{j_3}^{(i_3)}
+
{\bf 1}_{\{i_6=i_5\ne 0\}}
{\bf 1}_{\{j_6=j_5\}}
{\bf 1}_{\{i_1=i_3\ne 0\}}
{\bf 1}_{\{j_1=j_3\}}
\zeta_{j_2}^{(i_2)}
\zeta_{j_4}^{(i_4)}+
$$
$$
+
{\bf 1}_{\{i_6=i_5\ne 0\}}
{\bf 1}_{\{j_6=j_5\}}
{\bf 1}_{\{i_1=i_2\ne 0\}}
{\bf 1}_{\{j_1=j_2\}}
\zeta_{j_3}^{(i_3)}
\zeta_{j_4}^{(i_4)}-
$$
$$
-
{\bf 1}_{\{i_6=i_1\ne 0\}}
{\bf 1}_{\{j_6=j_1\}}
{\bf 1}_{\{i_2=i_5\ne 0\}}
{\bf 1}_{\{j_2=j_5\}}
{\bf 1}_{\{i_3=i_4\ne 0\}}
{\bf 1}_{\{j_3=j_4\}}-
$$
$$
-
{\bf 1}_{\{i_6=i_1\ne 0\}}
{\bf 1}_{\{j_6=j_1\}}
{\bf 1}_{\{i_2=i_4\ne 0\}}
{\bf 1}_{\{j_2=j_4\}}
{\bf 1}_{\{i_3=i_5\ne 0\}}
{\bf 1}_{\{j_3=j_5\}}-
$$
$$
-
{\bf 1}_{\{i_6=i_1\ne 0\}}
{\bf 1}_{\{j_6=j_1\}}
{\bf 1}_{\{i_2=i_3\ne 0\}}
{\bf 1}_{\{j_2=j_3\}}
{\bf 1}_{\{i_4=i_5\ne 0\}}
{\bf 1}_{\{j_4=j_5\}}-
$$
$$
-
{\bf 1}_{\{i_6=i_2\ne 0\}}
{\bf 1}_{\{j_6=j_2\}}
{\bf 1}_{\{i_1=i_5\ne 0\}}
{\bf 1}_{\{j_1=j_5\}}
{\bf 1}_{\{i_3=i_4\ne 0\}}
{\bf 1}_{\{j_3=j_4\}}-
$$
$$
-
{\bf 1}_{\{i_6=i_2\ne 0\}}
{\bf 1}_{\{j_6=j_2\}}
{\bf 1}_{\{i_1=i_4\ne 0\}}
{\bf 1}_{\{j_1=j_4\}}
{\bf 1}_{\{i_3=i_5\ne 0\}}
{\bf 1}_{\{j_3=j_5\}}-
$$
$$
-
{\bf 1}_{\{i_6=i_2\ne 0\}}
{\bf 1}_{\{j_6=j_2\}}
{\bf 1}_{\{i_1=i_3\ne 0\}}
{\bf 1}_{\{j_1=j_3\}}
{\bf 1}_{\{i_4=i_5\ne 0\}}
{\bf 1}_{\{j_4=j_5\}}-
$$
$$
-
{\bf 1}_{\{i_6=i_3\ne 0\}}
{\bf 1}_{\{j_6=j_3\}}
{\bf 1}_{\{i_1=i_5\ne 0\}}
{\bf 1}_{\{j_1=j_5\}}
{\bf 1}_{\{i_2=i_4\ne 0\}}
{\bf 1}_{\{j_2=j_4\}}-
$$
$$
-
{\bf 1}_{\{i_6=i_3\ne 0\}}
{\bf 1}_{\{j_6=j_3\}}
{\bf 1}_{\{i_1=i_4\ne 0\}}
{\bf 1}_{\{j_1=j_4\}}
{\bf 1}_{\{i_2=i_5\ne 0\}}
{\bf 1}_{\{j_2=j_5\}}-
$$
$$
-
{\bf 1}_{\{i_3=i_6\ne 0\}}
{\bf 1}_{\{j_3=j_6\}}
{\bf 1}_{\{i_1=i_2\ne 0\}}
{\bf 1}_{\{j_1=j_2\}}
{\bf 1}_{\{i_4=i_5\ne 0\}}
{\bf 1}_{\{j_4=j_5\}}-
$$
$$
-
{\bf 1}_{\{i_6=i_4\ne 0\}}
{\bf 1}_{\{j_6=j_4\}}
{\bf 1}_{\{i_1=i_5\ne 0\}}
{\bf 1}_{\{j_1=j_5\}}
{\bf 1}_{\{i_2=i_3\ne 0\}}
{\bf 1}_{\{j_2=j_3\}}-
$$
$$
-
{\bf 1}_{\{i_6=i_4\ne 0\}}
{\bf 1}_{\{j_6=j_4\}}
{\bf 1}_{\{i_1=i_3\ne 0\}}
{\bf 1}_{\{j_1=j_3\}}
{\bf 1}_{\{i_2=i_5\ne 0\}}
{\bf 1}_{\{j_2=j_5\}}-
$$
$$
-
{\bf 1}_{\{i_6=i_4\ne 0\}}
{\bf 1}_{\{j_6=j_4\}}
{\bf 1}_{\{i_1=i_2\ne 0\}}
{\bf 1}_{\{j_1=j_2\}}
{\bf 1}_{\{i_3=i_5\ne 0\}}
{\bf 1}_{\{j_3=j_5\}}-
$$
$$
-
{\bf 1}_{\{i_6=i_5\ne 0\}}
{\bf 1}_{\{j_6=j_5\}}
{\bf 1}_{\{i_1=i_4\ne 0\}}
{\bf 1}_{\{j_1=j_4\}}
{\bf 1}_{\{i_2=i_3\ne 0\}}
{\bf 1}_{\{j_2=j_3\}}-
$$
$$
-
{\bf 1}_{\{i_6=i_5\ne 0\}}
{\bf 1}_{\{j_6=j_5\}}
{\bf 1}_{\{i_1=i_2\ne 0\}}
{\bf 1}_{\{j_1=j_2\}}
{\bf 1}_{\{i_3=i_4\ne 0\}}
{\bf 1}_{\{j_3=j_4\}}-
$$
\begin{equation}
\label{a6}
\Biggl.-
{\bf 1}_{\{i_6=i_5\ne 0\}}
{\bf 1}_{\{j_6=j_5\}}
{\bf 1}_{\{i_1=i_3\ne 0\}}
{\bf 1}_{\{j_1=j_3\}}
{\bf 1}_{\{i_2=i_4\ne 0\}}
{\bf 1}_{\{j_2=j_4\}}\Biggr),
\end{equation}

\vspace{7mm}
\noindent
where ${\bf 1}_A$ is the indicator of the set $A$.

Thus, we obtain the following advantages and new possibilities
of the method of generalized multiple Fourier series (Theorem 1)
in comparison with the well known  methods of approximation 
of iterated stochastic integrals \cite{1988}-\cite{1994}, 
\cite{Zapad-xx-1}-\cite{Zapad-xx-11}.

1.\;There is an explicit formula (see (\ref{ppppa})) for calculation 
of expansion coefficients 
of the iterated Ito stochastic integral (\ref{ito}) with any
fixed multiplicity $k$ ($k\in\mathbb{N}$).

2.\;We have new possibilities for exact calculation 
and effective estimation of the mean-square 
approximation error
of iterated Ito stochastic integral \cite{2020}-\cite{2020xx1}, 
\cite{9a}, \cite{10},
\cite{15b} (see Sect.~3, 4).

3.\;Since the used
multiple Fourier series is a generalized in the sense
that it is built using various complete orthonormal
systems of functions in the space $L_2([t, T])$, then we 
have new possibilities 
for approximation --- we can 
use not only trigonometric functions as in \cite{1988}-\cite{1994},
\cite{Zapad-xx-1}, \cite{Zapad-xx-7}, \cite{Zapad-xx-8}, \cite{Zapad-xx-10}, 
\cite{Zapad-xx-11},
but Legendre polynomials.

4.\;As it turned out 
\cite{2020}-\cite{2020xx1}, \cite{13} (also see 
\cite{1998a}-\cite{arxiv-24}, \cite{1998b}-\cite{10},
\cite{15}, \cite{15d}-\cite{arxiv-12}, \cite{new-1}-\cite{new-4},
\cite{9999})-\cite{Mikh-1})
it is more convenient to work 
with the Legendre polynomials for constructing of approximations 
of the iterated Ito stochastic integrals (\ref{ito}). 
Approximations based on the Legendre polynomials essentially simpler 
than their analogues based on the trigonometric functions
(see \cite{1988}-\cite{1994},
\cite{Zapad-xx-1}-\cite{Zapad-xx-11}).
Another advantages of the application of Legendre polynomials 
in the framework of the mentioned problem are considered
in \cite{2020}-\cite{2020xx1} (Sect.~5.3), \cite{17}, \cite{arxiv-12}.

5.\;The approach to expansion of iterated 
stochastic integrals based on the Karhunen--Loeve expansion
of the Brownian bridge
process \cite{1988}-\cite{1994},
\cite{Zapad-xx-1}, \cite{Zapad-xx-7}, \cite{Zapad-xx-8}, \cite{Zapad-xx-10}, 
\cite{Zapad-xx-11} 
leads to 
iterated application of the operation of limit
transition (the operation of limit transition 
is implemented only once in Theorem 1)
starting from the 
second multiplicity (in the general case) 
and third multiplicity (for the case
$\psi_1(s), \psi_2(s), \psi_3(s)\equiv 1;$ 
$i_1, i_2, i_3=0,1,\ldots,m$)
of the iterated Ito stochastic integrals (\ref{ito}).
The same problem (iterated application of the operation of limit
transition) also appears in the method 
of expansion of iterated 
stochastic integrals based on the 
Wiener process series expansion using various complete 
orthonormal system of functions in the space $L_2([t, T])$
\cite{Zapad-xx-2}, \cite{Zapad-xx-3}.
Multiple series (the operation of limit transition 
is implemented only once) are more convenient 
for approximation than the iterated ones
(iterated application of the operation of limit
transition), 
since partial sums of multiple series converge for any possible case of  
convergence to infinity of their upper limits of summation 
(let us denote them as $p_1,\ldots, p_k$). 
For example, when
$p_1=\ldots=p_k=p\to\infty$. 
For iterated series, the condition $p_1=\ldots=p_k=p\to\infty$ obviously 
does not guarantee the convergence of this series.
However, 
in \cite{1992}
(Sect.~5.8, pp.~202--204), \cite{1994} (pp.~82-84),
\cite{Zapad-xx-1} (pp.~438-439),  
\cite{Zapad-xx-8} (pp.~263-264) 
the authors use (without rigorous proof)
the condition $p_1=p_2=p_3=p\to\infty$
within the frames of the mentioned approach
based on the Karhunen--Loeve expansion of the Brownian bridge
process \cite{1988} together with the Wong--Zakai approximation
\cite{W-Z-1}-\cite{Watanabe} (see discussions in
\cite{2020} (Sect.~2.18, 6.2), \cite{2020xx} (Sect.~2.6.2, 6.2),
\cite{2020xx1} (Sect.~2.6.2, 6.2),
\cite{11} (Sect.~11), \cite{12} (Sect.~6), 
\cite{13} (Sect.~8), \cite{14} (Sect.~6)
for detail).

6.\;Constructing the expansions of iterated
Ito stochastic integrals from Theorem 1, we 
saved all 
information about these integrals. That is why it is 
natural to expect that the mentioned expansions will converge
with probability 1 and in the mean of degree $2n$ ($n\in\mathbb{N}$).
The convergence with probability 1 in Theorem 1 
is proved \cite{2020}-\cite{2020xx1}, 
\cite{11}, \cite{13}, \cite{15b}, \cite{9999}
for complete orthonormal systems of Legendre polynomials 
and trigonometric functions
in the space $L_2([t,T])$.
Furthermore, the convergence 
in the mean of degree $2n$ ($n\in \mathbb{N}$) in Theorem 1 is proved in
\cite{2007b}-\cite{2020xx1},
\cite{5a}-\cite{7}, \cite{9a}-\cite{11}.

7.\;The versions of Theorem 1 for complete 
orthonormal with weight  
$r(t_1)\ldots r(t_k)$ systems of functions in 
the space $L_2([t,T]^k)$ ($k\in \mathbb{N}$) 
as well 
as for some other types of iterated stochastic 
integrals (iterated stochastic integrals 
with respect to martingale Poisson measures and 
iterated stochastic integrals with respect 
to martingales) were obtained in 
\cite{2020}-\cite{2020xx1}, \cite{xxx}
(also see \cite{2006}-\cite{2010}, 
\cite{5a}-\cite{7}, \cite{9a}, \cite{10}).

8.\;The adaptation of Theorem 1 for iterated Stratonovich
stochastic integrals of multiplicities 1 to 6 is realized in 
\cite{2020}-\cite{2020xx1}, \cite{5a}-\cite{7},
\cite{9a}, \cite{10}, \cite{12}, \cite{14}, \cite{15a},
\cite{15c} (see Theorems 3--7 below).

9.\;Application of Theorem 1 and Theorem 2 (see below) for the mean-square
approximation of iterated stochastic integrals 
with respect to the 
infinite-dimensional $Q$-Wiener process can be found
in the monographs \cite{2020}-\cite{2020xx1} (Chapter 7) and in 
\cite{new-1}-\cite{new-4}.

For further consideration, let us 
consider the generalization of formulas (\ref{a1})--(\ref{a6})                 
for the case of an arbitrary multiplicity $k$ $(k\in\mathbb{N})$ of 
the iterated Ito stochastic integral $J[\psi^{(k)}]_{T,t}$ defined by (\ref{ito})
as well as for the case of an arbitrary complete orthonormal 
systems of functions in $L_2([t, T])$ and $\psi_1(\tau),\ldots,\psi_k(\tau)\in L_2([t, T])$.
In order to do this, let us
introduce some notations. 
Consider the unordered
set $\{1, 2, \ldots, k\}$ 
and separate it into two parts:
the first part consists of $r$ unordered 
pairs (sequence order of these pairs is also unimportant) and the 
second one consists of the 
remaining $k-2r$ numbers.
So, we have

\begin{equation}
\label{leto5007}
(\{
\underbrace{\{g_1, g_2\}, \ldots, 
\{g_{2r-1}, g_{2r}\}}_{\small{\hbox{part 1}}}
\},
\{\underbrace{q_1, \ldots, q_{k-2r}}_{\small{\hbox{part 2}}}
\}),
\end{equation}

\vspace{4mm}
\noindent
where 

\vspace{-2mm}
$$
\{g_1, g_2, \ldots, 
g_{2r-1}, g_{2r}, q_1, \ldots, q_{k-2r}\}=\{1, 2, \ldots, k\},
$$

\vspace{4mm}
\noindent
braces   
mean an unordered 
set, and pa\-ren\-the\-ses mean an ordered set.

We will say that (\ref{leto5007}) is a partition 
and consider the sum with respect to all possible
partitions

\begin{equation}
\label{leto5008}
\sum_{\stackrel{(\{\{g_1, g_2\}, \ldots, 
\{g_{2r-1}, g_{2r}\}\}, \{q_1, \ldots, q_{k-2r}\})}
{{}_{\{g_1, g_2, \ldots, 
g_{2r-1}, g_{2r}, q_1, \ldots, q_{k-2r}\}=\{1, 2, \ldots, k\}}}}
a_{g_1 g_2, \ldots, 
g_{2r-1} g_{2r}, q_1 \ldots q_{k-2r}}.
\end{equation}

\vspace{4mm}

Below there are several examples of sums in the form (\ref{leto5008})

\vspace{2mm}
$$
\sum_{\stackrel{(\{g_1, g_2\})}{{}_{\{g_1, g_2\}=\{1, 2\}}}}
a_{g_1 g_2}=a_{12},
$$

\vspace{3mm}
$$
\sum_{\stackrel{(\{\{g_1, g_2\}, \{g_3, g_4\}\})}
{{}_{\{g_1, g_2, g_3, g_4\}=\{1, 2, 3, 4\}}}}
a_{g_1 g_2 g_3 g_4}=a_{1234} + a_{1324} + a_{2314},
$$

\vspace{3mm}
$$
\sum_{\stackrel{(\{g_1, g_2\}, \{q_1, q_{2}\})}
{{}_{\{g_1, g_2, q_1, q_{2}\}=\{1, 2, 3, 4\}}}}
a_{g_1 g_2, q_1 q_{2}}=
$$

$$
=a_{12,34}+a_{13,24}+a_{14,23}
+a_{23,14}+a_{24,13}+a_{34,12},
$$

\vspace{3mm}
$$
\sum_{\stackrel{(\{g_1, g_2\}, \{q_1, q_{2}, q_3\})}
{{}_{\{g_1, g_2, q_1, q_{2}, q_3\}=\{1, 2, 3, 4, 5\}}}}
a_{g_1 g_2, q_1 q_{2}q_3}
=
$$

$$
=a_{12,345}+a_{13,245}+a_{14,235}
+a_{15,234}+a_{23,145}+a_{24,135}+
$$
$$
+a_{25,134}+a_{34,125}+a_{35,124}+a_{45,123},
$$

\vspace{4mm}
$$
\sum_{\stackrel{(\{\{g_1, g_2\}, \{g_3, g_{4}\}\}, \{q_1\})}
{{}_{\{g_1, g_2, g_3, g_{4}, q_1\}=\{1, 2, 3, 4, 5\}}}}
a_{g_1 g_2, g_3 g_{4},q_1}
=
$$

$$
=
a_{12,34,5}+a_{13,24,5}+a_{14,23,5}+
a_{12,35,4}+a_{13,25,4}+a_{15,23,4}+
$$
$$
+a_{12,54,3}+a_{15,24,3}+a_{14,25,3}+a_{15,34,2}+a_{13,54,2}+a_{14,53,2}+
$$
$$
+
a_{52,34,1}+a_{53,24,1}+a_{54,23,1}.
$$

\vspace{5mm}

Now we can formulate the following generalization of Theorem~1.

\vspace{2mm}

{\bf Theorem~2}\ \cite{2020} (Sect.~1.11), \cite{11} (Sect.~15). 
{\it Suppose that
$\psi_1(\tau),\ldots,\psi_k(\tau)\in L_2([t, T])$ and
$\{\phi_j(x)\}_{j=0}^{\infty}$ is an arbitrary complete orthonormal system  
of functions in the space $L_2([t,T]).$
Then the following expansion

\vspace{1mm}
$$
J[\psi^{(k)}]_{T,t}=
\hbox{\vtop{\offinterlineskip\halign{
\hfil#\hfil\cr
{\rm l.i.m.}\cr
$\stackrel{}{{}_{p_1,\ldots,p_k\to \infty}}$\cr
}} }
\sum\limits_{j_1=0}^{p_1}\ldots
\sum\limits_{j_k=0}^{p_k}
C_{j_k\ldots j_1}\Biggl(
\prod_{l=1}^k\zeta_{j_l}^{(i_l)}+\sum\limits_{r=1}^{[k/2]}
(-1)^r \times
\Biggr.
$$

\vspace{1mm}
\begin{equation}
\label{tyyy111}
\times
\sum_{\stackrel{(\{\{g_1, g_2\}, \ldots, 
\{g_{2r-1}, g_{2r}\}\}, \{q_1, \ldots, q_{k-2r}\})}
{{}_{\{g_1, g_2, \ldots, 
g_{2r-1}, g_{2r}, q_1, \ldots, q_{k-2r}\}=\{1, 2, \ldots, k\}}}}
\prod\limits_{s=1}^r
{\bf 1}_{\{i_{g_{{}_{2s-1}}}=~i_{g_{{}_{2s}}}\ne 0\}}
\Biggl.{\bf 1}_{\{j_{g_{{}_{2s-1}}}=~j_{g_{{}_{2s}}}\}}
\prod_{l=1}^{k-2r}\zeta_{j_{q_l}}^{(i_{q_l})}\Biggr)
\end{equation}

\vspace{5mm}
\noindent
con\-verg\-ing in the mean-square sense is valid,
where $[x]$ is an integer part of a real number $x;$
another notations are the same as in Theorem~{\rm 1}.}

\vspace{2mm}

It should be noted that an analogue of Theorem 2 was considered 
in \cite{Rybakov1000}. 
Note that we use another notations 
\cite{2020} (Sect.~1.11), \cite{11} (Sect.~15)
in comparison with \cite{Rybakov1000}.
Moreover, the proof of an analogue of Theorem 2
from \cite{Rybakov1000} is somewhat different from the proof given in 
\cite{2020} (Sect.~1.11), \cite{11} (Sect.~15).

As we mentioned above, 
in a number of works 
\cite{2020}-\cite{2020xx1}, \cite{5a}-\cite{7},
\cite{9a}, \cite{10}, \cite{12}, \cite{14}, \cite{15a},
\cite{15c} 
Theorem 1 is adapted for the iterated  Stratonovich 
stochastic integrals
(\ref{str}) of multiplicities 2 to 6 (the case of multiplicity 1 is given by (\ref{a1})). 
Let as collect some old results in the following theorem.

\vspace{2mm}

{\bf Theorem 3} \cite{2020}-\cite{2020xx1}, \cite{5a}-\cite{7},
\cite{9a}, \cite{10}, \cite{12}, \cite{14}, \cite{15a},
\cite{15c}. {\it Suppose that 
$\{\phi_j(x)\}_{j=0}^{\infty}$ is a complete orthonormal system of 
Legendre polynomials or trigonometric functions in the space $L_2([t, T]).$
At the same time $\psi_2(\tau)$ is a continuously differentiable 
function on $[t, T]$ and $\psi_1(\tau),$ $\psi_3(\tau)$ are twice
continuously differentiable functions on $[t, T]$. Then 

\vspace{-1mm}
\begin{equation}
\label{a}
J^{*}[\psi^{(2)}]_{T,t}=
\hbox{\vtop{\offinterlineskip\halign{
\hfil#\hfil\cr
{\rm l.i.m.}\cr
$\stackrel{}{{}_{p_1,p_2\to \infty}}$\cr
}} }\sum_{j_1=0}^{p_1}\sum_{j_2=0}^{p_2}
C_{j_2j_1}\zeta_{j_1}^{(i_1)}\zeta_{j_2}^{(i_2)}\ \ (i_1,i_2=1,\ldots,m),
\end{equation}

\begin{equation}
\label{feto19000ab}
J^{*}[\psi^{(3)}]_{T,t}=
\hbox{\vtop{\offinterlineskip\halign{
\hfil#\hfil\cr
{\rm l.i.m.}\cr
$\stackrel{}{{}_{p_1,p_2,p_3\to \infty}}$\cr
}} }\sum_{j_1=0}^{p_1}\sum_{j_2=0}^{p_2}\sum_{j_3=0}^{p_3}
C_{j_3 j_2 j_1}\zeta_{j_1}^{(i_1)}\zeta_{j_2}^{(i_2)}\zeta_{j_3}^{(i_3)}\ \
(i_1,i_2,i_3=0, 1,\ldots,m),
\end{equation}

\begin{equation}
\label{feto19000a}
J^{*}[\psi^{(3)}]_{T,t}=
\hbox{\vtop{\offinterlineskip\halign{
\hfil#\hfil\cr
{\rm l.i.m.}\cr
$\stackrel{}{{}_{p\to \infty}}$\cr
}} }
\sum\limits_{j_1,j_2,j_3=0}^{p}
C_{j_3 j_2 j_1}\zeta_{j_1}^{(i_1)}\zeta_{j_2}^{(i_2)}\zeta_{j_3}^{(i_3)}\ \
(i_1,i_2,i_3=1,\ldots,m),
\end{equation}

\begin{equation}
\label{uu}
J^{*}[\psi^{(4)}]_{T,t}=
\hbox{\vtop{\offinterlineskip\halign{
\hfil#\hfil\cr
{\rm l.i.m.}\cr
$\stackrel{}{{}_{p\to \infty}}$\cr
}} }
\sum\limits_{j_1, \ldots, j_4=0}^{p}
C_{j_4 j_3 j_2 j_1}\zeta_{j_1}^{(i_1)}
\zeta_{j_2}^{(i_2)}\zeta_{j_3}^{(i_3)}\zeta_{j_4}^{(i_4)}\ \
(i_1,\ldots,i_4=0, 1,\ldots,m),
\end{equation}

\vspace{5mm}
\noindent
where $J^{*}[\psi^{(k)}]_{T,t}$ is defined by {\rm (\ref{str})}, and
$\psi_l(\tau)\equiv 1$ $(l=1,\ldots,4)$ in {\rm (\ref{feto19000ab})}, 
{\rm (\ref{uu});} another notations are the same as in Theorems {\rm 1, 2.}
}

Recently, a new approach to the expansion and mean-square 
approximation of iterated Stratonovich stochastic integrals has been obtained
\cite{2020} (Sect.~2.10--2.16), \cite{12} (Sect.~13--19), 
\cite{15a} (Sect.~5--11), \cite{arxiv-4} (Sect.~7--13), \cite{new-art-1xxy}
(Sect.~4--9).
Let us formulate four theorems that were obtained using this approach.

\vspace{2mm}

{\bf Theorem 4}\ \cite{2020}, \cite{12}, \cite{15a}, \cite{arxiv-4}, \cite{new-art-1xxy}.\
{\it Suppose 
that $\{\phi_j(x)\}_{j=0}^{\infty}$ is a complete orthonormal system of 
Legendre polynomials or trigonometric functions in the space $L_2([t, T]).$
Furthermore, let $\psi_1(\tau), \psi_2(\tau),$ $\psi_3(\tau)$ are continuously dif\-ferentiable 
nonrandom functions on $[t, T].$ 
Then, for the 
iterated Stra\-to\-no\-vich stochastic integral of third multiplicity

\vspace{-1mm}
$$
J^{*}[\psi^{(3)}]_{T,t}={\int\limits_t^{*}}^T\psi_3(t_3)
{\int\limits_t^{*}}^{t_3}\psi_2(t_2)
{\int\limits_t^{*}}^{t_2}\psi_1(t_1)
d{\bf w}_{t_1}^{(i_1)}
d{\bf w}_{t_2}^{(i_2)}d{\bf w}_{t_3}^{(i_3)}\ \ \ (i_1,i_2,i_3=0,1,\ldots,m)
$$

\vspace{3mm}
\noindent
the following 
relations

\vspace{-1mm}
\begin{equation}
\label{fin1}
J^{*}[\psi^{(3)}]_{T,t}
=\hbox{\vtop{\offinterlineskip\halign{
\hfil#\hfil\cr
{\rm l.i.m.}\cr
$\stackrel{}{{}_{p\to \infty}}$\cr
}} }
\sum\limits_{j_1, j_2, j_3=0}^{p}
C_{j_3 j_2 j_1}\zeta_{j_1}^{(i_1)}\zeta_{j_2}^{(i_2)}\zeta_{j_3}^{(i_3)},
\end{equation}

\vspace{1mm}
\begin{equation}
\label{fin2}
{\sf M}\left\{\left(
J^{*}[\psi^{(3)}]_{T,t}-
\sum\limits_{j_1, j_2, j_3=0}^{p}
C_{j_3 j_2 j_1}\zeta_{j_1}^{(i_1)}\zeta_{j_2}^{(i_2)}\zeta_{j_3}^{(i_3)}\right)^2\right\}
\le \frac{C}{p}
\end{equation}

\vspace{4mm}
\noindent
are fulfilled, where $i_1, i_2, i_3=0,1,\ldots,m$ in {\rm (\ref{fin1})} and 
$i_1, i_2, i_3=1,\ldots,m$ in {\rm (\ref{fin2})},
constant $C$ is independent of $p,$

\vspace{-1mm}
$$
C_{j_3 j_2 j_1}=\int\limits_t^T\psi_3(t_3)\phi_{j_3}(t_3)
\int\limits_t^{t_3}\psi_2(t_2)\phi_{j_2}(t_2)
\int\limits_t^{t_2}\psi_1(t_1)\phi_{j_1}(t_1)dt_1dt_2dt_3
$$

\vspace{4mm}
\noindent
and
$$
\zeta_{j}^{(i)}=
\int\limits_t^T \phi_{j}(\tau) d{\bf f}_{\tau}^{(i)}
$$ 

\vspace{2mm}
\noindent
are independent standard Gaussian random variables for various 
$i$ or $j$ {\rm (}in the case when $i\ne 0${\rm );} 
another notations are the same as in Theorems~{\rm 1, 2}.}

\vspace{2mm}

{\bf Theorem 5}\ \cite{2020}, \cite{12}, \cite{15a}, \cite{arxiv-4}, \cite{new-art-1xxy}.\ 
{\it Let
$\{\phi_j(x)\}_{j=0}^{\infty}$ be a complete orthonormal system of 
Legendre polynomials or trigonometric functions in the space $L_2([t, T]).$
Furthermore, let $\psi_1(\tau), \ldots,$ $\psi_4(\tau)$ be continuously dif\-ferentiable 
nonrandom functions on $[t, T].$ 
Then, for the 
iterated Stra\-to\-no\-vich stochastic integral of fourth multiplicity

\vspace{-1mm}
\begin{equation}
\label{fin0}
J^{*}[\psi^{(4)}]_{T,t}={\int\limits_t^{*}}^T\psi_4(t_4)
{\int\limits_t^{*}}^{t_4}\psi_3(t_3)
{\int\limits_t^{*}}^{t_3}\psi_2(t_2)
{\int\limits_t^{*}}^{t_2}\psi_1(t_1)
d{\bf w}_{t_1}^{(i_1)}
d{\bf w}_{t_2}^{(i_2)}d{\bf w}_{t_3}^{(i_3)}d{\bf w}_{t_4}^{(i_4)}
\end{equation}

\vspace{3mm}
\noindent
the following 
relations

\vspace{-1mm}
\begin{equation}
\label{fin3}
J^{*}[\psi^{(4)}]_{T,t}
=\hbox{\vtop{\offinterlineskip\halign{
\hfil#\hfil\cr
{\rm l.i.m.}\cr
$\stackrel{}{{}_{p\to \infty}}$\cr
}} }
\sum\limits_{j_1, j_2, j_3,j_4=0}^{p}
C_{j_4j_3 j_2 j_1}\zeta_{j_1}^{(i_1)}\zeta_{j_2}^{(i_2)}\zeta_{j_3}^{(i_3)}\zeta_{j_4}^{(i_4)},
\end{equation}

\vspace{1mm}

\begin{equation}
\label{fin4}
{\sf M}\left\{\left(
J^{*}[\psi^{(4)}]_{T,t}-
\sum\limits_{j_1, j_2, j_3, j_4=0}^{p}
C_{j_4 j_3 j_2 j_1}\zeta_{j_1}^{(i_1)}\zeta_{j_2}^{(i_2)}\zeta_{j_3}^{(i_3)}
\zeta_{j_4}^{(i_4)}
\right)^2\right\}
\le \frac{C}{p^{1-\varepsilon}}
\end{equation}

\vspace{4mm}
\noindent
are fulfilled, where $i_1, \ldots , i_4=0,1,\ldots,m$ in {\rm (\ref{fin0}),} {\rm (\ref{fin3})} 
and $i_1, \ldots, i_4=1,\ldots,m$ in {\rm (\ref{fin4}),}
constant $C$ does not depend on $p,$
$\varepsilon$ is an arbitrary
small positive real number 
for the case of complete orthonormal system of 
Legendre polynomials in the space $L_2([t, T])$
and $\varepsilon=0$ for the case of
complete orthonormal system of 
trigonometric functions in the space $L_2([t, T]),$

\vspace{-1mm}
$$
C_{j_4 j_3 j_2 j_1}=
\int\limits_t^T\psi_4(t_4)\phi_{j_4}(t_4)
\int\limits_t^{t_4}\psi_3(t_3)\phi_{j_3}(t_3)
\int\limits_t^{t_3}\psi_2(t_2)\phi_{j_2}(t_2)
\int\limits_t^{t_2}\psi_1(t_1)\phi_{j_1}(t_1)dt_1dt_2dt_3dt_4;
$$

\vspace{5mm}
\noindent
another notations are the same as in Theorem~{\rm 4}.}

\vspace{2mm}

{\bf Theorem 6}\ \cite{2020}, \cite{12}, \cite{15a}, \cite{arxiv-4}, \cite{new-art-1xxy}.\
{\it Assume 
that $\{\phi_j(x)\}_{j=0}^{\infty}$ is a complete orthonormal system of 
Legendre polynomials or trigonometric functions in the space $L_2([t, T])$
and $\psi_1(\tau), \ldots,$ $\psi_5(\tau)$ are continuously dif\-ferentiable 
nonrandom functions on $[t, T].$ 
Then, for the 
iterated Stra\-to\-no\-vich stochastic integral of fifth multiplicity

\vspace{-1mm}
\begin{equation}
\label{fin7}
J^{*}[\psi^{(5)}]_{T,t}={\int\limits_t^{*}}^T\psi_5(t_5)
\ldots
{\int\limits_t^{*}}^{t_2}\psi_1(t_1)
d{\bf w}_{t_1}^{(i_1)}
\ldots d{\bf w}_{t_5}^{(i_5)}
\end{equation}

\vspace{3mm}
\noindent
the following 
relations

\vspace{-1mm}
\begin{equation}
\label{fin8}
J^{*}[\psi^{(5)}]_{T,t}
=\hbox{\vtop{\offinterlineskip\halign{
\hfil#\hfil\cr
{\rm l.i.m.}\cr
$\stackrel{}{{}_{p\to \infty}}$\cr
}} }
\sum\limits_{j_1,\ldots,j_5=0}^{p}
C_{j_5 \ldots j_1}\zeta_{j_1}^{(i_1)}\ldots \zeta_{j_5}^{(i_5)},
\end{equation}

\vspace{1mm}

\begin{equation}
\label{fin9}
{\sf M}\left\{\left(
J^{*}[\psi^{(5)}]_{T,t}-
\sum\limits_{j_1, \ldots, j_5=0}^{p}
C_{j_5 \ldots j_1}\zeta_{j_1}^{(i_1)}\ldots
\zeta_{j_5}^{(i_5)}
\right)^2\right\}
\le \frac{C}{p^{1-\varepsilon}}
\end{equation}

\vspace{4mm}
\noindent
are fulfilled, where $i_1, \ldots , i_5=0,1,\ldots,m$ in {\rm (\ref{fin7}),} {\rm (\ref{fin8})} 
and $i_1, \ldots, i_5=1,\ldots,m$ in {\rm (\ref{fin9}),}
constant $C$ is independent of $p,$
$\varepsilon$ is an arbitrary
small positive real number 
for the case of complete orthonormal system of 
Legendre polynomials in the space $L_2([t, T])$
and $\varepsilon=0$ for the case of
complete orthonormal system of 
trigonometric functions in the space $L_2([t, T]),$

\vspace{-1mm}
$$
C_{j_5 \ldots j_1}=
\int\limits_t^T\psi_5(t_5)\phi_{j_5}(t_5)\ldots
\int\limits_t^{t_2}\psi_1(t_1)\phi_{j_1}(t_1)dt_1\ldots dt_5;
$$

\vspace{4mm}
\noindent
another notations are the same as in Theorems~{\rm 4, 5}.}

\vspace{2mm}

{\bf Theorem 7}\ \cite{2020}, \cite{12}, \cite{15a}, \cite{arxiv-4}.\
{\it Suppose that 
$\{\phi_j(x)\}_{j=0}^{\infty}$ is a complete orthonormal system of 
Legendre polynomials or trigonometric functions in the space $L_2([t, T]).$
Then, for the 
iterated Stratonovich stochastic integral of sixth multiplicity

$$
J_{T,t}^{*(i_1\ldots i_6)}={\int\limits_t^{*}}^T
\ldots
{\int\limits_t^{*}}^{t_2}
d{\bf w}_{t_1}^{(i_1)}
\ldots d{\bf w}_{t_6}^{(i_6)}
$$

\vspace{3mm}
\noindent
the following 
expansion 

\vspace{-1mm}
$$
J_{T,t}^{*(i_1\ldots i_6)}
=\hbox{\vtop{\offinterlineskip\halign{
\hfil#\hfil\cr
{\rm l.i.m.}\cr
$\stackrel{}{{}_{p\to \infty}}$\cr
}} }
\sum\limits_{j_1, \ldots, j_6=0}^{p}
C_{j_6 \ldots j_1}\zeta_{j_1}^{(i_1)}\ldots
\zeta_{j_6}^{(i_6)}
$$

\vspace{4mm}
\noindent
that converges in the mean-square sense is valid, where
$i_1, \ldots, i_6=0, 1,\ldots,m,$

$$
C_{j_6 \ldots j_1}=
\int\limits_t^T\phi_{j_6}(t_6)\ldots
\int\limits_t^{t_2}\phi_{j_1}(t_1)dt_1\ldots dt_6;
$$

\vspace{3mm}
\noindent
another notations are the same as in Theorems~{\rm 4--6}.}

\vspace{5mm}

\section{Estimate for the Mean-Square Approximation Error
in the Method of Approximation of Iterated Ito Stochastic 
Integrals Based on Generalized Multiple Fourier Series}

\vspace{5mm}

Assume that $J[\psi^{(k)}]_{T,t}^{p}$ is the approximation 
of (\ref{ito}), which is
the expression on the right-hand side of (\ref{tyyy111}) 
before passing to the limit 
$\hbox{\vtop{\offinterlineskip\halign{
\hfil#\hfil\cr
{\rm l.i.m.}\cr
$\stackrel{}{{}_{p_1,\ldots,p_k\to \infty}}$\cr
}} }
$ for the case $p_1=\ldots =p_k=p$, i.e.
                               
\vspace{1mm}

$$
J[\psi^{(k)}]_{T,t}^p=
\sum_{j_1,\ldots,j_k=0}^{p}
C_{j_k\ldots j_1}\Biggl(
\prod_{l=1}^k\zeta_{j_l}^{(i_l)}+\sum\limits_{r=1}^{[k/2]}
(-1)^r \times
\Biggr.
$$

\vspace{2mm}
$$
\times
\sum_{\stackrel{(\{\{g_1, g_2\}, \ldots, 
\{g_{2r-1}, g_{2r}\}\}, \{q_1, \ldots, q_{k-2r}\})}
{{}_{\{g_1, g_2, \ldots, 
g_{2r-1}, g_{2r}, q_1, \ldots, q_{k-2r}\}=\{1, 2, \ldots, k\}}}}
\prod\limits_{s=1}^r
{\bf 1}_{\{i_{g_{{}_{2s-1}}}=~i_{g_{{}_{2s}}}\ne 0\}}
\Biggl.{\bf 1}_{\{j_{g_{{}_{2s-1}}}=~j_{g_{{}_{2s}}}\}}
\prod_{l=1}^{k-2r}\zeta_{j_{q_l}}^{(i_{q_l})}\Biggr).
$$

\vspace{5mm}

Let us denote

\vspace{-1mm}
$$
{\sf M}\left\{\left(J[\psi^{(k)}]_{T,t}-
J[\psi^{(k)}]_{T,t}^{p}\right)^2\right\}\stackrel{{\rm def}}
{=}E_k^{p},
$$

\vspace{2mm}

$$
\left\Vert K\right\Vert^2_{L_2([t,T]^k)}=\int\limits_{[t,T]^k}
K^2(t_1,\ldots,t_k)dt_1\ldots dt_k\stackrel{{\rm def}}{=}I_k.
$$

\vspace{4mm}

When proving Theorems 1 and 2 
\cite{2020}-\cite{2020xx1} (also see \cite{7}-\cite{11}),
we have proved
the following estimate

\begin{equation}
\label{star00011}
E_k^p
\le k!\left(I_k
-\sum_{j_1,\ldots,j_k=0}^{p}C_{j_k\ldots j_1}^2\right),
\end{equation}

\vspace{4mm}
\noindent
where $i_1,\ldots,i_k=1,\ldots,m$ for $T-t\in (0,\infty)$ and 
$i_1,\ldots,i_k=0, 1,\ldots,m$ for $T-t\in (0, 1);$
another notations are the same as in Theorems 1, 2.

Combining the estimates (\ref{uslov}) and (\ref{star00011}), we obtain

\vspace{-1mm}
\begin{equation}
\label{sme}
k!\left(I_k
-\sum_{j_1,\ldots,j_k=0}^{p}C_{j_k\ldots j_1}^2\right)\le C(T-t)^{r+1},
\end{equation}

\vspace{3mm}
\noindent
where constant $C$ is independent of $T-t.$

It is not difficult to see that the multiplier factor $k!$ 
on the left-hand side of (\ref{sme}) leads to a significant
increase of the minimal natural number $p$
satisfying the estimate (\ref{sme}). For example,
for the numerical methods (\ref{al1})--(\ref{al4}) 
we will have the following multiplier
factors on the left-hand side of the inequality (\ref{sme}): $2!=2,$ 
$3!=6,$ $4!=24,$ $5!=120.$

As we will see in the next section, the mentioned 
problem can be partially overcome if we calculate 
the mean-square approximation error $E_k^{p}$ exactly.

\vspace{5mm}

\section{Exact Formulas for the Mean-Square Approximation Error
in the Method of Approximation of Iterated Ito Stochastic 
Integrals Based on Generalized Multiple Fourier Series}

\vspace{5mm}

This section is devoted
to 
exact 
expressions for the mean-square
approximation error in Theorems 1, 2 for iterated Ito stochastic integrals of
arbitrary multiplicity $k$ ($k\in \mathbb{N}$).

As it turned out, the value $E_k^{p}$
can be calculated exactly.

\vspace{2mm}

{\bf Theorem 8} \cite{2020} (Sect.~1.12), \cite{15b} (Sect.~6).
{\it Suppose that $\{\phi_j(x)\}_{j=0}^{\infty}$ 
is an arbitrary complete orthonormal system  
of functions in the space $L_2([t,T])$ and
$\psi_1(\tau),\ldots,\psi_k(\tau)\in L_2([t, T]),$  $i_1,\ldots, i_k=1,\ldots,m$.
Then

\vspace{-1mm}
\begin{equation}
\label{tttr11}
E_k^p=I_k- \sum_{j_1,\ldots, j_k=0}^{p}
C_{j_k\ldots j_1}
{\sf M}\left\{J[\psi^{(k)}]_{T,t}
\sum\limits_{(j_1,\ldots,j_k)}
\int\limits_t^T \phi_{j_k}(t_k)
\ldots
\int\limits_t^{t_{2}}\phi_{j_{1}}(t_{1})
d{\bf f}_{t_1}^{(i_1)}\ldots
d{\bf f}_{t_k}^{(i_k)}\right\},
\end{equation}

\vspace{4mm}
\noindent
where $i_1,\ldots,i_k = 1,\ldots,m;$
the expression 
$$
\sum\limits_{(j_1,\ldots,j_k)}
$$ 

\vspace{2mm}
\noindent
means the sum with respect to all
possible permutations 
$(j_1,\ldots,j_k)$. At the same time if 
$j_r$ swapped with $j_q$ in the permutation $(j_1,\ldots,j_k),$
then $i_r$ swapped with $i_q$ in the permutation
$(i_1,\ldots,i_k);$
another notations are the same as in Theorems {\rm 1, 2.}
}

\vspace{2mm}

Note that 

\vspace{-1mm}
$$
{\sf M}\left\{J[\psi^{(k)}]_{T,t}
\int\limits_t^T \phi_{j_k}(t_k)
\ldots
\int\limits_t^{t_{2}}\phi_{j_{1}}(t_{1})
d{\bf f}_{t_1}^{(i_1)}\ldots
d{\bf f}_{t_k}^{(i_k)}\right\}=C_{j_k\ldots j_1}.
$$

\vspace{4mm}

Then we can obtain the following particular cases 
of Theorem 8 for $k=1,\ldots,5$ and $i_1,\ldots,i_5=1,\ldots,m$
\cite{2020}-\cite{2020xx1}, \cite{8}, \cite{10}, \cite{15b}.

\vspace{5mm}

\centerline{\bf The case $k=1$}

\vspace{4mm}

$$
E_1^p
=I_1
-\sum_{j_1=0}^p
C_{j_1}^2.
$$

\vspace{8mm}

\centerline{\bf The case $k=2$}

\vspace{4mm}

(I).\ $i_1\ne i_2$:
\begin{equation}
\label{kruto1}
E_2^p
=I_2
-\sum_{j_1,j_2=0}^p
C_{j_2j_1}^2.
\end{equation}

\vspace{2mm}

(II).\ $i_1=i_2:$
\begin{equation}
\label{kruto2}
E_2^p
=I_2
-\sum_{j_1,j_2=0}^p
C_{j_2j_1}^2-
\sum_{j_1,j_2=0}^p
C_{j_2j_1}C_{j_1j_2}.
\end{equation}

\vspace{8mm}

\centerline{\bf The case $k=3$}

\vspace{4mm}

(I).\ $i_1\ne i_2, i_1\ne i_3, i_2\ne i_3:$
\begin{equation}
\label{kruto3}
E_3^p=I_3
-\sum_{j_1,j_2,j_3=0}^p C_{j_3j_2j_1}^2.
\end{equation}

\vspace{2mm}

(II).\ $i_1=i_2=i_3:$
\begin{equation}
\label{kruto4}
E_3^p = I_3 - \sum_{j_1,j_2,j_3=0}^{p}
C_{j_3j_2j_1}\Biggl(\sum\limits_{(j_1,j_2,j_3)}
C_{j_3j_2j_1}\Biggr).
\end{equation}

\vspace{2mm}

(III).1.\ $i_1=i_2\ne i_3:$
\begin{equation}
\label{kruto5}
E_3^p=I_3
-\sum_{j_1,j_2,j_3=0}^p C_{j_3j_2j_1}^2-
\sum_{j_1,j_2,j_3=0}^p C_{j_3j_1j_2}C_{j_3j_2j_1}.
\end{equation}

\vspace{2mm}

(III).2.\ $i_1\ne i_2=i_3:$
\begin{equation}
\label{kruto6}
E_3^p=I_3-
\sum_{j_1,j_2,j_3=0}^p C_{j_3j_2j_1}^2-
\sum_{j_1,j_2,j_3=0}^p C_{j_2j_3j_1}C_{j_3j_2j_1}.
\end{equation}

\vspace{2mm}

(III).3.\ $i_1=i_3\ne i_2:$
\begin{equation}
\label{kruto7}
E_3^p=I_3
-\sum_{j_1,j_2,j_3=0}^p C_{j_3j_2j_1}^2-
\sum_{j_1,j_2,j_3=0}^p C_{j_3j_2j_1}C_{j_1j_2j_3}.
\end{equation}

\vspace{8mm}

\centerline{\bf The case $k=4$}

\vspace{4mm}

(I).\ $i_1,\ldots,i_4$ are pairwise different:
\begin{equation}
\label{kruto8}
E_4^p= I_4 - \sum_{j_1,\ldots,j_4=0}^{p}C_{j_4\ldots j_1}^2.
\end{equation}

\vspace{4mm}

(II).\ $i_1=i_2=i_3=i_4$:
\begin{equation}
\label{kruto9}
E_4^p = I_4 -
 \sum_{j_1,\ldots,j_4=0}^{p}
C_{j_4\ldots j_1}\Biggl(\sum\limits_{(j_1,\ldots,j_4)}
C_{j_4\ldots j_1}\Biggr).
\end{equation}

\vspace{4mm}

(III).1.\ $i_1=i_2\ne i_3, i_4;\ i_3\ne i_4:$
\begin{equation}
\label{kruto10}
E^p_4 = I_4 - \sum_{j_1,\ldots,j_4=0}^{p}
C_{j_4\ldots j_1}\Biggl(\sum\limits_{(j_1,j_2)}
C_{j_4\ldots j_1}\Biggr).
\end{equation}

\vspace{4mm}

(III).2.\ $i_1=i_3\ne i_2, i_4;\ i_2\ne i_4:$
\begin{equation}
\label{kruto11}
E^p_4 = I_4 - \sum_{j_1,\ldots,j_4=0}^{p}
C_{j_4\ldots j_1}\Biggl(\sum\limits_{(j_1,j_3)}
C_{j_4\ldots j_1}\Biggr).
\end{equation}

\vspace{4mm}

(III).3.\ $i_1=i_4\ne i_2, i_3;\ i_2\ne i_3:$
\begin{equation}
\label{kruto12}
E^p_4 = I_4 - \sum_{j_1,\ldots,j_4=0}^{p}
C_{j_4\ldots j_1}\Biggl(\sum\limits_{(j_1,j_4)}
C_{j_4\ldots j_1}\Biggr).
\end{equation}

\vspace{4mm}

(III).4.\ $i_2=i_3\ne i_1, i_4;\ i_1\ne i_4:$
\begin{equation}
\label{kruto13}
E^p_4 = I_4 - \sum_{j_1,\ldots,j_4=0}^{p}
C_{j_4\ldots j_1}\Biggl(\sum\limits_{(j_2,j_3)}
C_{j_4\ldots j_1}\Biggr).
\end{equation}

\vspace{4mm}

(III).5.\ $i_2=i_4\ne i_1, i_3;\ i_1\ne i_3:$
\begin{equation}
\label{kruto14}
E^p_4 = I_4 - \sum_{j_1,\ldots,j_4=0}^{p}
C_{j_4\ldots j_1}\Biggl(\sum\limits_{(j_2,j_4)}
C_{j_4\ldots j_1}\Biggr).
\end{equation}

\vspace{4mm}

(III).6.\ $i_3=i_4\ne i_1, i_2;\ i_1\ne i_2:$
\begin{equation}
\label{kruto15}
E^p_4 = I_4 - \sum_{j_1,\ldots,j_4=0}^{p}
C_{j_4\ldots j_1}\Biggl(\sum\limits_{(j_3,j_4)}
C_{j_4\ldots j_1}\Biggr).
\end{equation}

\vspace{4mm}

(IV).1.\ $i_1=i_2=i_3\ne i_4$:
\begin{equation}
\label{kruto16}
E_4^p = I_4 -
 \sum_{j_1,\ldots,j_4=0}^{p}
C_{j_4\ldots j_1}\Biggl(\sum\limits_{(j_1,j_2,j_3)}
C_{j_4\ldots j_1}\Biggr).
\end{equation}

\vspace{4mm}

(IV).2.\ $i_2=i_3=i_4\ne i_1$:
\begin{equation}
\label{kruto17}
E_4^p = I_4 -
 \sum_{j_1,\ldots,j_4=0}^{p}
C_{j_4\ldots j_1}\Biggl(\sum\limits_{(j_2,j_3,j_4)}
C_{j_4\ldots j_1}\Biggr).
\end{equation}

\vspace{4mm}

(IV).3.\ $i_1=i_2=i_4\ne i_3$:
\begin{equation}
\label{kruto18}
E_4^p = I_4 -
 \sum_{j_1,\ldots,j_4=0}^{p}
C_{j_4\ldots j_1}\Biggl(\sum\limits_{(j_1,j_2,j_4)}
C_{j_4\ldots j_1}\Biggr).
\end{equation}

\vspace{4mm}

(IV).4.\ $i_1=i_3=i_4\ne i_2$:
\begin{equation}
\label{kruto19}
E_4^p = I_4 -
 \sum_{j_1,\ldots,j_4=0}^{p}
C_{j_4\ldots j_1}\Biggl(\sum\limits_{(j_1,j_3,j_4)}
C_{j_4\ldots j_1}\Biggr).
\end{equation}

\vspace{4mm}

(V).1.\ $i_1=i_2\ne i_3=i_4$:
\begin{equation}
\label{kruto20}
E^p_4 = I_4 - \sum_{j_1,\ldots,j_4=0}^{p}
C_{j_4\ldots j_1}\Biggl(\sum\limits_{(j_1,j_2)}\Biggl(
\sum\limits_{(j_3,j_4)}
C_{j_4\ldots j_1}\Biggr)\Biggr).
\end{equation}

\vspace{4mm}

(V).2.\ $i_1=i_3\ne i_2=i_4$:
\begin{equation}
\label{kruto21}
E^p_4 = I_4 - \sum_{j_1,\ldots,j_4=0}^{p}
C_{j_4\ldots j_1}\Biggl(\sum\limits_{(j_1,j_3)}\Biggl(
\sum\limits_{(j_2,j_4)}
C_{j_4\ldots j_1}\Biggr)\Biggr).
\end{equation}

\vspace{4mm}

(V).3.\ $i_1=i_4\ne i_2=i_3$:
\begin{equation}
\label{kruto22}
E^p_4 = I_4 - \sum_{j_1,\ldots,j_4=0}^{p}
C_{j_4\ldots j_1}\Biggl(\sum\limits_{(j_1,j_4)}\Biggl(
\sum\limits_{(j_2,j_3)}
C_{j_4\ldots j_1}\Biggr)\Biggr).
\end{equation}

\vspace{8mm}

\centerline{\bf The case $k=5$}

\vspace{3mm}

(I).\ $i_1,\ldots,i_5$ are pairwise different:
\begin{equation}
\label{kruto23}
E_5^p= I_5 - \sum_{j_1,\ldots,j_5=0}^{p}C_{j_5\ldots j_1}^2.
\end{equation}

\vspace{3mm}

(II).\ $i_1=i_2=i_3=i_4=i_5$:
\begin{equation}
\label{kruto24}
E_5^p = I_5 - \sum_{j_1,\ldots,j_5=0}^{p}
C_{j_5\ldots j_1}\Biggl(\sum\limits_{(j_1,\ldots,j_5)}
C_{j_5\ldots j_1}\Biggr).
\end{equation}

\vspace{6mm}

(III).1.\ $i_1=i_2\ne i_3,i_4,i_5$\ ($i_3,i_4,i_5$ are pairwise different):
\begin{equation}
\label{kruto25}
E^p_5 = I_5 - \sum_{j_1,\ldots,j_5=0}^{p}
C_{j_5\ldots j_1}\Biggl(\sum\limits_{(j_1,j_2)}
C_{j_5\ldots j_1}\Biggr).
\end{equation}

\vspace{6mm}

(III).2.\ $i_1=i_3\ne i_2,i_4,i_5$\ ($i_2,i_4,i_5$ are pairwise different):
\begin{equation}
\label{kruto26}
E^p_5 = I_5 - \sum_{j_1,\ldots,j_5=0}^{p}
C_{j_5\ldots j_1}\Biggl(\sum\limits_{(j_1,j_3)}
C_{j_5\ldots j_1}\Biggr).
\end{equation}

\vspace{6mm}

(III).3.\ $i_1=i_4\ne i_2,i_3,i_5$\ ($i_2,i_3,i_5$ are pairwise different):
\begin{equation}
\label{kruto27}
E^p_5 = I_5 - \sum_{j_1,\ldots,j_5=0}^{p}
C_{j_5\ldots j_1}\Biggl(\sum\limits_{(j_1,j_4)}
C_{j_5\ldots j_1}\Biggr).
\end{equation}

\vspace{6mm}

(III).4.\ $i_1=i_5\ne i_2,i_3,i_4$\ ($i_2,i_3,i_4$  are pairwise different):
\begin{equation}
\label{kruto28}
E^p_5 = I_5 - \sum_{j_1,\ldots,j_5=0}^{p}
C_{j_5\ldots j_1}\Biggl(\sum\limits_{(j_1,j_5)}
C_{j_5\ldots j_1}\Biggr).
\end{equation}

\vspace{6mm}

(III).5.\ $i_2=i_3\ne i_1,i_4,i_5$\ ($i_1,i_4,i_5$ are pairwise different):
\begin{equation}
\label{kruto29}
E^p_5 = I_5 - \sum_{j_1,\ldots,j_5=0}^{p}
C_{j_5\ldots j_1}\Biggl(\sum\limits_{(j_2,j_3)}
C_{j_5\ldots j_1}\Biggr).
\end{equation}

\vspace{6mm}

(III).6.\ $i_2=i_4\ne i_1,i_3,i_5$\ ($i_1,i_3,i_5$ are pairwise different):
\begin{equation}
\label{kruto30}
E^p_5 = I_5 - \sum_{j_1,\ldots,j_5=0}^{p}
C_{j_5\ldots j_1}\Biggl(\sum\limits_{(j_2,j_4)}
C_{j_5\ldots j_1}\Biggr).
\end{equation}

\vspace{6mm}

(III).7.\ $i_2=i_5\ne i_1,i_3,i_4$\ ($i_1,i_3,i_4$ are pairwise different):
\begin{equation}
\label{kruto31}
E^p_5 = I_5 - \sum_{j_1,\ldots,j_5=0}^{p}
C_{j_5\ldots j_1}\Biggl(\sum\limits_{(j_2,j_5)}
C_{j_5\ldots j_1}\Biggr).
\end{equation}

\vspace{6mm}

(III).8.\ $i_3=i_4\ne i_1,i_2,i_5$\ ($i_1,i_2,i_5$  are pairwise different):
\begin{equation}
\label{kruto32}
E^p_5 = I_5 - \sum_{j_1,\ldots,j_5=0}^{p}
C_{j_5\ldots j_1}\Biggl(\sum\limits_{(j_3,j_4)}
C_{j_5\ldots j_1}\Biggr).
\end{equation}

\vspace{6mm}

(III).9.\ $i_3=i_5\ne i_1,i_2,i_4$\ ($i_1,i_2,i_4$  are pairwise different):
\begin{equation}
\label{kruto33}
E^p_5 = I_5 - \sum_{j_1,\ldots,j_5=0}^{p}
C_{j_5\ldots j_1}\Biggl(\sum\limits_{(j_3,j_5)}
C_{j_5\ldots j_1}\Biggr).
\end{equation}

\vspace{6mm}

(III).10.\ $i_4=i_5\ne i_1,i_2,i_3$\ ($i_1,i_2,i_3$  are pairwise different):
\begin{equation}
\label{kruto34}
E^p_5 = I_5 - \sum_{j_1,\ldots,j_5=0}^{p}
C_{j_5\ldots j_1}\Biggl(\sum\limits_{(j_4,j_5)}
C_{j_5\ldots j_1}\Biggr).
\end{equation}

\vspace{6mm}

(IV).1.\ $i_1=i_2=i_3\ne i_4, i_5$\ $(i_4\ne i_5$):
\begin{equation}
\label{kruto35}
E^p_5 = I_5 - \sum_{j_1,\ldots,j_5=0}^{p}
C_{j_5\ldots j_1}\Biggl(\sum\limits_{(j_1,j_2,j_3)}
C_{j_5\ldots j_1}\Biggr).
\end{equation}

\vspace{6mm}

(IV).2.\ $i_1=i_2=i_4\ne i_3, i_5$\  $(i_3\ne i_5$):
\begin{equation}
\label{kruto36}
E^p_5 = I_5 - \sum_{j_1,\ldots,j_5=0}^{p}
C_{j_5\ldots j_1}\Biggl(\sum\limits_{(j_1,j_2,j_4)}
C_{j_5\ldots j_1}\Biggr).
\end{equation}

\vspace{6mm}

(IV).3.\ $i_1=i_2=i_5\ne i_3, i_4$\  $(i_3\ne i_4$):
\begin{equation}
\label{kruto37}
E^p_5 = I_5 - \sum_{j_1,\ldots,j_5=0}^{p}
C_{j_5\ldots j_1}\Biggl(\sum\limits_{(j_1,j_2,j_5)}
C_{j_5\ldots j_1}\Biggr).
\end{equation}

\vspace{6mm}

(IV).4.\ $i_2=i_3=i_4\ne i_1, i_5$\  $(i_1\ne i_5$):
\begin{equation}
\label{kruto38}
E^p_5 = I_5 - \sum_{j_1,\ldots,j_5=0}^{p}
C_{j_5\ldots j_1}\Biggl(\sum\limits_{(j_2,j_3,j_4)}
C_{j_5\ldots j_1}\Biggr).
\end{equation}

\vspace{6mm}
(IV).5.\ $i_2=i_3=i_5\ne i_1, i_4$\  $(i_1\ne i_4$):
\begin{equation}
\label{kruto39}
E^p_5 = I_5 - \sum_{j_1,\ldots,j_5=0}^{p}
C_{j_5\ldots j_1}\Biggl(\sum\limits_{(j_2,j_3,j_5)}
C_{j_5\ldots j_1}\Biggr).
\end{equation}

\vspace{6mm}

(IV).6.\ $i_2=i_4=i_5\ne i_1, i_3$\  $(i_1\ne i_3$):
\begin{equation}
\label{kruto40}
E^p_5 = I_5 - \sum_{j_1,\ldots,j_5=0}^{p}
C_{j_5\ldots j_1}\Biggl(\sum\limits_{(j_2,j_4,j_5)}
C_{j_5\ldots j_1}\Biggr).
\end{equation}

\vspace{6mm}

(IV).7.\ $i_3=i_4=i_5\ne i_1, i_2$\  $(i_1\ne i_2$):
\begin{equation}
\label{kruto41}
E^p_5 = I_5 - \sum_{j_1,\ldots,j_5=0}^{p}
C_{j_5\ldots j_1}\Biggl(\sum\limits_{(j_3,j_4,j_5)}
C_{j_5\ldots j_1}\Biggr).
\end{equation}

\vspace{6mm}

(IV).8.\ $i_1=i_3=i_5\ne i_2, i_4$\  $(i_2\ne i_4$):
\begin{equation}
\label{kruto42}
E^p_5 = I_5 - \sum_{j_1,\ldots,j_5=0}^{p}
C_{j_5\ldots j_1}\Biggl(\sum\limits_{(j_1,j_3,j_5)}
C_{j_5\ldots j_1}\Biggr).
\end{equation}

\vspace{6mm}

(IV).9.\ $i_1=i_3=i_4\ne i_2, i_5$\  $(i_2\ne i_5$):
\begin{equation}
\label{kruto43}
E^p_5 = I_5 - \sum_{j_1,\ldots,j_5=0}^{p}
C_{j_5\ldots j_1}\Biggl(\sum\limits_{(j_1,j_3,j_4)}
C_{j_5\ldots j_1}\Biggr).
\end{equation}

\vspace{6mm}

(IV).10.\ $i_1=i_4=i_5\ne i_2, i_3$\  $(i_2\ne i_3$):
\begin{equation}
\label{kruto44}
E^p_5 = I_5 - \sum_{j_1,\ldots,j_5=0}^{p}
C_{j_5\ldots j_1}\Biggl(\sum\limits_{(j_1,j_4,j_5)}
C_{j_5\ldots j_1}\Biggr).
\end{equation}

\vspace{6mm}

(V).1.\ $i_1=i_2=i_3=i_4\ne i_5$:
\begin{equation}
\label{kruto45}
E^p_5 = I_5 - \sum_{j_1,\ldots,j_5=0}^{p}
C_{j_5\ldots j_1}\Biggl(\sum\limits_{(j_1,j_2,j_3,j_4)}
C_{j_5\ldots j_1}\Biggr).
\end{equation}

\vspace{6mm}

(V).2.\ $i_1=i_2=i_3=i_5\ne i_4$:
\begin{equation}
\label{kruto46}
E^p_5 = I_5 - \sum_{j_1,\ldots,j_5=0}^{p}
C_{j_5\ldots j_1}\Biggl(\sum\limits_{(j_1,j_2,j_3,j_5)}
C_{j_5\ldots j_1}\Biggr).
\end{equation}

\vspace{6mm}

(V).3.\ $i_1=i_2=i_4=i_5\ne i_3$:
\begin{equation}
\label{kruto47}
E^p_5 = I_5 - \sum_{j_1,\ldots,j_5=0}^{p}
C_{j_5\ldots j_1}\Biggl(\sum\limits_{(j_1,j_2,j_4,j_5)}
C_{j_5\ldots j_1}\Biggr).
\end{equation}

\vspace{6mm}

(V).4.\ $i_1=i_3=i_4=i_5\ne i_2$:
\begin{equation}
\label{kruto48}
E^p_5 = I_5 - \sum_{j_1,\ldots,j_5=0}^{p}
C_{j_5\ldots j_1}\Biggl(\sum\limits_{(j_1,j_3,j_4,j_5)}
C_{j_5\ldots j_1}\Biggr).
\end{equation}

\vspace{6mm}

(V).5.\ $i_2=i_3=i_4=i_5\ne i_1$:
\begin{equation}
\label{kruto49}
E^p_5 = I_5 - \sum_{j_1,\ldots,j_5=0}^{p}
C_{j_5\ldots j_1}\Biggl(\sum\limits_{(j_2,j_3,j_4,j_5)}
C_{j_5\ldots j_1}\Biggr).
\end{equation}

\vspace{6mm}

(VI).1.\ $i_5\ne i_1=i_2\ne i_3=i_4\ne i_5$:
\begin{equation}
\label{kruto50}
E^p_5 = I_5 - \sum_{j_1,\ldots,j_5=0}^{p}
C_{j_5\ldots j_1}\Biggl(\sum\limits_{(j_1,j_2)}\Biggl(
\sum\limits_{(j_3,j_4)}
C_{j_5\ldots j_1}\Biggr)\Biggr).
\end{equation}

\vspace{6mm}

(VI).2.\ $i_5\ne i_1=i_3\ne i_2=i_4\ne i_5$:
\begin{equation}
\label{kruto51}
E^p_5 = I_5 - \sum_{j_1,\ldots,j_5=0}^{p}
C_{j_5\ldots j_1}\Biggl(\sum\limits_{(j_1,j_3)}\Biggl(
\sum\limits_{(j_2,j_4)}
C_{j_5\ldots j_1}\Biggr)\Biggr).
\end{equation}

\vspace{6mm}

(VI).3.\ $i_5\ne i_1=i_4\ne i_2=i_3\ne i_5$:
\begin{equation}
\label{kruto52}
E^p_5 = I_5 - \sum_{j_1,\ldots,j_5=0}^{p}
C_{j_5\ldots j_1}\Biggl(\sum\limits_{(j_1,j_4)}\Biggl(
\sum\limits_{(j_2,j_3)}
C_{j_5\ldots j_1}\Biggr)\Biggr).
\end{equation}

\vspace{6mm}

(VI).4.\ $i_4\ne i_1=i_2\ne i_3=i_5\ne i_4$:
\begin{equation}
\label{kruto53}
E^p_5 = I_5 - \sum_{j_1,\ldots,j_5=0}^{p}
C_{j_5\ldots j_1}\Biggl(\sum\limits_{(j_1,j_2)}\Biggl(
\sum\limits_{(j_3,j_5)}
C_{j_5\ldots j_1}\Biggr)\Biggr).
\end{equation}

\vspace{6mm}

(VI).5.\ $i_4\ne i_1=i_5\ne i_2=i_3\ne i_4$:
\begin{equation}
\label{kruto54}
E^p_5 = I_5 - \sum_{j_1,\ldots,j_5=0}^{p}
C_{j_5\ldots j_1}\Biggl(\sum\limits_{(j_1,j_5)}\Biggl(
\sum\limits_{(j_2,j_3)}
C_{j_5\ldots j_1}\Biggr)\Biggr).
\end{equation}

\vspace{6mm}

(VI).6.\ $i_4\ne i_2=i_5\ne i_1=i_3\ne i_4$:
\begin{equation}
\label{kruto55}
E^p_5 = I_5 - \sum_{j_1,\ldots,j_5=0}^{p}
C_{j_5\ldots j_1}\Biggl(\sum\limits_{(j_2,j_5)}\Biggl(
\sum\limits_{(j_1,j_3)}
C_{j_5\ldots j_1}\Biggr)\Biggr).
\end{equation}

\vspace{6mm}

(VI).7.\ $i_3\ne i_2=i_5\ne i_1=i_4\ne i_3$:
\begin{equation}
\label{kruto56}
E^p_5 = I_5 - \sum_{j_1,\ldots,j_5=0}^{p}
C_{j_5\ldots j_1}\Biggl(\sum\limits_{(j_2,j_5)}\Biggl(
\sum\limits_{(j_1,j_4)}
C_{j_5\ldots j_1}\Biggr)\Biggr).
\end{equation}

\vspace{6mm}

(VI).8.\ $i_3\ne i_1=i_2\ne i_4=i_5\ne i_3$:
\begin{equation}
\label{kruto57}
E^p_5 = I_5 - \sum_{j_1,\ldots,j_5=0}^{p}
C_{j_5\ldots j_1}\Biggl(\sum\limits_{(j_1,j_2)}\Biggl(
\sum\limits_{(j_4,j_5)}
C_{j_5\ldots j_1}\Biggr)\Biggr).
\end{equation}

\vspace{6mm}

(VI).9.\ $i_3\ne i_2=i_4\ne i_1=i_5\ne i_3$:
\begin{equation}
\label{kruto58}
E^p_5 = I_5 - \sum_{j_1,\ldots,j_5=0}^{p}
C_{j_5\ldots j_1}\Biggl(\sum\limits_{(j_2,j_4)}\Biggl(
\sum\limits_{(j_1,j_5)}
C_{j_5\ldots j_1}\Biggr)\Biggr).
\end{equation}

\vspace{6mm}

(VI).10.\ $i_2\ne i_1=i_4\ne i_3=i_5\ne i_2$:
\begin{equation}
\label{kruto59}
E^p_5 = I_5 - \sum_{j_1,\ldots,j_5=0}^{p}
C_{j_5\ldots j_1}\Biggl(\sum\limits_{(j_1,j_4)}\Biggl(
\sum\limits_{(j_3,j_5)}
C_{j_5\ldots j_1}\Biggr)\Biggr).
\end{equation}

\vspace{6mm}

(VI).11.\ $i_2\ne i_1=i_3\ne i_4=i_5\ne i_2$:
\begin{equation}
\label{kruto60}
E^p_5 = I_5 - \sum_{j_1,\ldots,j_5=0}^{p}
C_{j_5\ldots j_1}\Biggl(\sum\limits_{(j_1,j_3)}\Biggl(
\sum\limits_{(j_4,j_5)}
C_{j_5\ldots j_1}\Biggr)\Biggr).
\end{equation}

\vspace{6mm}

(VI).12.\ $i_2\ne i_1=i_5\ne i_3=i_4\ne i_2$:
\begin{equation}
\label{kruto61}
E^p_5 = I_5 - \sum_{j_1,\ldots,j_5=0}^{p}
C_{j_5\ldots j_1}\Biggl(\sum\limits_{(j_1,j_5)}\Biggl(
\sum\limits_{(j_3,j_4)}
C_{j_5\ldots j_1}\Biggr)\Biggr).
\end{equation}

\vspace{6mm}

(VI).13.\ $i_1\ne i_2=i_3\ne i_4=i_5\ne i_1$:
\begin{equation}
\label{kruto62}
E^p_5 = I_5 - \sum_{j_1,\ldots,j_5=0}^{p}
C_{j_5\ldots j_1}\Biggl(\sum\limits_{(j_2,j_3)}\Biggl(
\sum\limits_{(j_4,j_5)}
C_{j_5\ldots j_1}\Biggr)\Biggr).
\end{equation}

\vspace{6mm}

(VI).14.\ $i_1\ne i_2=i_4\ne i_3=i_5\ne i_1$:
\begin{equation}
\label{kruto63}
E^p_5 = I_5 - \sum_{j_1,\ldots,j_5=0}^{p}
C_{j_5\ldots j_1}\Biggl(\sum\limits_{(j_2,j_4)}\Biggl(
\sum\limits_{(j_3,j_5)}
C_{j_5\ldots j_1}\Biggr)\Biggr).
\end{equation}

\vspace{6mm}

(VI).15.\ $i_1\ne i_2=i_5\ne i_3=i_4\ne i_1$:
\begin{equation}
\label{kruto64}
E^p_5 = I_5 - \sum_{j_1,\ldots,j_5=0}^{p}
C_{j_5\ldots j_1}\Biggl(\sum\limits_{(j_2,j_5)}\Biggl(
\sum\limits_{(j_3,j_4)}
C_{j_5\ldots j_1}\Biggr)\Biggr).
\end{equation}

\vspace{6mm}

(VII).1.\ $i_1=i_2=i_3\ne i_4=i_5$:
\begin{equation}
\label{kruto65}
E^p_5 = I_5 - \sum_{j_1,\ldots,j_5=0}^{p}
C_{j_5\ldots j_1}\Biggl(\sum\limits_{(j_4,j_5)}\Biggl(
\sum\limits_{(j_1,j_2,j_3)}
C_{j_5\ldots j_1}\Biggr)\Biggr).
\end{equation}

\vspace{6mm}

(VII).2.\ $i_1=i_2=i_4\ne i_3=i_5$:
\begin{equation}
\label{kruto66}
E^p_5 = I_5 - \sum_{j_1,\ldots,j_5=0}^{p}
C_{j_5\ldots j_1}\Biggl(\sum\limits_{(j_3,j_5)}\Biggl(
\sum\limits_{(j_1,j_2,j_4)}
C_{j_5\ldots j_1}\Biggr)\Biggr).
\end{equation}

\vspace{6mm}

(VII).3.\ $i_1=i_2=i_5\ne i_3=i_4$:
\begin{equation}
\label{kruto67}
E_p = I - \sum_{j_1,\ldots,j_5=0}^{p}
C_{j_5\ldots j_1}\Biggl(\sum\limits_{(j_3,j_4)}\Biggl(
\sum\limits_{(j_1,j_2,j_5)}
C_{j_5\ldots j_1}\Biggr)\Biggr).
\end{equation}

\vspace{6mm}

(VII).4.\ $i_2=i_3=i_4\ne i_1=i_5$:
\begin{equation}
\label{kruto68}
E^p_5 = I_5 - \sum_{j_1,\ldots,j_5=0}^{p}
C_{j_5\ldots j_1}\Biggl(\sum\limits_{(j_1,j_5)}\Biggl(
\sum\limits_{(j_2,j_3,j_4)}
C_{j_5\ldots j_1}\Biggr)\Biggr).
\end{equation}

\vspace{6mm}

(VII).5.\ $i_2=i_3=i_5\ne i_1=i_4$:
\begin{equation}
\label{kruto69}
E^p_5 = I_5 - \sum_{j_1,\ldots,j_5=0}^{p}
C_{j_5\ldots j_1}\Biggl(\sum\limits_{(j_1,j_4)}\Biggl(
\sum\limits_{(j_2,j_3,j_5)}
C_{j_5\ldots j_1}\Biggr)\Biggr).
\end{equation}

\vspace{6mm}

(VII).6.\ $i_2=i_4=i_5\ne i_1=i_3$:
\begin{equation}
\label{kruto70}
E^p_5 = I_5 - \sum_{j_1,\ldots,j_5=0}^{p}
C_{j_5\ldots j_1}\Biggl(\sum\limits_{(j_1,j_3)}\Biggl(
\sum\limits_{(j_2,j_4,j_5)}
C_{j_5\ldots j_1}\Biggr)\Biggr).
\end{equation}

\vspace{6mm}

(VII).7.\ $i_3=i_4=i_5\ne i_1=i_2$:
\begin{equation}
\label{kruto71}
E^p_5 = I_5 - \sum_{j_1,\ldots,j_5=0}^{p}
C_{j_5\ldots j_1}\Biggl(\sum\limits_{(j_1,j_2)}\Biggl(
\sum\limits_{(j_3,j_4,j_5)}
C_{j_5\ldots j_1}\Biggr)\Biggr).
\end{equation}

\vspace{6mm}

(VII).8.\ $i_1=i_3=i_5\ne i_2=i_4$:
\begin{equation}
\label{kruto72}
E^p_5 = I_5 - \sum_{j_1,\ldots,j_5=0}^{p}
C_{j_5\ldots j_1}\Biggl(\sum\limits_{(j_2,j_4)}\Biggl(
\sum\limits_{(j_1,j_3,j_5)}
C_{j_5\ldots j_1}\Biggr)\Biggr).
\end{equation}

\vspace{6mm}

(VII).9.\ $i_1=i_3=i_4\ne i_2=i_5$:
\begin{equation}
\label{kruto73}
E^p_5 = I_5 - \sum_{j_1,\ldots,j_5=0}^{p}
C_{j_5\ldots j_1}\Biggl(\sum\limits_{(j_2,j_5)}\Biggl(
\sum\limits_{(j_1,j_3,j_4)}
C_{j_5\ldots j_1}\Biggr)\Biggr).
\end{equation}

\vspace{6mm}

(VII).10.\ $i_1=i_4=i_5\ne i_2=i_3$:
\begin{equation}
\label{kruto74}
E^p_5 = I_5 - \sum_{j_1,\ldots,j_5=0}^{p}
C_{j_5\ldots j_1}\Biggl(\sum\limits_{(j_2,j_3)}\Biggl(
\sum\limits_{(j_1,j_4,j_5)}
C_{j_5\ldots j_1}\Biggr)\Biggr).
\end{equation}

\vspace{5mm}

Obviously, the above formulas do not contain multiplier
factors $2!,$ $3!$, $4!,$ and $5!$ in contrast to 
the estimate
(\ref{star00011}). However, the number of the mentioned
conditions is quite large, which is inconvenient 
for practical calculations.

In the papers \cite{Kuz-Kuz} and \cite{Mikh-1},
it was proposed the hypothesis that all the formulas 
(\ref{kruto1})-(\ref{kruto74})
can be replaced by the following equalities

\vspace{1mm}
$$
E^p_2 = I_2 - \sum_{j_1,j_2=0}^{p}
C_{j_2j_1}^2,
$$

\vspace{3mm}
$$
E^p_3 = I_3 - \sum_{j_1,j_2,j_3=0}^{p}
C_{j_3j_2j_1}^2,
$$

\vspace{3mm}
$$
E^p_4 = I_4 - \sum_{j_1,\ldots,j_4=0}^{p}
C_{j_4 \ldots j_1}^2,
$$

\vspace{3mm}
$$
E^p_5 = I_5 - \sum_{j_1,\ldots,j_5=0}^{p}
C_{j_5 \ldots j_1}^2,
$$

\vspace{5mm}
\noindent
where $i_1,\ldots,i_5=1,\ldots,m.$

At that, the specified replacement will not lead 
to a noticeable loss of the mean-square accuracy  
of approximation of iterated Ito stochastic integrals
from the family (\ref{qqq1x}).

This paper is devoted to the detailed 
confirmation of the hypothesis from
\cite{Kuz-Kuz}, \cite{Mikh-1}
for the case of multiple Fourier--Legendre series.

It should be noted that unlike the method based on Theorems 1 and 2,
existing approaches to the mean-square approximation 
of iterated stochastic integrals 
(see, for example, \cite{1988}-\cite{1994},
\cite{Zapad-xx-1}-\cite{Zapad-xx-11}) do not allow choosing
different numbers $p$ for approximations of different 
iterated stochastic integrals. Moreover, the noted
approaches \cite{1988}-\cite{1994},
\cite{Zapad-xx-1}-\cite{Zapad-xx-11} exclude the possibility
for obtaining of approximate and exact expressions
for the mean-square approximation error similar to the formulas
(\ref{star00011}), (\ref{tttr11}).

\vspace{5mm}

\section{Approximations of Iterated 
Ito Stochastic 
Integrals from the Numerical Schemes (\ref{al1})--(\ref{al4})
Using Legendre Polynomials}

\vspace{5mm}

In this section, we consider the
approximations of the iterated Ito 
stochastic integrals (\ref{qqq1x})
of multiplicities 1 to 5 based on 
Theorems 1, 2 and multiple Fourier--Legendre series.

The numerical schemes (\ref{al1})--(\ref{al4})
contain the following set (see (\ref{qqq1x}))
of iterated Ito 
stochastic integrals 

\begin{equation}
\label{sm10}
I_{(0)T,t}^{(i_1)},\ \ \ I_{(1)T,t}^{(i_1)},\ \ \ I_{(2)T,t}^{(i_1)},\ \ \ 
I_{(00)T,t}^{(i_1 i_2)},\ \ \ I_{(10)T,t}^{(i_1 i_2)},\ \ \ 
I_{(01)T,t}^{(i_1 i_2)},\ \ \ I_{(000)T,t}^{(i_1 i_2 i_3)},\ \ \ 
I_{(0000)T,t}^{(i_1 i_2 i_3 i_4)},\ \ \
\end{equation}

\vspace{1mm}
\begin{equation}
\label{sm11}
I_{(00000)T,t}^{(i_1 i_2 i_3 i_4 i_5)},\ \ \
I_{(100)T,t}^{(i_1 i_2 i_3)},\ \ \ I_{(010)T,t}^{(i_1 i_2 i_3)},\ \ \ 
I_{(001)T,t}^{(i_1 i_2 i_3)}.
\end{equation}

\vspace{5mm}

Let us consider 
the complete orthonormal system of Legendre polynomials in the 
space $L_2([t,T])$

\begin{equation}
\label{4009}
\phi_j(x)=\sqrt{\frac{2j+1}{T-t}}P_j\left(\left(
x-\frac{T+t}{2}\right)\frac{2}{T-t}\right),\ \ \ j=0, 1, 2,\ldots,
\end{equation}

\vspace{4mm}
\noindent
where $P_j(x)$ is the Legendre polynomial

$$
P_j(x)=\frac{1}{2^j j!} \frac{d^j}{dx^j}\left(x^2-1\right)^j.
$$

\vspace{4mm}

Using Theorems 1, 2 and the system of functions (\ref{4009}),
we obtain the following formulas for numerical 
modeling of the stochastic integrals
(\ref{sm10}), (\ref{sm11})
\cite{1998a}-\cite{Mikh-1}

\vspace{2mm}
$$
I_{(0)T,t}^{(i_1)}=\sqrt{T-t}\zeta_0^{(i_1)},
$$

\vspace{3mm}
$$
I_{(1)T,t}^{(i_1)}=-\frac{(T-t)^{3/2}}{2}\left(\zeta_0^{(i_1)}+
\frac{1}{\sqrt{3}}\zeta_1^{(i_1)}\right),
$$

\vspace{3mm}
$$
I_{(2)T,t}^{(i_1)}=\frac{(T-t)^{5/2}}{3}\left(\zeta_0^{(i_1)}+
\frac{\sqrt{3}}{2}\zeta_1^{(i_1)}+
\frac{1}{2\sqrt{5}}\zeta_2^{(i_1)}\right),
$$

\vspace{4mm}
\begin{equation}
\label{nach1}
I_{(00)T,t}^{(i_1 i_2)q}=
\frac{T-t}{2}\left(\zeta_0^{(i_1)}\zeta_0^{(i_2)}+\sum_{i=1}^{q}
\frac{1}{\sqrt{4i^2-1}}\left(
\zeta_{i-1}^{(i_1)}\zeta_{i}^{(i_2)}-
\zeta_i^{(i_1)}\zeta_{i-1}^{(i_2)}\right) - {\bf 1}_{\{i_1=i_2\}}\right),
\end{equation}

\vspace{8mm}

$$
I_{(000)T,t}^{(i_1i_2i_3)q_1}
=
\sum_{j_1,j_2,j_3=0}^{q_1}
C_{j_3j_2j_1}^{000}
\Biggl(
\zeta_{j_1}^{(i_1)}\zeta_{j_2}^{(i_2)}\zeta_{j_3}^{(i_3)}
-{\bf 1}_{\{i_1=i_2\}}
{\bf 1}_{\{j_1=j_2\}}
\zeta_{j_3}^{(i_3)}-
\Biggr.
$$

\vspace{2mm}
$$
\Biggl.
-{\bf 1}_{\{i_2=i_3\}}
{\bf 1}_{\{j_2=j_3\}}
\zeta_{j_1}^{(i_1)}-
{\bf 1}_{\{i_1=i_3\}}
{\bf 1}_{\{j_1=j_3\}}
\zeta_{j_2}^{(i_2)}\Biggr),
$$

\vspace{7mm}

\begin{equation}
\label{ogo1}
I_{(10)T,t}^{(i_1 i_2)q_2}=
\sum_{j_1,j_2=0}^{q_2}
C_{j_2j_1}^{10}\Biggl(\zeta_{j_1}^{(i_1)}\zeta_{j_2}^{(i_2)}
-{\bf 1}_{\{i_1=i_2\}}
{\bf 1}_{\{j_1=j_2\}}\Biggr),
\end{equation}

\vspace{3mm}
\begin{equation}
\label{ogo2}
I_{(01)T,t}^{(i_1 i_2)q_2}=
\sum_{j_1,j_2=0}^{q_2}
C_{j_2j_1}^{01}\Biggl(\zeta_{j_1}^{(i_1)}\zeta_{j_2}^{(i_2)}
-{\bf 1}_{\{i_1=i_2\}}
{\bf 1}_{\{j_1=j_2\}}\Biggr),
\end{equation}

\vspace{10mm}

$$
I_{(0000)T,t}^{(i_1 i_2 i_3 i_4)q_3}
=
\sum_{j_1,j_2,j_3,j_4=0}^{q_3}
C_{j_4 j_3 j_2 j_1}^{0000}\Biggl(
\zeta_{j_1}^{(i_1)}\zeta_{j_2}^{(i_2)}\zeta_{j_3}^{(i_3)}\zeta_{j_4}^{(i_4)}
-\Biggr.
$$
$$
-
{\bf 1}_{\{i_1=i_2\}}
{\bf 1}_{\{j_1=j_2\}}
\zeta_{j_3}^{(i_3)}
\zeta_{j_4}^{(i_4)}
-
{\bf 1}_{\{i_1=i_3\}}
{\bf 1}_{\{j_1=j_3\}}
\zeta_{j_2}^{(i_2)}
\zeta_{j_4}^{(i_4)}-
$$
$$
-
{\bf 1}_{\{i_1=i_4\}}
{\bf 1}_{\{j_1=j_4\}}
\zeta_{j_2}^{(i_2)}
\zeta_{j_3}^{(i_3)}
-
{\bf 1}_{\{i_2=i_3\}}
{\bf 1}_{\{j_2=j_3\}}
\zeta_{j_1}^{(i_1)}
\zeta_{j_4}^{(i_4)}-
$$
$$
-
{\bf 1}_{\{i_2=i_4\}}
{\bf 1}_{\{j_2=j_4\}}
\zeta_{j_1}^{(i_1)}
\zeta_{j_3}^{(i_3)}
-
{\bf 1}_{\{i_3=i_4\}}
{\bf 1}_{\{j_3=j_4\}}
\zeta_{j_1}^{(i_1)}
\zeta_{j_2}^{(i_2)}+
$$
$$
+
{\bf 1}_{\{i_1=i_2\}}
{\bf 1}_{\{j_1=j_2\}}
{\bf 1}_{\{i_3=i_4\}}
{\bf 1}_{\{j_3=j_4\}}+
$$
$$
+
{\bf 1}_{\{i_1=i_3\}}
{\bf 1}_{\{j_1=j_3\}}
{\bf 1}_{\{i_2=i_4\}}
{\bf 1}_{\{j_2=j_4\}}+
$$
$$
+\Biggl.
{\bf 1}_{\{i_1=i_4\}}
{\bf 1}_{\{j_1=j_4\}}
{\bf 1}_{\{i_2=i_3\}}
{\bf 1}_{\{j_2=j_3\}}\Biggr),
$$

\vspace{10mm}

$$
I_{(00000)T,t}^{(i_1 i_2 i_3 i_4 i_5)q_4}
=
\sum_{j_1,j_2,j_3,j_4,j_5=0}^{q_4}
C_{j_5 j_4 j_3 j_2 j_1}^{00000}\Biggl(
\prod_{l=1}^5\zeta_{j_l}^{(i_l)}
-\Biggr.
$$
$$
-
{\bf 1}_{\{i_1=i_2\}}
{\bf 1}_{\{j_1=j_2\}}
\zeta_{j_3}^{(i_3)}
\zeta_{j_4}^{(i_4)}
\zeta_{j_5}^{(i_5)}-
{\bf 1}_{\{i_1=i_3\}}
{\bf 1}_{\{j_1=j_3\}}
\zeta_{j_2}^{(i_2)}
\zeta_{j_4}^{(i_4)}
\zeta_{j_5}^{(i_5)}-
$$
$$
-
{\bf 1}_{\{i_1=i_4\}}
{\bf 1}_{\{j_1=j_4\}}
\zeta_{j_2}^{(i_2)}
\zeta_{j_3}^{(i_3)}
\zeta_{j_5}^{(i_5)}-
{\bf 1}_{\{i_1=i_5\}}
{\bf 1}_{\{j_1=j_5\}}
\zeta_{j_2}^{(i_2)}
\zeta_{j_3}^{(i_3)}
\zeta_{j_4}^{(i_4)}-
$$
$$
-
{\bf 1}_{\{i_2=i_3\}}
{\bf 1}_{\{j_2=j_3\}}
\zeta_{j_1}^{(i_1)}
\zeta_{j_4}^{(i_4)}
\zeta_{j_5}^{(i_5)}-
{\bf 1}_{\{i_2=i_4\}}
{\bf 1}_{\{j_2=j_4\}}
\zeta_{j_1}^{(i_1)}
\zeta_{j_3}^{(i_3)}
\zeta_{j_5}^{(i_5)}-
$$
$$
-
{\bf 1}_{\{i_2=i_5\}}
{\bf 1}_{\{j_2=j_5\}}
\zeta_{j_1}^{(i_1)}
\zeta_{j_3}^{(i_3)}
\zeta_{j_4}^{(i_4)}
-{\bf 1}_{\{i_3=i_4\}}
{\bf 1}_{\{j_3=j_4\}}
\zeta_{j_1}^{(i_1)}
\zeta_{j_2}^{(i_2)}
\zeta_{j_5}^{(i_5)}-
$$
$$
-
{\bf 1}_{\{i_3=i_5\}}
{\bf 1}_{\{j_3=j_5\}}
\zeta_{j_1}^{(i_1)}
\zeta_{j_2}^{(i_2)}
\zeta_{j_4}^{(i_4)}
-{\bf 1}_{\{i_4=i_5\}}
{\bf 1}_{\{j_4=j_5\}}
\zeta_{j_1}^{(i_1)}
\zeta_{j_2}^{(i_2)}
\zeta_{j_3}^{(i_3)}+
$$
$$
+
{\bf 1}_{\{i_1=i_2\}}
{\bf 1}_{\{j_1=j_2\}}
{\bf 1}_{\{i_3=i_4\}}
{\bf 1}_{\{j_3=j_4\}}\zeta_{j_5}^{(i_5)}+
{\bf 1}_{\{i_1=i_2\}}
{\bf 1}_{\{j_1=j_2\}}
{\bf 1}_{\{i_3=i_5\}}
{\bf 1}_{\{j_3=j_5\}}\zeta_{j_4}^{(i_4)}+
$$
$$
+
{\bf 1}_{\{i_1=i_2\}}
{\bf 1}_{\{j_1=j_2\}}
{\bf 1}_{\{i_4=i_5\}}
{\bf 1}_{\{j_4=j_5\}}\zeta_{j_3}^{(i_3)}+
{\bf 1}_{\{i_1=i_3\}}
{\bf 1}_{\{j_1=j_3\}}
{\bf 1}_{\{i_2=i_4\}}
{\bf 1}_{\{j_2=j_4\}}\zeta_{j_5}^{(i_5)}+
$$
$$
+
{\bf 1}_{\{i_1=i_3\}}
{\bf 1}_{\{j_1=j_3\}}
{\bf 1}_{\{i_2=i_5\}}
{\bf 1}_{\{j_2=j_5\}}\zeta_{j_4}^{(i_4)}+
{\bf 1}_{\{i_1=i_3\}}
{\bf 1}_{\{j_1=j_3\}}
{\bf 1}_{\{i_4=i_5\}}
{\bf 1}_{\{j_4=j_5\}}\zeta_{j_2}^{(i_2)}+
$$
$$
+
{\bf 1}_{\{i_1=i_4\}}
{\bf 1}_{\{j_1=j_4\}}
{\bf 1}_{\{i_2=i_3\}}
{\bf 1}_{\{j_2=j_3\}}\zeta_{j_5}^{(i_5)}+
{\bf 1}_{\{i_1=i_4\}}
{\bf 1}_{\{j_1=j_4\}}
{\bf 1}_{\{i_2=i_5\}}
{\bf 1}_{\{j_2=j_5\}}\zeta_{j_3}^{(i_3)}+
$$
$$
+
{\bf 1}_{\{i_1=i_4\}}
{\bf 1}_{\{j_1=j_4\}}
{\bf 1}_{\{i_3=i_5\}}
{\bf 1}_{\{j_3=j_5\}}\zeta_{j_2}^{(i_2)}+
{\bf 1}_{\{i_1=i_5\}}
{\bf 1}_{\{j_1=j_5\}}
{\bf 1}_{\{i_2=i_3\}}
{\bf 1}_{\{j_2=j_3\}}\zeta_{j_4}^{(i_4)}+
$$
$$
+
{\bf 1}_{\{i_1=i_5\}}
{\bf 1}_{\{j_1=j_5\}}
{\bf 1}_{\{i_2=i_4\}}
{\bf 1}_{\{j_2=j_4\}}\zeta_{j_3}^{(i_3)}+
{\bf 1}_{\{i_1=i_5\}}
{\bf 1}_{\{j_1=j_5\}}
{\bf 1}_{\{i_3=i_4\}}
{\bf 1}_{\{j_3=j_4\}}\zeta_{j_2}^{(i_2)}+
$$
$$
+
{\bf 1}_{\{i_2=i_3\}}
{\bf 1}_{\{j_2=j_3\}}
{\bf 1}_{\{i_4=i_5\}}
{\bf 1}_{\{j_4=j_5\}}\zeta_{j_1}^{(i_1)}+
{\bf 1}_{\{i_2=i_4\}}
{\bf 1}_{\{j_2=j_4\}}
{\bf 1}_{\{i_3=i_5\}}
{\bf 1}_{\{j_3=j_5\}}\zeta_{j_1}^{(i_1)}+
$$
$$
+\Biggl.
{\bf 1}_{\{i_2=i_5\}}
{\bf 1}_{\{j_2=j_5\}}
{\bf 1}_{\{i_3=i_4\}}
{\bf 1}_{\{j_3=j_4\}}\zeta_{j_1}^{(i_1)}\Biggr),
$$

\vspace{9mm}

$$
I_{(001)T,t}^{(i_1i_2i_3)q_5}
=
\sum_{j_1,j_2,j_3=0}^{q_5}
C_{j_3j_2j_1}^{001}\Biggl(
\zeta_{j_1}^{(i_1)}\zeta_{j_2}^{(i_2)}\zeta_{j_3}^{(i_3)}
-{\bf 1}_{\{i_1=i_2\}}
{\bf 1}_{\{j_1=j_2\}}
\zeta_{j_3}^{(i_3)}-
\Biggr.
$$

\vspace{2mm}
$$
\Biggl.
-{\bf 1}_{\{i_2=i_3\}}
{\bf 1}_{\{j_2=j_3\}}
\zeta_{j_1}^{(i_1)}-
{\bf 1}_{\{i_1=i_3\}}
{\bf 1}_{\{j_1=j_3\}}
\zeta_{j_2}^{(i_2)}\Biggr),
$$

\vspace{8mm}

$$
I_{(010)T,t}^{(i_1i_2i_3)q_6}
=
\sum_{j_1,j_2,j_3=0}^{q_6}
C_{j_3j_2j_1}^{010}\Biggl(
\zeta_{j_1}^{(i_1)}\zeta_{j_2}^{(i_2)}\zeta_{j_3}^{(i_3)}
-{\bf 1}_{\{i_1=i_2\}}
{\bf 1}_{\{j_1=j_2\}}
\zeta_{j_3}^{(i_3)}-
\Biggr.
$$

\vspace{2mm}
$$
\Biggl.
-{\bf 1}_{\{i_2=i_3\}}
{\bf 1}_{\{j_2=j_3\}}
\zeta_{j_1}^{(i_1)}-
{\bf 1}_{\{i_1=i_3\}}
{\bf 1}_{\{j_1=j_3\}}
\zeta_{j_2}^{(i_2)}\Biggr),
$$

\vspace{8mm}

$$
I_{(100)T,t}^{(i_1i_2i_3)q_{7}}
=
\sum_{j_1,j_2,j_3=0}^{q_{7}}
C_{j_3j_2j_1}^{100}\Biggl(
\zeta_{j_1}^{(i_1)}\zeta_{j_2}^{(i_2)}\zeta_{j_3}^{(i_3)}
-{\bf 1}_{\{i_1=i_2\}}
{\bf 1}_{\{j_1=j_2\}}
\zeta_{j_3}^{(i_3)}-
\Biggr.
$$

\vspace{2mm}
$$
\Biggl.
-{\bf 1}_{\{i_2=i_3\}}
{\bf 1}_{\{j_2=j_3\}}
\zeta_{j_1}^{(i_1)}-
{\bf 1}_{\{i_1=i_3\}}
{\bf 1}_{\{j_1=j_3\}}
\zeta_{j_2}^{(i_2)}\Biggr),
$$

\vspace{8mm}
\noindent
where ${\bf 1}_A$ is the indicator of the set $A,$

\vspace{2mm}
\begin{equation}
\label{ma1}
C_{j_3j_2j_1}^{000}
=\frac{\sqrt{(2j_1+1)(2j_2+1)(2j_3+1)}}{8}(T-t)^{3/2}\bar
C_{j_3j_2j_1}^{000},
\end{equation}

\vspace{2mm}
\begin{equation}
\label{ma2}
C_{j_2j_1}^{01}
=\frac{\sqrt{(2j_1+1)(2j_2+1)}}{8}(T-t)^{2}\bar
C_{j_2j_1}^{01},
\end{equation}

\vspace{2mm}
\begin{equation}
\label{ma3}
C_{j_2j_1}^{10}
=\frac{\sqrt{(2j_1+1)(2j_2+1)}}{8}(T-t)^{2}\bar
C_{j_2j_1}^{10},
\end{equation}

\vspace{2mm}
\begin{equation}
\label{ma4}
C_{j_4j_3j_2j_1}^{0000}=
\frac{\sqrt{(2j_1+1)(2j_2+1)(2j_3+1)(2j_4+1)}}{16}(T-t)^{2}\bar
C_{j_4j_3j_2j_1}^{0000},
\end{equation}

\vspace{2mm}
\begin{equation}
\label{ma5}
C_{j_3j_2j_1}^{001}
=\frac{\sqrt{(2j_1+1)(2j_2+1)(2j_3+1)}}{16}(T-t)^{5/2}\bar
C_{j_3j_2j_1}^{001},
\end{equation}

\vspace{2mm}
\begin{equation}
\label{ma6}
C_{j_3j_2j_1}^{010}
=\frac{\sqrt{(2j_1+1)(2j_2+1)(2j_3+1)}}{16}(T-t)^{5/2}\bar
C_{j_3j_2j_1}^{010},
\end{equation}

\vspace{2mm}
\begin{equation}
\label{ma7}
C_{j_3j_2j_1}^{100}
=\frac{\sqrt{(2j_1+1)(2j_2+1)(2j_3+1)}}{16}(T-t)^{5/2}\bar
C_{j_3j_2j_1}^{100},
\end{equation}

\vspace{2mm}
\begin{equation}
\label{ma8}
C_{j_5j_4 j_3 j_2 j_1}^{00000}=
\frac{\sqrt{(2j_1+1)(2j_2+1)(2j_3+1)(2j_4+1)(2j_5+1)}}{32}(T-t)^{5/2}\bar
C_{j_5j_4 j_3 j_2 j_1}^{00000},
\end{equation}

\vspace{5mm}
\noindent
where
\begin{equation}
\label{ma14}
\bar C_{j_3j_2j_1}^{000}=\int\limits_{-1}^{1}P_{j_3}(z)
\int\limits_{-1}^{z}P_{j_2}(y)
\int\limits_{-1}^{y}
P_{j_1}(x)dx dy dz,
\end{equation}

\vspace{2mm}
\begin{equation}
\label{ma15}
\bar C_{j_2j_1}^{01}=-
\int\limits_{-1}^{1}(1+y)P_{j_2}(y)
\int\limits_{-1}^{y}
P_{j_1}(x)dx dy,
\end{equation}

\vspace{2mm}
\begin{equation}
\label{ma16}
\bar C_{j_2j_1}^{10}=-
\int\limits_{-1}^{1}P_{j_2}(y)
\int\limits_{-1}^{y}
(1+x)P_{j_1}(x)dx dy,
\end{equation}

\vspace{2mm}
\begin{equation}
\label{ma17}
\bar C_{j_4j_3j_2j_1}^{0000}=\int\limits_{-1}^{1}P_{j_4}(u)
\int\limits_{-1}^{u}P_{j_3}(z)
\int\limits_{-1}^{z}P_{j_2}(y)
\int\limits_{-1}^{y}
P_{j_1}(x)dx dy dz du,
\end{equation}

\vspace{2mm}
\begin{equation}
\label{ma18}
\bar C_{j_3j_2j_1}^{001}=-
\int\limits_{-1}^{1}P_{j_3}(z)(z+1)
\int\limits_{-1}^{z}P_{j_2}(y)
\int\limits_{-1}^{y}
P_{j_1}(x)dx dy dz,
\end{equation}

\vspace{2mm}
\begin{equation}
\label{ma19}
\bar C_{j_3j_2j_1}^{010}=-
\int\limits_{-1}^{1}P_{j_3}(z)
\int\limits_{-1}^{z}P_{j_2}(y)(y+1)
\int\limits_{-1}^{y}
P_{j_1}(x)dx dy dz,
\end{equation}

\vspace{2mm}
\begin{equation}
\label{ma20}
\bar C_{j_3j_2j_1}^{100}=-
\int\limits_{-1}^{1}P_{j_3}(z)
\int\limits_{-1}^{z}P_{j_2}(y)
\int\limits_{-1}^{y}
P_{j_1}(x)(x+1)dx dy dz,
\end{equation}

\vspace{2mm}
\begin{equation}
\label{ma21}
\bar C_{j_5j_4 j_3 j_2 j_1}^{00000}=
\int\limits_{-1}^{1}P_{j_5}(v)
\int\limits_{-1}^{v}P_{j_4}(u)
\int\limits_{-1}^{u}P_{j_3}(z)
\int\limits_{-1}^{z}P_{j_2}(y)
\int\limits_{-1}^{y}
P_{j_1}(x)dx dy dz du dv;
\end{equation}

\vspace{6mm}
\noindent
another notations are the same as in Theorems 1, 2.

\vspace{5mm}

\section{Optimization of Approximations of Iterated 
Ito Stochastic 
Integrals from the Numerical Schemes (\ref{al1})--(\ref{al4})}

\vspace{5mm}

This section is devoted to the optimization
of approximations of iterated 
Ito stochastic 
integrals from the numerical schemes (\ref{al1})--(\ref{al4}).
More precisely, we discuss how minimize the numbers 
$q, q_1, q_2,\ldots, q_{7}$ from Sect.~5.

Let us combine the relations (\ref{nach1}), (\ref{ma1})-(\ref{ma21}) with 
(\ref{kruto1})-(\ref{kruto74}). Thus, we have the following
formulas.

\vspace{5mm}

\centerline{\bf The case $k=2$ for the integral $I_{(00)T,t}^{(i_1 i_2)}$}

\vspace{5mm}

{\bf 2}.1.a.\ $i_1\ne i_2$:
\begin{equation}
\label{prod1}
E_2^p
=(T-t)^2\left(\frac{1}{2}-\frac{1}{4}-\frac{1}{2}\sum_{i=1}^p
\frac{1}{4i^2-1}\right)
=
\frac{(T-t)^2}{2}\left(\frac{1}{2}-\sum_{i=1}^p
\frac{1}{4i^2-1}\right).
\end{equation}

\vspace{4mm}

{\bf 2}.2.a.\ $i_1=i_2:$
$$
E_2^p
=(T-t)^2\left(\frac{1}{2}-\frac{1}{16}
\sum_{j_1,j_2=0}^p (2j_1+1)(2j_2+1)
\bar C_{j_2j_1}^{00}\left(\sum\limits_{(j_1,j_2)}\bar C_{j_2j_1}^{00}
\right)\right)\equiv 0.
$$

\vspace{7mm}

\centerline{\bf The case $k=2$ for the integral $I_{(01)T,t}^{(i_1 i_2)}$}

\vspace{5mm}

{\bf 2}.1.b.\ $i_1\ne i_2$:
\begin{equation}
\label{prod3}
E_2^p
=(T-t)^4\left(\frac{1}{4}-\frac{1}{64}
\sum_{j_1,j_2=0}^p (2j_1+1)(2j_2+1)
\left(\bar C_{j_2j_1}^{01}\right)^2
\right).
\end{equation}

\vspace{4mm}

{\bf 2}.2.b.\ $i_1=i_2:$
\begin{equation}
\label{prod4}
E_2^p
=(T-t)^4\left(\frac{1}{4}-\frac{1}{64}
\sum_{j_1,j_2=0}^p (2j_1+1)(2j_2+1)
\bar C_{j_2j_1}^{01}\left(\sum\limits_{(j_1,j_2)}\bar C_{j_2j_1}^{01}
\right)\right).
\end{equation}

\vspace{7mm}

\centerline{\bf The case $k=2$ for the integral $I_{(10)T,t}^{(i_1 i_2)}$}

\vspace{5mm}

{\bf 2}.1.c.\ $i_1\ne i_2$:
\begin{equation}
\label{prod5}
E_2^p
=(T-t)^4\left(\frac{1}{12}-\frac{1}{64}
\sum_{j_1,j_2=0}^p (2j_1+1)(2j_2+1)
\left(\bar C_{j_2j_1}^{10}\right)^2
\right).
\end{equation}

\vspace{4mm}

{\bf 2}.2.c.\ $i_1=i_2:$
\begin{equation}
\label{prod6}
E_2^p
=(T-t)^4\left(\frac{1}{12}-\frac{1}{64}
\sum_{j_1,j_2=0}^p (2j_1+1)(2j_2+1)
\bar C_{j_2j_1}^{10}\left(\sum\limits_{(j_1,j_2)}\bar C_{j_2j_1}^{10}
\right)\right).
\end{equation}

\vspace{10mm}

\centerline{\bf The case $k=3$
for the integral $I_{(000)T,t}^{(i_1 i_2 i_3)}$}

\vspace{5mm}

{\bf 3}.1.a.\ $i_1\ne i_2, i_1\ne i_3, i_2\ne i_3:$
\begin{equation}
\label{prod7}
E^p_3 = (T-t)^3\left(\frac{1}{6}-\frac{1}{64}
\sum_{j_1,j_2,j_3=0}^{p}(2j_1+1)(2j_2+1)(2j_3+1)
\left(\bar C_{j_3j_2j_1}^{000}\right)^2\right).
\end{equation}

\vspace{4mm}

{\bf 3}.2.a.\ $i_1=i_2=i_3:$
$$
E^p_3 = (T-t)^3\left(\frac{1}{6}-\frac{1}{64}
\sum_{j_1,j_2,j_3=0}^{p}(2j_1+1)(2j_2+1)(2j_3+1)
\bar C_{j_3j_2j_1}^{000}
\left(\sum\limits_{(j_1,j_2,j_3)}
\bar C_{j_3j_2j_1}^{000}\right)\right)\equiv 0.
$$

\vspace{4mm}

{\bf 3}.3.1.a.\ $i_1=i_2\ne i_3:$
\begin{equation}
\label{prod8}
E^p_3 = (T-t)^3\left(\frac{1}{6}-\frac{1}{64}
\sum_{j_1,j_2,j_3=0}^{p}(2j_1+1)(2j_2+1)(2j_3+1)
\left(\left(\bar C_{j_3j_2j_1}^{000}\right)^2+
\bar C_{j_3j_1j_2}^{000}\bar C_{j_3j_2j_1}^{000}\right)\right).
\end{equation}

\vspace{4mm}

{\bf 3}.3.2.a.\ $i_1\ne i_2=i_3:$
\begin{equation}
\label{prod9}
E^p_3 = (T-t)^3\left(\frac{1}{6}-\frac{1}{64}
\sum_{j_1,j_2,j_3=0}^{p}(2j_1+1)(2j_2+1)(2j_3+1)
\left(\left(\bar C_{j_3j_2j_1}^{000}\right)^2+
\bar C_{j_2j_3j_1}^{000}\bar C_{j_3j_2j_1}^{000}\right)\right).
\end{equation}

\vspace{4mm}

{\bf 3}.3.3.a.\ $i_1=i_3\ne i_2:$
\begin{equation}
\label{prod10}
E^p_3 = (T-t)^3\left(\frac{1}{6}-\frac{1}{64}
\sum_{j_1,j_2,j_3=0}^{p}(2j_1+1)(2j_2+1)(2j_3+1)
\left(\left(\bar C_{j_3j_2j_1}^{000}\right)^2+
\bar C_{j_3j_2j_1}^{000}\bar C_{j_1j_2j_3}^{000}\right)\right).
\end{equation}

\vspace{10mm}

\centerline{\bf The case $k=3$
for the integral $I_{(001)T,t}^{(i_1 i_2 i_3)}$}

\vspace{5mm}

{\bf 3}.1.b.\ $i_1\ne i_2, i_1\ne i_3, i_2\ne i_3:$
\begin{equation}
\label{prod11}
E^p_3 = (T-t)^5\left(\frac{1}{10}-\frac{1}{256}
\sum_{j_1,j_2,j_3=0}^{p}(2j_1+1)(2j_2+1)(2j_3+1)
\left(\bar C_{j_3j_2j_1}^{001}\right)^2\right).
\end{equation}

\vspace{4mm}

{\bf 3}.2.b.\ $i_1=i_2=i_3:$
\begin{equation}
\label{prod12}
E^p_3 = (T-t)^5\left(\frac{1}{10}-\frac{1}{256}
\sum_{j_1,j_2,j_3=0}^{p}(2j_1+1)(2j_2+1)(2j_3+1)
\bar C_{j_3j_2j_1}^{001}
\left(\sum\limits_{(j_1,j_2,j_3)}
\bar C_{j_3j_2j_1}^{001}\right)\right).
\end{equation}

\vspace{4mm}

{\bf 3}.3.1.b.\ $i_1=i_2\ne i_3:$
\begin{equation}
\label{prod13}
E^p_3 = (T-t)^5\left(\frac{1}{10}-\frac{1}{256}
\sum_{j_1,j_2,j_3=0}^{p}(2j_1+1)(2j_2+1)(2j_3+1)
\left(\left(\bar C_{j_3j_2j_1}^{001}\right)^2+
\bar C_{j_3j_1j_2}^{001}\bar C_{j_3j_2j_1}^{001}\right)\right).
\end{equation}

\vspace{4mm}

{\bf 3}.3.2.b.\ $i_1\ne i_2=i_3:$
\begin{equation}
\label{prod14}
E^p_3 = (T-t)^5\left(\frac{1}{10}-\frac{1}{256}
\sum_{j_1,j_2,j_3=0}^{p}(2j_1+1)(2j_2+1)(2j_3+1)
\left(\left(\bar C_{j_3j_2j_1}^{001}\right)^2+
\bar C_{j_2j_3j_1}^{001}\bar C_{j_3j_2j_1}^{001}\right)\right).
\end{equation}

\vspace{4mm}

{\bf 3}.3.3.b.\ $i_1=i_3\ne i_2:$
\begin{equation}
\label{prod15}
E^p_3 = (T-t)^5\left(\frac{1}{10}-\frac{1}{256}
\sum_{j_1,j_2,j_3=0}^{p}(2j_1+1)(2j_2+1)(2j_3+1)
\left(\left(\bar C_{j_3j_2j_1}^{001}\right)^2+
\bar C_{j_3j_2j_1}^{001}\bar C_{j_1j_2j_3}^{001}\right)\right).
\end{equation}

\vspace{10mm}

\centerline{\bf The case $k=3$
for the integral $I_{(010)T,t}^{(i_1 i_2 i_3)}$}

\vspace{5mm}

{\bf 3}.1.c.\ $i_1\ne i_2, i_1\ne i_3, i_2\ne i_3:$
\begin{equation}
\label{prod16}
E^p_3 = (T-t)^5\left(\frac{1}{20}-\frac{1}{256}
\sum_{j_1,j_2,j_3=0}^{p}(2j_1+1)(2j_2+1)(2j_3+1)
\left(\bar C_{j_3j_2j_1}^{010}\right)^2\right).
\end{equation}

\vspace{4mm}

{\bf 3}.2.c.\ $i_1=i_2=i_3:$
\begin{equation}
\label{prod17}
E^p_3 = (T-t)^5\left(\frac{1}{20}-\frac{1}{256}
\sum_{j_1,j_2,j_3=0}^{p}(2j_1+1)(2j_2+1)(2j_3+1)
\bar C_{j_3j_2j_1}^{010}
\left(\sum\limits_{(j_1,j_2,j_3)}
\bar C_{j_3j_2j_1}^{010}\right)\right).
\end{equation}

\vspace{4mm}

{\bf 3}.3.1.c.\ $i_1=i_2\ne i_3:$
\begin{equation}
\label{prod18}
E^p_3 = (T-t)^5\left(\frac{1}{20}-\frac{1}{256}
\sum_{j_1,j_2,j_3=0}^{p}(2j_1+1)(2j_2+1)(2j_3+1)
\left(\left(\bar C_{j_3j_2j_1}^{010}\right)^2+
\bar C_{j_3j_1j_2}^{010}\bar C_{j_3j_2j_1}^{010}\right)\right).
\end{equation}

\vspace{4mm}

{\bf 3}.3.2.c.\ $i_1\ne i_2=i_3:$
\begin{equation}
\label{prod19}
E^p_3 = (T-t)^5\left(\frac{1}{20}-\frac{1}{256}
\sum_{j_1,j_2,j_3=0}^{p}(2j_1+1)(2j_2+1)(2j_3+1)
\left(\left(\bar C_{j_3j_2j_1}^{010}\right)^2+
\bar C_{j_2j_3j_1}^{010}\bar C_{j_3j_2j_1}^{010}\right)\right).
\end{equation}

\vspace{4mm}

{\bf 3}.3.3.c.\ $i_1=i_3\ne i_2:$
\begin{equation}
\label{prod20}
E^p_3 = (T-t)^5\left(\frac{1}{20}-\frac{1}{256}
\sum_{j_1,j_2,j_3=0}^{p}(2j_1+1)(2j_2+1)(2j_3+1)
\left(\left(\bar C_{j_3j_2j_1}^{010}\right)^2+
\bar C_{j_3j_2j_1}^{010}\bar C_{j_1j_2j_3}^{010}\right)\right).
\end{equation}

\vspace{10mm}

\centerline{\bf The case $k=3$
for the integral $I_{(100)T,t}^{(i_1 i_2 i_3)}$}

\vspace{5mm}

{\bf 3}.1.d.\ $i_1\ne i_2, i_1\ne i_3, i_2\ne i_3:$
\begin{equation}
\label{prod21}
E^p_3 = (T-t)^5\left(\frac{1}{60}-\frac{1}{256}
\sum_{j_1,j_2,j_3=0}^{p}(2j_1+1)(2j_2+1)(2j_3+1)
\left(\bar C_{j_3j_2j_1}^{100}\right)^2\right).
\end{equation}

\vspace{4mm}

{\bf 3}.2.d.\ $i_1=i_2=i_3:$
\begin{equation}
\label{prod22}
E^p_3 = (T-t)^5\left(\frac{1}{60}-\frac{1}{256}
\sum_{j_1,j_2,j_3=0}^{p}(2j_1+1)(2j_2+1)(2j_3+1)
\bar C_{j_3j_2j_1}^{100}
\left(\sum\limits_{(j_1,j_2,j_3)}
\bar C_{j_3j_2j_1}^{100}\right)\right).
\end{equation}

\vspace{4mm}

{\bf 3}.3.1.d.\ $i_1=i_2\ne i_3:$
\begin{equation}
\label{prod23}
E^p_3 = (T-t)^5\left(\frac{1}{60}-\frac{1}{256}
\sum_{j_1,j_2,j_3=0}^{p}(2j_1+1)(2j_2+1)(2j_3+1)
\left(\left(\bar C_{j_3j_2j_1}^{100}\right)^2+
\bar C_{j_3j_1j_2}^{100}\bar C_{j_3j_2j_1}^{100}\right)\right).
\end{equation}

\vspace{4mm}

{\bf 3}.3.2.d.\ $i_1\ne i_2=i_3:$
\begin{equation}
\label{prod24}
E^p_3 = (T-t)^5\left(\frac{1}{60}-\frac{1}{256}
\sum_{j_1,j_2,j_3=0}^{p}(2j_1+1)(2j_2+1)(2j_3+1)
\left(\left(\bar C_{j_3j_2j_1}^{100}\right)^2+
\bar C_{j_2j_3j_1}^{100}\bar C_{j_3j_2j_1}^{100}\right)\right).
\end{equation}

\vspace{4mm}

{\bf 3}.3.3.d.\ $i_1=i_3\ne i_2:$
\begin{equation}
\label{prod25}
E^p_3 = (T-t)^5\left(\frac{1}{60}-\frac{1}{256}
\sum_{j_1,j_2,j_3=0}^{p}(2j_1+1)(2j_2+1)(2j_3+1)
\left(\left(\bar C_{j_3j_2j_1}^{100}\right)^2+
\bar C_{j_3j_2j_1}^{100}\bar C_{j_1j_2j_3}^{100}\right)\right).
\end{equation}

\vspace{10mm}

\centerline{\bf The case $k=4$
for the integral $I_{(0000)T,t}^{(i_1 i_2 i_3 i_4)}$}

\vspace{5mm}

{\bf 4}.1.\ $i_1,\ldots,i_4$ are pairwise different:
\begin{equation}
\label{prod26}
E^p_4 = (T-t)^4\left(\frac{1}{24}-\frac{1}{256}
\sum_{j_1,\ldots,j_4=0}^{p}\left(\prod_{l=1}^4 (2j_l+1)\right)\cdot
\left(\bar C_{j_4\ldots j_1}^{0000}\right)^2\right).
\end{equation}

\vspace{4mm}

{\bf 4}.2.\ $i_1=i_2=i_3=i_4$:
$$
E^p_4 = (T-t)^4\left(\frac{1}{24}-\frac{1}{256}
\sum_{j_1,\ldots,j_4=0}^{p}\left(\prod_{l=1}^4 (2j_l+1)\right)\cdot
\bar C_{j_4\ldots j_1}^{0000}\left(
\sum\limits_{(j_1,\ldots,j_4)}
\bar C_{j_4\ldots j_1}^{0000}\right)\right)\equiv 0.
$$

\vspace{4mm}

{\bf 4}.3.1.\ $i_1=i_2\ne i_3, i_4;\ i_3\ne i_4:$
\begin{equation}
\label{prod27}
E^p_4 = (T-t)^4\left(\frac{1}{24}-\frac{1}{256}
\sum_{j_1,\ldots,j_4=0}^{p}\left(\prod_{l=1}^4 (2j_l+1)\right)\cdot
\bar C_{j_4\ldots j_1}^{0000}\left(
\sum\limits_{(j_1,j_2)}
\bar C_{j_4\ldots j_1}^{0000}\right)\right).
\end{equation}

\vspace{4mm}

{\bf 4}.3.2.\ $i_1=i_3\ne i_2, i_4;\ i_2\ne i_4:$
\begin{equation}
\label{prod28}
E^p_4 = (T-t)^4\left(\frac{1}{24}-\frac{1}{256}
\sum_{j_1,\ldots,j_4=0}^{p}\left(\prod_{l=1}^4 (2j_l+1)\right)\cdot
\bar C_{j_4\ldots j_1}^{0000}\left(
\sum\limits_{(j_1,j_3)}
\bar C_{j_4\ldots j_1}^{0000}\right)\right).
\end{equation}

\vspace{4mm}

{\bf 4}.3.3.\ $i_1=i_4\ne i_2, i_3;\ i_2\ne i_3:$
\begin{equation}
\label{prod29}
E^p_4 = (T-t)^4\left(\frac{1}{24}-\frac{1}{256}
\sum_{j_1,\ldots,j_4=0}^{p}\left(\prod_{l=1}^4 (2j_l+1)\right)\cdot
\bar C_{j_4\ldots j_1}^{0000}\left(
\sum\limits_{(j_1,j_4)}
\bar C_{j_4\ldots j_1}^{0000}\right)\right).
\end{equation}

\vspace{4mm}

{\bf 4}.3.4.\ $i_2=i_3\ne i_1, i_4;\ i_1\ne i_4:$
\begin{equation}
\label{prod30}
E^p_4 = (T-t)^4\left(\frac{1}{24}-\frac{1}{256}
\sum_{j_1,\ldots,j_4=0}^{p}\left(\prod_{l=1}^4 (2j_l+1)\right)\cdot
\bar C_{j_4\ldots j_1}^{0000}\left(
\sum\limits_{(j_2,j_3)}
\bar C_{j_4\ldots j_1}^{0000}\right)\right).
\end{equation}

\vspace{4mm}

{\bf 4}.3.5.\ $i_2=i_4\ne i_1, i_3;\ i_1\ne i_3:$
\begin{equation}
\label{prod31}
E^p_4 = (T-t)^4\left(\frac{1}{24}-\frac{1}{256}
\sum_{j_1,\ldots,j_4=0}^{p}\left(\prod_{l=1}^4 (2j_l+1)\right)\cdot
\bar C_{j_4\ldots j_1}^{0000}\left(
\sum\limits_{(j_2,j_4)}
\bar C_{j_4\ldots j_1}^{0000}\right)\right).
\end{equation}

\vspace{4mm}

{\bf 4}.3.6.\ $i_3=i_4\ne i_1, i_2;\ i_1\ne i_2:$
\begin{equation}
\label{prod32}
E^p_4 = (T-t)^4\left(\frac{1}{24}-\frac{1}{256}
\sum_{j_1,\ldots,j_4=0}^{p}\left(\prod_{l=1}^4 (2j_l+1)\right)\cdot
\bar C_{j_4\ldots j_1}^{0000}\left(
\sum\limits_{(j_3,j_4)}
\bar C_{j_4\ldots j_1}^{0000}\right)\right).
\end{equation}

\vspace{4mm}

{\bf 4}.4.1.\ $i_1=i_2=i_3\ne i_4$:
\begin{equation}
\label{prod33}
E^p_4 = (T-t)^4\left(\frac{1}{24}-\frac{1}{256}
\sum_{j_1,\ldots,j_4=0}^{p}\left(\prod_{l=1}^4 (2j_l+1)\right)\cdot
\bar C_{j_4\ldots j_1}^{0000}\left(
\sum\limits_{(j_1,j_2,j_3)}
\bar C_{j_4\ldots j_1}^{0000}\right)\right).
\end{equation}

\vspace{4mm}

{\bf 4}.4.2.\ $i_2=i_3=i_4\ne i_1$:
\begin{equation}
\label{prod34}
E^p_4 = (T-t)^4\left(\frac{1}{24}-\frac{1}{256}
\sum_{j_1,\ldots,j_4=0}^{p}\left(\prod_{l=1}^4 (2j_l+1)\right)\cdot
\bar C_{j_4\ldots j_1}^{0000}\left(
\sum\limits_{(j_2,j_3,j_4)}
\bar C_{j_4\ldots j_1}^{0000}\right)\right).
\end{equation}

\vspace{4mm}

{\bf 4}.4.3.\ $i_1=i_2=i_4\ne i_3$:
\begin{equation}
\label{prod35}
E^p_4 = (T-t)^4\left(\frac{1}{24}-\frac{1}{256}
\sum_{j_1,\ldots,j_4=0}^{p}\left(\prod_{l=1}^4 (2j_l+1)\right)\cdot
\bar C_{j_4\ldots j_1}^{0000}\left(
\sum\limits_{(j_1,j_2,j_4)}
\bar C_{j_4\ldots j_1}^{0000}\right)\right).
\end{equation}

\vspace{4mm}

{\bf 4}.4.4.\ $i_1=i_3=i_4\ne i_2$:
\begin{equation}
\label{prod36}
E^p_4 = (T-t)^4\left(\frac{1}{24}-\frac{1}{256}
\sum_{j_1,\ldots,j_4=0}^{p}\left(\prod_{l=1}^4 (2j_l+1)\right)\cdot
\bar C_{j_4\ldots j_1}^{0000}\left(
\sum\limits_{(j_1,j_3,j_4)}
\bar C_{j_4\ldots j_1}^{0000}\right)\right).
\end{equation}

\vspace{4mm}

{\bf 4}.5.1.\ $i_1=i_2\ne i_3=i_4$:
\begin{equation}
\label{prod37}
E^p_4 = (T-t)^4\left(\frac{1}{24}-\frac{1}{256}
\sum_{j_1,\ldots,j_4=0}^{p}\left(\prod_{l=1}^4 (2j_l+1)\right)\cdot
\bar C_{j_4\ldots j_1}^{0000}\left(\sum\limits_{(j_1,j_2)}\left(
\sum\limits_{(j_3,j_4)}
\bar C_{j_4\ldots j_1}^{0000}\right)\right)\right).
\end{equation}

\vspace{4mm}

{\bf 4}.5.2.\ $i_1=i_3\ne i_2=i_4$:
\begin{equation}
\label{prod38}
E^p_4 = (T-t)^4\left(\frac{1}{24}-\frac{1}{256}
\sum_{j_1,\ldots,j_4=0}^{p}\left(\prod_{l=1}^4 (2j_l+1)\right)\cdot
\bar C_{j_4\ldots j_1}^{0000}\left(\sum\limits_{(j_1,j_3)}\left(
\sum\limits_{(j_2,j_4)}
\bar C_{j_4\ldots j_1}^{0000}\right)\right)\right).
\end{equation}

\vspace{4mm}

{\bf 4}.5.3.\ $i_1=i_4\ne i_2=i_3$:
\begin{equation}
\label{prod39}
E^p_4 = (T-t)^4\left(\frac{1}{24}-\frac{1}{256}
\sum_{j_1,\ldots,j_4=0}^{p}\left(\prod_{l=1}^4 (2j_l+1)\right)\cdot
\bar C_{j_4\ldots j_1}^{0000}\left(\sum\limits_{(j_1,j_4)}\left(
\sum\limits_{(j_2,j_3)}
\bar C_{j_4\ldots j_1}^{0000}\right)\right)\right).
\end{equation}

\vspace{10mm}

\centerline{\bf The case $k=5$
for the integral $I_{(00000)T,t}^{(i_1 i_2 i_3 i_4 i_5)}$}

\vspace{5mm}

{\bf 5}.1.\ $i_1,\ldots,i_5$ are pairwise different:
\begin{equation}
\label{prod40}
E^p_5 = (T-t)^5\left(\frac{1}{120}-\frac{1}{32^2}
\sum_{j_1,\ldots,j_5=0}^{p}\left(\prod_{l=1}^5 (2j_l+1)\right)\cdot
\left(\bar C_{j_5\ldots j_1}^{00000}\right)^2\right).
\end{equation}

\vspace{4mm}

{\bf 5}.2.\ $i_1=i_2=i_3=i_4=i_5$:
$$
E^p_5 = (T-t)^5\left(\frac{1}{120}-\frac{1}{32^2}
\sum_{j_1,\ldots,j_5=0}^{p}\left(\prod_{l=1}^5 (2j_l+1)\right)\cdot
\bar C_{j_5\ldots j_1}^{00000}\left(
\sum\limits_{(j_1,\ldots,j_5)}
\bar C_{j_5\ldots j_1}^{00000}\right)\right)\equiv 0.
$$

\vspace{7mm}

{\bf 5}.3.1.\ $i_1=i_2\ne i_3,i_4,i_5$\ ($i_3,i_4,i_5$ are pairwise different):
\begin{equation}
\label{prod41}
E^p_5 = (T-t)^5\left(\frac{1}{120}-\frac{1}{32^2}
\sum_{j_1,\ldots,j_5=0}^{p}\left(\prod_{l=1}^5 (2j_l+1)\right)\cdot
\bar C_{j_5\ldots j_1}^{00000}\left(
\sum\limits_{(j_1,j_2)}
\bar C_{j_5\ldots j_1}^{00000}\right)\right).
\end{equation}

\vspace{7mm}

{\bf 5}.3.2.\ $i_1=i_3\ne i_2,i_4,i_5$\ ($i_2,i_4,i_5$ are pairwise different):
\begin{equation}
\label{prod42}
E^p_5 = (T-t)^5\left(\frac{1}{120}-\frac{1}{32^2}
\sum_{j_1,\ldots,j_5=0}^{p}\left(\prod_{l=1}^5 (2j_l+1)\right)\cdot
\bar C_{j_5\ldots j_1}^{00000}\left(
\sum\limits_{(j_1,j_3)}
\bar C_{j_5\ldots j_1}^{00000}\right)\right).
\end{equation}

\vspace{7mm}

{\bf 5}.3.3.\ $i_1=i_4\ne i_2,i_3,i_5$\ ($i_2,i_3,i_5$ are pairwise different):
\begin{equation}
\label{prod43}
E^p_5 = (T-t)^5\left(\frac{1}{120}-\frac{1}{32^2}
\sum_{j_1,\ldots,j_5=0}^{p}\left(\prod_{l=1}^5 (2j_l+1)\right)\cdot
\bar C_{j_5\ldots j_1}^{00000}\left(
\sum\limits_{(j_1,j_4)}
\bar C_{j_5\ldots j_1}^{00000}\right)\right).
\end{equation}

\vspace{7mm}

{\bf 5}.3.4.\ $i_1=i_5\ne i_2,i_3,i_4$\ ($i_2,i_3,i_4$  are pairwise different):
\begin{equation}
\label{prod44}
E^p_5 = (T-t)^5\left(\frac{1}{120}-\frac{1}{32^2}
\sum_{j_1,\ldots,j_5=0}^{p}\left(\prod_{l=1}^5 (2j_l+1)\right)\cdot
\bar C_{j_5\ldots j_1}^{00000}\left(
\sum\limits_{(j_1,j_5)}
\bar C_{j_5\ldots j_1}^{00000}\right)\right).
\end{equation}

\vspace{7mm}

{\bf 5}.3.5.\ $i_2=i_3\ne i_1,i_4,i_5$\ ($i_1,i_4,i_5$ are pairwise different):
\begin{equation}
\label{prod45}
E^p_5 = (T-t)^5\left(\frac{1}{120}-\frac{1}{32^2}
\sum_{j_1,\ldots,j_5=0}^{p}\left(\prod_{l=1}^5 (2j_l+1)\right)\cdot
\bar C_{j_5\ldots j_1}^{00000}\left(
\sum\limits_{(j_2,j_3)}
\bar C_{j_5\ldots j_1}^{00000}\right)\right).
\end{equation}

\vspace{7mm}

{\bf 5}.3.6.\ $i_2=i_4\ne i_1,i_3,i_5$\ ($i_1,i_3,i_5$ are pairwise different):
\begin{equation}
\label{prod46}
E^p_5 = (T-t)^5\left(\frac{1}{120}-\frac{1}{32^2}
\sum_{j_1,\ldots,j_5=0}^{p}\left(\prod_{l=1}^5 (2j_l+1)\right)\cdot
\bar C_{j_5\ldots j_1}^{00000}\left(
\sum\limits_{(j_2,j_4)}
\bar C_{j_5\ldots j_1}^{00000}\right)\right).
\end{equation}

\vspace{7mm}

{\bf 5}.3.7.\ $i_2=i_5\ne i_1,i_3,i_4$\ ($i_1,i_3,i_4$ are pairwise different):
\begin{equation}
\label{prod47}
E^p_5 = (T-t)^5\left(\frac{1}{120}-\frac{1}{32^2}
\sum_{j_1,\ldots,j_5=0}^{p}\left(\prod_{l=1}^5 (2j_l+1)\right)\cdot
\bar C_{j_5\ldots j_1}^{00000}\left(
\sum\limits_{(j_2,j_5)}
\bar C_{j_5\ldots j_1}^{00000}\right)\right).
\end{equation}

\vspace{7mm}

{\bf 5}.3.8.\ $i_3=i_4\ne i_1,i_2,i_5$\ ($i_1,i_2,i_5$  are pairwise different):
\begin{equation}
\label{prod48}
E^p_5 = (T-t)^5\left(\frac{1}{120}-\frac{1}{32^2}
\sum_{j_1,\ldots,j_5=0}^{p}\left(\prod_{l=1}^5 (2j_l+1)\right)\cdot
\bar C_{j_5\ldots j_1}^{00000}\left(
\sum\limits_{(j_3,j_4)}
\bar C_{j_5\ldots j_1}^{00000}\right)\right).
\end{equation}

\vspace{7mm}

{\bf 5}.3.9.\ $i_3=i_5\ne i_1,i_2,i_4$\ ($i_1,i_2,i_4$  are pairwise different):
\begin{equation}
\label{prod49}
E^p_5 = (T-t)^5\left(\frac{1}{120}-\frac{1}{32^2}
\sum_{j_1,\ldots,j_5=0}^{p}\left(\prod_{l=1}^5 (2j_l+1)\right)\cdot
\bar C_{j_5\ldots j_1}^{00000}\left(
\sum\limits_{(j_3,j_5)}
\bar C_{j_5\ldots j_1}^{00000}\right)\right).
\end{equation}

\vspace{7mm}

{\bf 5}.3.10.\ $i_4=i_5\ne i_1,i_2,i_3$\ ($i_1,i_2,i_3$  are pairwise different):
\begin{equation}
\label{prod50}
E^p_5 = (T-t)^5\left(\frac{1}{120}-\frac{1}{32^2}
\sum_{j_1,\ldots,j_5=0}^{p}\left(\prod_{l=1}^5 (2j_l+1)\right)\cdot
\bar C_{j_5\ldots j_1}^{00000}\left(
\sum\limits_{(j_4,j_5)}
\bar C_{j_5\ldots j_1}^{00000}\right)\right).
\end{equation}

\vspace{7mm}

{\bf 5}.4.1.\ $i_1=i_2=i_3\ne i_4, i_5$\ $(i_4\ne i_5$):
\begin{equation}
\label{prod51}
E^p_5 = (T-t)^5\left(\frac{1}{120}-\frac{1}{32^2}
\sum_{j_1,\ldots,j_5=0}^{p}\left(\prod_{l=1}^5 (2j_l+1)\right)\cdot
\bar C_{j_5\ldots j_1}^{00000}\left(
\sum\limits_{(j_1,j_2,j_3)}
\bar C_{j_5\ldots j_1}^{00000}\right)\right).
\end{equation}

\vspace{7mm}

{\bf 5}.4.2.\ $i_1=i_2=i_4\ne i_3, i_5$\  $(i_3\ne i_5$):
\begin{equation}
\label{prod52}
E^p_5 = (T-t)^5\left(\frac{1}{120}-\frac{1}{32^2}
\sum_{j_1,\ldots,j_5=0}^{p}\left(\prod_{l=1}^5 (2j_l+1)\right)\cdot
\bar C_{j_5\ldots j_1}^{00000}\left(
\sum\limits_{(j_1,j_2,j_4)}
\bar C_{j_5\ldots j_1}^{00000}\right)\right).
\end{equation}

\vspace{7mm}

{\bf 5}.4.3.\ $i_1=i_2=i_5\ne i_3, i_4$\  $(i_3\ne i_4$):
\begin{equation}
\label{prod53}
E^p_5 = (T-t)^5\left(\frac{1}{120}-\frac{1}{32^2}
\sum_{j_1,\ldots,j_5=0}^{p}\left(\prod_{l=1}^5 (2j_l+1)\right)\cdot
\bar C_{j_5\ldots j_1}^{00000}\left(
\sum\limits_{(j_1,j_2,j_5)}
\bar C_{j_5\ldots j_1}^{00000}\right)\right).
\end{equation}

\vspace{7mm}

{\bf 5}.4.4.\ $i_2=i_3=i_4\ne i_1, i_5$\  $(i_1\ne i_5$):
\begin{equation}
\label{prod54}
E^p_5 = (T-t)^5\left(\frac{1}{120}-\frac{1}{32^2}
\sum_{j_1,\ldots,j_5=0}^{p}\left(\prod_{l=1}^5 (2j_l+1)\right)\cdot
\bar C_{j_5\ldots j_1}^{00000}\left(
\sum\limits_{(j_2,j_3,j_4)}
\bar C_{j_5\ldots j_1}^{00000}\right)\right).
\end{equation}

\vspace{7mm}
{\bf 5}.4.5.\ $i_2=i_3=i_5\ne i_1, i_4$\  $(i_1\ne i_4$):
\begin{equation}
\label{prod55}
E^p_5 = (T-t)^5\left(\frac{1}{120}-\frac{1}{32^2}
\sum_{j_1,\ldots,j_5=0}^{p}\left(\prod_{l=1}^5 (2j_l+1)\right)\cdot
\bar C_{j_5\ldots j_1}^{00000}\left(
\sum\limits_{(j_2,j_3,j_5)}
\bar C_{j_5\ldots j_1}^{00000}\right)\right).
\end{equation}

\vspace{7mm}

{\bf 5}.4.6.\ $i_2=i_4=i_5\ne i_1, i_3$\  $(i_1\ne i_3$):
\begin{equation}
\label{prod56}
E^p_5 = (T-t)^5\left(\frac{1}{120}-\frac{1}{32^2}
\sum_{j_1,\ldots,j_5=0}^{p}\left(\prod_{l=1}^5 (2j_l+1)\right)\cdot
\bar C_{j_5\ldots j_1}^{00000}\left(
\sum\limits_{(j_2,j_4,j_5)}
\bar C_{j_5\ldots j_1}^{00000}\right)\right).
\end{equation}

\vspace{7mm}

{\bf 5}.4.7.\ $i_3=i_4=i_5\ne i_1, i_2$\  $(i_1\ne i_2$):
\begin{equation}
\label{prod57}
E^p_5 = (T-t)^5\left(\frac{1}{120}-\frac{1}{32^2}
\sum_{j_1,\ldots,j_5=0}^{p}\left(\prod_{l=1}^5 (2j_l+1)\right)\cdot
\bar C_{j_5\ldots j_1}^{00000}\left(
\sum\limits_{(j_3,j_4,j_5)}
\bar C_{j_5\ldots j_1}^{00000}\right)\right).
\end{equation}

\vspace{7mm}

{\bf 5}.4.8.\ $i_1=i_3=i_5\ne i_2, i_4$\  $(i_2\ne i_4$):
\begin{equation}
\label{prod58}
E^p_5 = (T-t)^5\left(\frac{1}{120}-\frac{1}{32^2}
\sum_{j_1,\ldots,j_5=0}^{p}\left(\prod_{l=1}^5 (2j_l+1)\right)\cdot
\bar C_{j_5\ldots j_1}^{00000}\left(
\sum\limits_{(j_1,j_3,j_5)}
\bar C_{j_5\ldots j_1}^{00000}\right)\right).
\end{equation}

\vspace{7mm}

{\bf 5}.4.9.\ $i_1=i_3=i_4\ne i_2, i_5$\  $(i_2\ne i_5$):
\begin{equation}
\label{prod59}
E^p_5 = (T-t)^5\left(\frac{1}{120}-\frac{1}{32^2}
\sum_{j_1,\ldots,j_5=0}^{p}\left(\prod_{l=1}^5 (2j_l+1)\right)\cdot
\bar C_{j_5\ldots j_1}^{00000}\left(
\sum\limits_{(j_1,j_3,j_4)}
\bar C_{j_5\ldots j_1}^{00000}\right)\right).
\end{equation}

\vspace{7mm}

{\bf 5}.4.10.\ $i_1=i_4=i_5\ne i_2, i_3$\  $(i_2\ne i_3$):
\begin{equation}
\label{prod60}
E^p_5 = (T-t)^5\left(\frac{1}{120}-\frac{1}{32^2}
\sum_{j_1,\ldots,j_5=0}^{p}\left(\prod_{l=1}^5 (2j_l+1)\right)\cdot
\bar C_{j_5\ldots j_1}^{00000}\left(
\sum\limits_{(j_1,j_4,j_5)}
\bar C_{j_5\ldots j_1}^{00000}\right)\right).
\end{equation}

\vspace{7mm}

{\bf 5}.5.1.\ $i_1=i_2=i_3=i_4\ne i_5$:
\begin{equation}
\label{prod61}
E^p_5 = (T-t)^5\left(\frac{1}{120}-\frac{1}{32^2}
\sum_{j_1,\ldots,j_5=0}^{p}\left(\prod_{l=1}^5 (2j_l+1)\right)\cdot
\bar C_{j_5\ldots j_1}^{00000}\left(
\sum\limits_{(j_1,j_2,j_3,j_4)}
\bar C_{j_5\ldots j_1}^{00000}\right)\right).
\end{equation}

\vspace{7mm}

{\bf 5}.5.2.\ $i_1=i_2=i_3=i_5\ne i_4$:
\begin{equation}
\label{prod62}
E^p_5 = (T-t)^5\left(\frac{1}{120}-\frac{1}{32^2}
\sum_{j_1,\ldots,j_5=0}^{p}\left(\prod_{l=1}^5 (2j_l+1)\right)\cdot
\bar C_{j_5\ldots j_1}^{00000}\left(
\sum\limits_{(j_1,j_2,j_3,j_5)}
\bar C_{j_5\ldots j_1}^{00000}\right)\right).
\end{equation}

\vspace{7mm}

{\bf 5}.5.3.\ $i_1=i_2=i_4=i_5\ne i_3$:
\begin{equation}
\label{prod63}
E^p_5 = (T-t)^5\left(\frac{1}{120}-\frac{1}{32^2}
\sum_{j_1,\ldots,j_5=0}^{p}\left(\prod_{l=1}^5 (2j_l+1)\right)\cdot
\bar C_{j_5\ldots j_1}^{00000}\left(
\sum\limits_{(j_1,j_2,j_4,j_5)}
\bar C_{j_5\ldots j_1}^{00000}\right)\right).
\end{equation}

\vspace{7mm}

{\bf 5}.5.4.\ $i_1=i_3=i_4=i_5\ne i_2$:
\begin{equation}
\label{prod64}
E^p_5 = (T-t)^5\left(\frac{1}{120}-\frac{1}{32^2}
\sum_{j_1,\ldots,j_5=0}^{p}\left(\prod_{l=1}^5 (2j_l+1)\right)\cdot
\bar C_{j_5\ldots j_1}^{00000}\left(
\sum\limits_{(j_1,j_3,j_4,j_5)}
\bar C_{j_5\ldots j_1}^{00000}\right)\right).
\end{equation}

\vspace{7mm}

{\bf 5}.5.5.\ $i_2=i_3=i_4=i_5\ne i_1$:
\begin{equation}
\label{prod65}
E^p_5 = (T-t)^5\left(\frac{1}{120}-\frac{1}{32^2}
\sum_{j_1,\ldots,j_5=0}^{p}\left(\prod_{l=1}^5 (2j_l+1)\right)\cdot
\bar C_{j_5\ldots j_1}^{00000}\left(
\sum\limits_{(j_2,j_3,j_4,j_5)}
\bar C_{j_5\ldots j_1}^{00000}\right)\right).
\end{equation}

\vspace{7mm}

{\bf 5}.6.1.\ $i_5\ne i_1=i_2\ne i_3=i_4\ne i_5$:
\begin{equation}
\label{prod66}
E^p_5 = (T-t)^5\left(\frac{1}{120}-\frac{1}{32^2}
\sum_{j_1,\ldots,j_5=0}^{p}\left(\prod_{l=1}^5 (2j_l+1)\right)\cdot
\bar C_{j_5\ldots j_1}^{00000}\left(\sum\limits_{(j_1,j_2)}\left(
\sum\limits_{(j_3,j_4)}
\bar C_{j_5\ldots j_1}^{00000}\right)\right)\right).
\end{equation}

\vspace{7mm}

{\bf 5}.6.2.\ $i_5\ne i_1=i_3\ne i_2=i_4\ne i_5$:
\begin{equation}
\label{prod67}
E^p_5 = (T-t)^5\left(\frac{1}{120}-\frac{1}{32^2}
\sum_{j_1,\ldots,j_5=0}^{p}\left(\prod_{l=1}^5 (2j_l+1)\right)\cdot
\bar C_{j_5\ldots j_1}^{00000}\left(\sum\limits_{(j_1,j_3)}\left(
\sum\limits_{(j_2,j_4)}
\bar C_{j_5\ldots j_1}^{00000}\right)\right)\right).
\end{equation}

\vspace{7mm}

{\bf 5}.6.3.\ $i_5\ne i_1=i_4\ne i_2=i_3\ne i_5$:
\begin{equation}
\label{prod68}
E^p_5 = (T-t)^5\left(\frac{1}{120}-\frac{1}{32^2}
\sum_{j_1,\ldots,j_5=0}^{p}\left(\prod_{l=1}^5 (2j_l+1)\right)\cdot
\bar C_{j_5\ldots j_1}^{00000}\left(\sum\limits_{(j_1,j_4)}\left(
\sum\limits_{(j_2,j_3)}
\bar C_{j_5\ldots j_1}^{00000}\right)\right)\right).
\end{equation}

\vspace{7mm}

{\bf 5}.6.4.\ $i_4\ne i_1=i_2\ne i_3=i_5\ne i_4$:
\begin{equation}
\label{prod69}
E^p_5 = (T-t)^5\left(\frac{1}{120}-\frac{1}{32^2}
\sum_{j_1,\ldots,j_5=0}^{p}\left(\prod_{l=1}^5 (2j_l+1)\right)\cdot
\bar C_{j_5\ldots j_1}^{00000}\left(\sum\limits_{(j_1,j_2)}\left(
\sum\limits_{(j_3,j_5)}
\bar C_{j_5\ldots j_1}^{00000}\right)\right)\right).
\end{equation}

\vspace{7mm}

{\bf 5}.6.5.\ $i_4\ne i_1=i_5\ne i_2=i_3\ne i_4$:
\begin{equation}
\label{prod70}
E^p_5 = (T-t)^5\left(\frac{1}{120}-\frac{1}{32^2}
\sum_{j_1,\ldots,j_5=0}^{p}\left(\prod_{l=1}^5 (2j_l+1)\right)\cdot
\bar C_{j_5\ldots j_1}^{00000}\left(\sum\limits_{(j_1,j_5)}\left(
\sum\limits_{(j_2,j_3)}
\bar C_{j_5\ldots j_1}^{00000}\right)\right)\right).
\end{equation}

\vspace{7mm}

{\bf 5}.6.6.\ $i_4\ne i_2=i_5\ne i_1=i_3\ne i_4$:
\begin{equation}
\label{prod71}
E^p_5 = (T-t)^5\left(\frac{1}{120}-\frac{1}{32^2}
\sum_{j_1,\ldots,j_5=0}^{p}\left(\prod_{l=1}^5 (2j_l+1)\right)\cdot
\bar C_{j_5\ldots j_1}^{00000}\left(\sum\limits_{(j_2,j_5)}\left(
\sum\limits_{(j_1,j_3)}
\bar C_{j_5\ldots j_1}^{00000}\right)\right)\right).
\end{equation}

\vspace{7mm}

{\bf 5}.6.7.\ $i_3\ne i_2=i_5\ne i_1=i_4\ne i_3$:
\begin{equation}
\label{prod72}
E^p_5 = (T-t)^5\left(\frac{1}{120}-\frac{1}{32^2}
\sum_{j_1,\ldots,j_5=0}^{p}\left(\prod_{l=1}^5 (2j_l+1)\right)\cdot
\bar C_{j_5\ldots j_1}^{00000}\left(\sum\limits_{(j_2,j_5)}\left(
\sum\limits_{(j_1,j_4)}
\bar C_{j_5\ldots j_1}^{00000}\right)\right)\right).
\end{equation}

\vspace{7mm}

{\bf 5}.6.8.\ $i_3\ne i_1=i_2\ne i_4=i_5\ne i_3$:
\begin{equation}
\label{prod73}
E^p_5 = (T-t)^5\left(\frac{1}{120}-\frac{1}{32^2}
\sum_{j_1,\ldots,j_5=0}^{p}\left(\prod_{l=1}^5 (2j_l+1)\right)\cdot
\bar C_{j_5\ldots j_1}^{00000}\left(\sum\limits_{(j_1,j_2)}\left(
\sum\limits_{(j_4,j_5)}
\bar C_{j_5\ldots j_1}^{00000}\right)\right)\right).
\end{equation}

\vspace{7mm}

{\bf 5}.6.9.\ $i_3\ne i_2=i_4\ne i_1=i_5\ne i_3$:
\begin{equation}
\label{prod74}
E^p_5 = (T-t)^5\left(\frac{1}{120}-\frac{1}{32^2}
\sum_{j_1,\ldots,j_5=0}^{p}\left(\prod_{l=1}^5 (2j_l+1)\right)\cdot
\bar C_{j_5\ldots j_1}^{00000}\left(\sum\limits_{(j_2,j_4)}\left(
\sum\limits_{(j_1,j_5)}
\bar C_{j_5\ldots j_1}^{00000}\right)\right)\right).
\end{equation}

\vspace{7mm}

{\bf 5}.6.10.\ $i_2\ne i_1=i_4\ne i_3=i_5\ne i_2$:
\begin{equation}
\label{prod75}
E^p_5 = (T-t)^5\left(\frac{1}{120}-\frac{1}{32^2}
\sum_{j_1,\ldots,j_5=0}^{p}\left(\prod_{l=1}^5 (2j_l+1)\right)\cdot
\bar C_{j_5\ldots j_1}^{00000}\left(\sum\limits_{(j_1,j_4)}\left(
\sum\limits_{(j_3,j_5)}
\bar C_{j_5\ldots j_1}^{00000}\right)\right)\right).
\end{equation}

\vspace{7mm}

{\bf 5}.6.11.\ $i_2\ne i_1=i_3\ne i_4=i_5\ne i_2$:
\begin{equation}
\label{prod76}
E^p_5 = (T-t)^5\left(\frac{1}{120}-\frac{1}{32^2}
\sum_{j_1,\ldots,j_5=0}^{p}\left(\prod_{l=1}^5 (2j_l+1)\right)\cdot
\bar C_{j_5\ldots j_1}^{00000}\left(\sum\limits_{(j_1,j_3)}\left(
\sum\limits_{(j_4,j_5)}
\bar C_{j_5\ldots j_1}^{00000}\right)\right)\right).
\end{equation}

\vspace{7mm}

{\bf 5}.6.12.\ $i_2\ne i_1=i_5\ne i_3=i_4\ne i_2$:
\begin{equation}
\label{prod77}
E^p_5 = (T-t)^5\left(\frac{1}{120}-\frac{1}{32^2}
\sum_{j_1,\ldots,j_5=0}^{p}\left(\prod_{l=1}^5 (2j_l+1)\right)\cdot
\bar C_{j_5\ldots j_1}^{00000}\left(\sum\limits_{(j_1,j_5)}\left(
\sum\limits_{(j_3,j_4)}
\bar C_{j_5\ldots j_1}^{00000}\right)\right)\right).
\end{equation}

\vspace{7mm}

{\bf 5}.6.13.\ $i_1\ne i_2=i_3\ne i_4=i_5\ne i_1$:
\begin{equation}
\label{prod78}
E^p_5 = (T-t)^5\left(\frac{1}{120}-\frac{1}{32^2}
\sum_{j_1,\ldots,j_5=0}^{p}\left(\prod_{l=1}^5 (2j_l+1)\right)\cdot
\bar C_{j_5\ldots j_1}^{00000}\left(\sum\limits_{(j_2,j_3)}\left(
\sum\limits_{(j_4,j_5)}
\bar C_{j_5\ldots j_1}^{00000}\right)\right)\right).
\end{equation}

\vspace{7mm}

{\bf 5}.6.14.\ $i_1\ne i_2=i_4\ne i_3=i_5\ne i_1$:
\begin{equation}
\label{prod79}
E^p_5 = (T-t)^5\left(\frac{1}{120}-\frac{1}{32^2}
\sum_{j_1,\ldots,j_5=0}^{p}\left(\prod_{l=1}^5 (2j_l+1)\right)\cdot
\bar C_{j_5\ldots j_1}^{00000}\left(\sum\limits_{(j_2,j_4)}\left(
\sum\limits_{(j_3,j_5)}
\bar C_{j_5\ldots j_1}^{00000}\right)\right)\right).
\end{equation}

\vspace{7mm}

{\bf 5}.6.15.\ $i_1\ne i_2=i_5\ne i_3=i_4\ne i_1$:
\begin{equation}
\label{prod80}
E^p_5 = (T-t)^5\left(\frac{1}{120}-\frac{1}{32^2}
\sum_{j_1,\ldots,j_5=0}^{p}\left(\prod_{l=1}^5 (2j_l+1)\right)\cdot
\bar C_{j_5\ldots j_1}^{00000}\left(\sum\limits_{(j_2,j_5)}\left(
\sum\limits_{(j_3,j_4)}
\bar C_{j_5\ldots j_1}^{00000}\right)\right)\right).
\end{equation}

\vspace{7mm}

{\bf 5}.7.1.\ $i_1=i_2=i_3\ne i_4=i_5$:
\begin{equation}
\label{prod81}
E^p_5 = (T-t)^5\left(\frac{1}{120}-\frac{1}{32^2}
\sum_{j_1,\ldots,j_5=0}^{p}\left(\prod_{l=1}^5 (2j_l+1)\right)\cdot
\bar C_{j_5\ldots j_1}^{00000}\left(\sum\limits_{(j_4,j_5)}\left(
\sum\limits_{(j_1,j_2,j_3)}
\bar C_{j_5\ldots j_1}^{00000}\right)\right)\right).
\end{equation}

\vspace{7mm}

{\bf 5}.7.2.\ $i_1=i_2=i_4\ne i_3=i_5$:
\begin{equation}
\label{prod82}
E^p_5 = (T-t)^5\left(\frac{1}{120}-\frac{1}{32^2}
\sum_{j_1,\ldots,j_5=0}^{p}\left(\prod_{l=1}^5 (2j_l+1)\right)\cdot
\bar C_{j_5\ldots j_1}^{00000}\left(\sum\limits_{(j_3,j_5)}\left(
\sum\limits_{(j_1,j_2,j_4)}
\bar C_{j_5\ldots j_1}^{00000}\right)\right)\right).
\end{equation}

\vspace{7mm}

{\bf 5}.7.3.\ $i_1=i_2=i_5\ne i_3=i_4$:
\begin{equation}
\label{prod83}
E^p_5 = (T-t)^5\left(\frac{1}{120}-\frac{1}{32^2}
\sum_{j_1,\ldots,j_5=0}^{p}\left(\prod_{l=1}^5 (2j_l+1)\right)\cdot
\bar C_{j_5\ldots j_1}^{00000}\left(\sum\limits_{(j_3,j_4)}\left(
\sum\limits_{(j_1,j_2,j_5)}
\bar C_{j_5\ldots j_1}^{00000}\right)\right)\right).
\end{equation}

\vspace{8mm}

{\bf 5}.7.4.\ $i_2=i_3=i_4\ne i_1=i_5$:
\begin{equation}
\label{prod84}
E^p_5 = (T-t)^5\left(\frac{1}{120}-\frac{1}{32^2}
\sum_{j_1,\ldots,j_5=0}^{p}\left(\prod_{l=1}^5 (2j_l+1)\right)\cdot
\bar C_{j_5\ldots j_1}^{00000}\left(\sum\limits_{(j_1,j_5)}\left(
\sum\limits_{(j_2,j_3,j_4)}
\bar C_{j_5\ldots j_1}^{00000}\right)\right)\right).
\end{equation}

\vspace{8mm}

{\bf 5}.7.5.\ $i_2=i_3=i_5\ne i_1=i_4$:
\begin{equation}
\label{prod85}
E^p_5 = (T-t)^5\left(\frac{1}{120}-\frac{1}{32^2}
\sum_{j_1,\ldots,j_5=0}^{p}\left(\prod_{l=1}^5 (2j_l+1)\right)\cdot
\bar C_{j_5\ldots j_1}^{00000}\left(\sum\limits_{(j_1,j_4)}\left(
\sum\limits_{(j_2,j_3,j_5)}
\bar C_{j_5\ldots j_1}^{00000}\right)\right)\right).
\end{equation}

\vspace{8mm}

{\bf 5}.7.6.\ $i_2=i_4=i_5\ne i_1=i_3$:
\begin{equation}
\label{prod86}
E^p_5 = (T-t)^5\left(\frac{1}{120}-\frac{1}{32^2}
\sum_{j_1,\ldots,j_5=0}^{p}\left(\prod_{l=1}^5 (2j_l+1)\right)\cdot
\bar C_{j_5\ldots j_1}^{00000}\left(\sum\limits_{(j_1,j_3)}\left(
\sum\limits_{(j_2,j_4,j_5)}
\bar C_{j_5\ldots j_1}^{00000}\right)\right)\right).
\end{equation}

\vspace{8mm}

{\bf 5}.7.7.\ $i_3=i_4=i_5\ne i_1=i_2$:
\begin{equation}
\label{prod87}
E^p_5 = (T-t)^5\left(\frac{1}{120}-\frac{1}{32^2}
\sum_{j_1,\ldots,j_5=0}^{p}\left(\prod_{l=1}^5 (2j_l+1)\right)\cdot
\bar C_{j_5\ldots j_1}^{00000}\left(\sum\limits_{(j_1,j_2)}\left(
\sum\limits_{(j_3,j_4,j_5)}
\bar C_{j_5\ldots j_1}^{00000}\right)\right)\right).
\end{equation}

\vspace{8mm}

{\bf 5}.7.8.\ $i_1=i_3=i_5\ne i_2=i_4$:
\begin{equation}
\label{prod88}
E^p_5 = (T-t)^5\left(\frac{1}{120}-\frac{1}{32^2}
\sum_{j_1,\ldots,j_5=0}^{p}\left(\prod_{l=1}^5 (2j_l+1)\right)\cdot
\bar C_{j_5\ldots j_1}^{00000}\left(\sum\limits_{(j_2,j_4)}\left(
\sum\limits_{(j_1,j_3,j_5)}
\bar C_{j_5\ldots j_1}^{00000}\right)\right)\right).
\end{equation}

\vspace{8mm}

{\bf 5}.7.9.\ $i_1=i_3=i_4\ne i_2=i_5$:
\begin{equation}
\label{prod89}
E^p_5 = (T-t)^5\left(\frac{1}{120}-\frac{1}{32^2}
\sum_{j_1,\ldots,j_5=0}^{p}\left(\prod_{l=1}^5 (2j_l+1)\right)\cdot
\bar C_{j_5\ldots j_1}^{00000}\left(\sum\limits_{(j_2,j_5)}\left(
\sum\limits_{(j_1,j_3,j_4)}
\bar C_{j_5\ldots j_1}^{00000}\right)\right)\right).
\end{equation}

\vspace{8mm}

{\bf 5}.7.10.\ $i_1=i_4=i_5\ne i_2=i_3$:
\begin{equation}
\label{prod90}
E^p_5 = (T-t)^5\left(\frac{1}{120}-\frac{1}{32^2}
\sum_{j_1,\ldots,j_5=0}^{p}\left(\prod_{l=1}^5 (2j_l+1)\right)\cdot
\bar C_{j_5\ldots j_1}^{00000}\left(\sum\limits_{(j_2,j_3)}\left(
\sum\limits_{(j_1,j_4,j_5)}
\bar C_{j_5\ldots j_1}^{00000}\right)\right)\right).
\end{equation}

\vspace{5mm}

Denote $q(\alpha)$ the numbers $p$
from the formulas (\ref{prod1})-(\ref{prod90}), 
where $\alpha$ are the numbers of the cases corresponding 
to the formulas (\ref{prod1})-(\ref{prod90}). 
For example, $q({\bf 2}.2.b)$ is the number $p$
from the formula (\ref{prod4}), $q({\bf 5}.7.10)$
is the number $p$
from the formula (\ref{prod90}), etc. 

Let                                              

\vspace{-2mm}
\begin{equation}
\label{tred1}
E^p_2\le (T-t)^4,\ \ \ E^p_3\le (T-t)^4
\end{equation}

\vspace{6mm}
\noindent
where $E^p_2$ is defined by (\ref{prod1}) and
$E^p_3$ is defined by (\ref{prod7})-(\ref{prod10}). 

Let

\vspace{-2mm}
\begin{equation}
\label{tred2}
E^p_2\le (T-t)^5,\ \ \ E^p_3\le (T-t)^5,\ \ \ E^p_4\le (T-t)^5
\end{equation}

\vspace{6mm}
\noindent
where $E^p_2$ is defined by (\ref{prod1})-(\ref{prod6}),
$E^p_3$ is defined by (\ref{prod7})-(\ref{prod10}),
and $E^p_4$ is defined by (\ref{prod26})-(\ref{prod39}).

Let

\vspace{-2mm}
\begin{equation}
\label{tred3}
E^p_2\le (T-t)^6,\ \ \ E^p_3\le (T-t)^6,\ \ \ E^p_4\le (T-t)^6,\ \ \
E^p_5\le (T-t)^6,
\end{equation}

\vspace{6mm}
\noindent
where $E^p_2$ is defined by (\ref{prod1})-(\ref{prod6}),
$E^p_3$ is defined by (\ref{prod7})-(\ref{prod25}),
$E^p_4$ is defined by (\ref{prod26})-(\ref{prod39}),
and $E^p_5$ is defined by (\ref{prod40})-(\ref{prod90}).

Note that the conditions (\ref{tred1})-(\ref{tred3})
are particular cases of (\ref{uslov}) for $r=3, 4,$ and $5.$

Let us show by numerical experiments (see Tables 1--13) that in most situations
the following inequalities are fullfilled (under conditions
(\ref{tred1})-(\ref{tred3}))

\begin{equation}
\label{som1}
q({\bf 2}.1.b)\ge q({\bf 2}.2.b),\ \ \ 
q({\bf 2}.1.c)\ge q({\bf 2}.2.c),
\end{equation}

\begin{equation}
\label{som2}
q({\bf 3}.1.a)\ge q({\bf 3}.3.1.a),\ 
q({\bf 3}.3.2.a),\ 
q({\bf 3}.3.3.a),
\end{equation}

\begin{equation}
\label{som3}
q({\bf 3}.1.b)\ge q({\bf 3}.2.b),\ q({\bf 3}.3.1.b),\ 
q({\bf 3}.3.2.b),\ 
q({\bf 3}.3.3.b),
\end{equation}

\begin{equation}
\label{som4}
q({\bf 3}.1.c)\ge q({\bf 3}.2.c),\ q({\bf 3}.3.1.c),\ 
q({\bf 3}.3.2.c),\ 
q({\bf 3}.3.3.c),
\end{equation}

\begin{equation}
\label{som5}
q({\bf 3}.1.d)\ge q({\bf 3}.2.d),\ q({\bf 3}.3.1.d),\ 
q({\bf 3}.3.2.d),\ 
q({\bf 3}.3.3.d),
\end{equation}

\begin{equation}
\label{som6}
q({\bf 4}.1)\ge q({\bf 4}.3.1),\ \ldots,\ q({\bf 4}.3.6),\ 
q({\bf 4}.4.1),\ \ldots,\ q({\bf 4}.4.4),\ 
q({\bf 4}.5.1),\ \ldots,\ q({\bf 4}.5.3),\ 
\end{equation}

\begin{equation}
\label{som7}
q({\bf 5}.1)\ge q({\bf 5}.3.1),\ \ldots,\ q({\bf 5}.3.10),\ 
q({\bf 5}.4.1),\ \ldots,\ q({\bf 5}.4.10),\ 
q({\bf 5}.5.1),\ \ldots,\ q({\bf 5}.5.5),\ 
\end{equation}

\begin{equation}
\label{som8}
q({\bf 5}.1)\ge q({\bf 5}.6.1),\ \ldots,\ q({\bf 5}.6.15),\ 
q({\bf 5}.7.1),\ \ldots,\ q({\bf 5}.7.10),
\end{equation}

\vspace{5mm}
\noindent
where all numbers in the inequalities (\ref{som1})-(\ref{som8}) are minimal
natural numbers satisfying the conditions 
(\ref{tred1})-(\ref{tred3}).

\begin{table}
\centering
\caption{Stochastic integral $I_{(000)T,t}^{(i_1i_2i_3)}$. The condition
(\ref{tred1}).}
\label{tab:1}      
\begin{tabular}{p{1.7cm}p{1.7cm}p{1.7cm}p{1.7cm}p{1.7cm}p{1.7cm}p{1.7cm}}
\hline\noalign{\smallskip}
$T-t$&$0.011$&$0.008$&$0.0045$&$0.0035$&$0.0027$&$0.0025$\\
\noalign{\smallskip}\hline\noalign{\smallskip}
$q({\bf 3}.1.a)$&$\fbox{12}$&$\fbox{16}$&$\fbox{28}$&$\fbox{36}$&$\fbox{47}$&$\fbox{50}$\\
$q({\bf 3}.3.1.a)$&$6$&$8$&$14$&$18$&$23$&$25$\\
$q({\bf 3}.3.2.a)$&$6$&$8$&$14$&$18$&$23$&$25$\\
$q({\bf 3}.3.3.a)$&$12$&$16$&$28$&$36$&$47$&$\fbox{51}$\\
\noalign{\smallskip}\hline\noalign{\smallskip}
\end{tabular}
\vspace{6mm}
\end{table}

\begin{table}
\centering
\caption{Stochastic integrals $I_{(01)T,t}^{(i_1i_2)}$, 
$I_{(10)T,t}^{(i_1i_2)}$. The condition (\ref{tred2}).}
\label{tab:2}      
\begin{tabular}{p{3.0cm}p{3.0cm}p{3.0cm}p{3.0cm}}
\hline\noalign{\smallskip}
$T-t$&0.010&0.005&0.0025\\
\noalign{\smallskip}\hline\noalign{\smallskip}
$q({\bf 2}.1.b)$&$\fbox{4}$&$\fbox{8}$&$\fbox{16}$\\
$q({\bf 2}.2.b)$&1&1&1\\
$q({\bf 2}.1.c)$&$\fbox{4}$&$\fbox{8}$&$\fbox{16}$\\
$q({\bf 2}.2.c)$&1&1&1\\
\noalign{\smallskip}\hline\noalign{\smallskip}
\end{tabular}
\vspace{6mm}
\end{table}

\begin{table}
\centering
\caption{$T-t=0.01.$ The condition (\ref{tred3}).}
\label{tab:3}      
\begin{tabular}{p{3.0cm}p{3.0cm}p{3.0cm}}
\hline\noalign{\smallskip}
$I_{(001)T,t}^{(i_1i_2i_3)}$&$I_{(010)T,t}^{(i_1i_2i_3)}$&$I_{(100)T,t}^{(i_1i_2i_3)}$\\
\noalign{\smallskip}\hline\noalign{\smallskip}
$q({\bf 3}.1.b)=\fbox{6}$&$q({\bf 3}.1.c)=\fbox{4}$&$q({\bf 3}.1.d)=\fbox{2}$\\
$q({\bf 3}.2.b)=0$&$q({\bf 3}.2.c)=0$&$q({\bf 3}.2.d)=0$\\
$q({\bf 3}.3.1.b)=3$&$q({\bf 3}.3.1.c)=3$&$q({\bf 3}.3.1.d)=1$\\
$q({\bf 3}.3.2.b)=3$&$q({\bf 3}.3.2.c)=1$&$q({\bf 3}.3.2.d)=1$\\
$q({\bf 3}.3.3.b)=6$&$q({\bf 3}.3.3.c)=4$&$q({\bf 3}.3.3.d)=2$\\
\noalign{\smallskip}\hline\noalign{\smallskip}
\end{tabular}
\vspace{9mm}
\end{table}

\begin{table}
\centering
\caption{Stochastic integral $I_{(0000)T,t}^{(i_1i_2i_3i_4)}$. The condition (\ref{tred2}).}
\label{tab:4}      
\begin{tabular}{p{1.7cm}p{1.7cm}p{1.7cm}p{1.7cm}p{1.7cm}p{1.7cm}}
\hline\noalign{\smallskip}
$T-t$&$0.011$&$0.008$&$0.0045$&$0.0042$&$0.0040$\\
\noalign{\smallskip}\hline\noalign{\smallskip}
$q({\bf 4}.1)$&$\fbox{6}$&$\fbox{8}$&$\fbox{14}$&$\fbox{15}$&$\fbox{16}$\\
$q({\bf 4}.3.1)$&4&5&10&11&11\\
$q({\bf 4}.3.2)$&6&8&14&15&16\\
$q({\bf 4}.3.3)$&6&8&14&15&16\\
$q({\bf 4}.3.4)$&3&5&9&9&10\\
$q({\bf 4}.3.5)$&6&8&14&15&16\\
$q({\bf 4}.3.6)$&4&5&10&11&11\\
$q({\bf 4}.4.1)$&2&3&4&5&5\\
$q({\bf 4}.4.2)$&2&3&4&5&5\\
$q({\bf 4}.4.3)$&4&6&10&11&11\\
$q({\bf 4}.4.4)$&4&6&10&11&11\\
$q({\bf 4}.5.1)$&2&3&5&6&6\\
$q({\bf 4}.5.2)$&6&8&14&15&16\\
$q({\bf 4}.5.3)$&3&5&9&9&10\\
\noalign{\smallskip}\hline\noalign{\smallskip}
\end{tabular}
\vspace{9mm}
\end{table}

\begin{table}
\centering
\caption{$T-t=0.011.$ Stochastic integral 
$I_{(00000)T,t}^{(i_1i_2i_3i_4i_5)}.$ The condition (\ref{tred3}).}
\label{tab:5}      
\begin{tabular}{p{1.9cm}p{2.2cm}p{2.2cm}p{2.2cm}p{2.2cm}p{2.2cm}}
\hline\noalign{\smallskip}
$q({\bf 5}.1)=\fbox{0}$&$q({\bf 5}.3.1)=0$&$q({\bf 5}.4.1)=0$&$q({\bf 5}.5.1)=0$&$q({\bf 5}.6.1)=0$&$q({\bf 5}.7.1)=0$\\
                &$q({\bf 5}.3.2)=0$&$q({\bf 5}.4.2)=0$&$q({\bf 5}.5.2)=0$&$q({\bf 5}.6.2)=0$&$q({\bf 5}.7.2)=0$\\
                &$q({\bf 5}.3.3)=0$&$q({\bf 5}.4.3)=0$&$q({\bf 5}.5.3)=0$&$q({\bf 5}.6.3)=0$&$q({\bf 5}.7.3)=0$\\
                &$q({\bf 5}.3.4)=0$&$q({\bf 5}.4.4)=0$&$q({\bf 5}.5.4)=0$&$q({\bf 5}.6.4)=0$&$q({\bf 5}.7.4)=0$\\
                &$q({\bf 5}.3.5)=0$&$q({\bf 5}.4.5)=0$&$q({\bf 5}.5.5)=0$&$q({\bf 5}.6.5)=0$&$q({\bf 5}.7.5)=0$\\
                &$q({\bf 5}.3.6)=0$&$q({\bf 5}.4.6)=0$&                  &$q({\bf 5}.6.6)=0$&$q({\bf 5}.7.6)=0$\\
                &$q({\bf 5}.3.7)=0$&$q({\bf 5}.4.7)=0$&                  &$q({\bf 5}.6.7)=0$&$q({\bf 5}.7.7)=0$\\
                &$q({\bf 5}.3.8)=0$&$q({\bf 5}.4.8)=0$&                  &$q({\bf 5}.6.8)=0$&$q({\bf 5}.7.8)=0$\\
                &$q({\bf 5}.3.9)=0$&$q({\bf 5}.4.9)=0$&                  &$q({\bf 5}.6.9)=0$&$q({\bf 5}.7.9)=0$\\
                &$q({\bf 5}.3.10)=0$&$q({\bf 5}.4.10)=0$&                 &$q({\bf 5}.6.10)=0$&$q({\bf 5}.7.10)=0$\\
                &                  &                &                   &$q({\bf 5}.6.11)=0$&                  \\
                &                  &                &                   &$q({\bf 5}.6.12)=0$&                  \\
                &                  &                &                   &$q({\bf 5}.6.13)=0$&                  \\
                &                  &                &                   &$q({\bf 5}.6.14)=0$&                  \\
                &                  &                &                   &$q({\bf 5}.6.15)=0$&                  \\
\noalign{\smallskip}\hline\noalign{\smallskip}
\end{tabular}
\vspace{14mm}
\end{table}

\begin{table}
\centering
\caption{$T-t=0.008.$ Stochastic integral 
$I_{(00000)T,t}^{(i_1i_2i_3i_4i_5)}.$ The condition (\ref{tred3}).}
\label{tab:6}      
\begin{tabular}{p{1.9cm}p{2.2cm}p{2.2cm}p{2.2cm}p{2.2cm}p{2.2cm}}
\hline\noalign{\smallskip}
$q({\bf 5}.1)=\fbox{1}$&$q({\bf 5}.3.1)=1$&$q({\bf 5}.4.1)=0$&$q({\bf 5}.5.1)=0$&$q({\bf 5}.6.1)=1$&$q({\bf 5}.7.1)=0$\\
                &$q({\bf 5}.3.2)=1$&$q({\bf 5}.4.2)=0$&$q({\bf 5}.5.2)=0$&$q({\bf 5}.6.2)=1$&$q({\bf 5}.7.2)=0$\\
                &$q({\bf 5}.3.3)=1$&$q({\bf 5}.4.3)=0$&$q({\bf 5}.5.3)=0$&$q({\bf 5}.6.3)=1$&$q({\bf 5}.7.3)=0$\\
                &$q({\bf 5}.3.4)=1$&$q({\bf 5}.4.4)=0$&$q({\bf 5}.5.4)=0$&$q({\bf 5}.6.4)=1$&$q({\bf 5}.7.4)=0$\\
                &$q({\bf 5}.3.5)=1$&$q({\bf 5}.4.5)=0$&$q({\bf 5}.5.5)=0$&$q({\bf 5}.6.5)=1$&$q({\bf 5}.7.5)=0$\\
                &$q({\bf 5}.3.6)=1$&$q({\bf 5}.4.6)=0$&                  &$q({\bf 5}.6.6)=1$&$q({\bf 5}.7.6)=0$\\
                &$q({\bf 5}.3.7)=1$&$q({\bf 5}.4.7)=0$&                  &$q({\bf 5}.6.7)=1$&$q({\bf 5}.7.7)=0$\\
                &$q({\bf 5}.3.8)=1$&$q({\bf 5}.4.8)=0$&                  &$q({\bf 5}.6.8)=1$&$q({\bf 5}.7.8)=0$\\
                &$q({\bf 5}.3.9)=1$&$q({\bf 5}.4.9)=0$&                  &$q({\bf 5}.6.9)=1$&$q({\bf 5}.7.9)=0$\\
                &$q({\bf 5}.3.10)=1$&$q({\bf 5}.4.10)=0$&                 &$q({\bf 5}.6.10)=1$&$q({\bf 5}.7.10)=0$\\
                &                  &                &                   &$q({\bf 5}.6.11)=1$&                  \\
                &                  &                &                   &$q({\bf 5}.6.12)=1$&                  \\
                &                  &                &                   &$q({\bf 5}.6.13)=1$&                  \\
                &                  &                &                   &$q({\bf 5}.6.14)=1$&                  \\
                &                  &                &                   &$q({\bf 5}.6.15)=1$&                  \\
\noalign{\smallskip}\hline\noalign{\smallskip}
\end{tabular}
\vspace{9mm}
\end{table}

\begin{table}
\centering
\caption{$T-t=0.0045.$ Stochastic integral 
$I_{(00000)T,t}^{(i_1i_2i_3i_4i_5)}.$ The condition (\ref{tred3}).}
\label{tab:7}      
\begin{tabular}{p{1.9cm}p{2.2cm}p{2.2cm}p{2.2cm}p{2.2cm}p{2.2cm}}
\hline\noalign{\smallskip}
$q({\bf 5}.1)=\fbox{4}$&$q({\bf 5}.3.1)=4$&$q({\bf 5}.4.1)=2$&$q({\bf 5}.5.1)=1$&$q({\bf 5}.6.1)=2$&$q({\bf 5}.7.1)=1$\\
                &$q({\bf 5}.3.2)=4$&$q({\bf 5}.4.2)=3$&$q({\bf 5}.5.2)=2$&$q({\bf 5}.6.2)=4$&$q({\bf 5}.7.2)=3$\\
                &$q({\bf 5}.3.3)=4$&$q({\bf 5}.4.3)=3$&$q({\bf 5}.5.3)=2$&$q({\bf 5}.6.3)=3$&$q({\bf 5}.7.3)=2$\\
                &$q({\bf 5}.3.4)=4$&$q({\bf 5}.4.4)=2$&$q({\bf 5}.5.4)=2$&$q({\bf 5}.6.4)=3$&$q({\bf 5}.7.4)=2$\\
                &$q({\bf 5}.3.5)=3$&$q({\bf 5}.4.5)=3$&$q({\bf 5}.5.5)=1$&$q({\bf 5}.6.5)=3$&$q({\bf 5}.7.5)=3$\\
                &$q({\bf 5}.3.6)=4$&$q({\bf 5}.4.6)=3$&                  &$q({\bf 5}.6.6)=4$&$q({\bf 5}.7.6)=3$\\
                &$q({\bf 5}.3.7)=4$&$q({\bf 5}.4.7)=2$&                  &$q({\bf 5}.6.7)=4$&$q({\bf 5}.7.7)=1$\\
                &$q({\bf 5}.3.8)=3$&$q({\bf 5}.4.8)=4$&                  &$q({\bf 5}.6.8)=2$&$q({\bf 5}.7.8)=4$\\
                &$q({\bf 5}.3.9)=4$&$q({\bf 5}.4.9)=3$&                  &$q({\bf 5}.6.9)=4$&$q({\bf 5}.7.9)=3$\\
                &$q({\bf 5}.3.10)=3$&$q({\bf 5}.4.10)=3$&                 &$q({\bf 5}.6.10)=4$&$q({\bf 5}.7.10)=2$\\
                &                  &                &                   &$q({\bf 5}.6.11)=3$&                  \\
                &                  &                &                   &$q({\bf 5}.6.12)=3$&                  \\
                &                  &                &                   &$q({\bf 5}.6.13)=2$&                  \\
                &                  &                &                   &$q({\bf 5}.6.14)=4$&                  \\
                &                  &                &                   &$q({\bf 5}.6.15)=3$&                  \\
\noalign{\smallskip}\hline\noalign{\smallskip}
\end{tabular}
\vspace{14mm}
\end{table}

\begin{table}
\centering
\caption{$T-t=0.0042.$ Stochastic integral 
$I_{(00000)T,t}^{(i_1i_2i_3i_4i_5)}.$ The condition (\ref{tred3}).}
\label{tab:8}      
\begin{tabular}{p{1.9cm}p{2.2cm}p{2.2cm}p{2.2cm}p{2.2cm}p{2.2cm}}
\hline\noalign{\smallskip}
$q({\bf 5}.1)=\fbox{5}$&$q({\bf 5}.3.1)=5$&$q({\bf 5}.4.1)=2$&$q({\bf 5}.5.1)=1$&$q({\bf 5}.6.1)=2$&$q({\bf 5}.7.1)=1$\\
                &$q({\bf 5}.3.2)=5$&$q({\bf 5}.4.2)=4$&$q({\bf 5}.5.2)=2$&$q({\bf 5}.6.2)=4$&$q({\bf 5}.7.2)=3$\\
                &$q({\bf 5}.3.3)=5$&$q({\bf 5}.4.3)=4$&$q({\bf 5}.5.3)=2$&$q({\bf 5}.6.3)=3$&$q({\bf 5}.7.3)=2$\\
                &$q({\bf 5}.3.4)=5$&$q({\bf 5}.4.4)=2$&$q({\bf 5}.5.4)=2$&$q({\bf 5}.6.4)=4$&$q({\bf 5}.7.4)=2$\\
                &$q({\bf 5}.3.5)=3$&$q({\bf 5}.4.5)=3$&$q({\bf 5}.5.5)=1$&$q({\bf 5}.6.5)=3$&$q({\bf 5}.7.5)=3$\\
                &$q({\bf 5}.3.6)=4$&$q({\bf 5}.4.6)=4$&                  &$q({\bf 5}.6.6)=5$&$q({\bf 5}.7.6)=3$\\
                &$q({\bf 5}.3.7)=5$&$q({\bf 5}.4.7)=2$&                  &$q({\bf 5}.6.7)=5$&$q({\bf 5}.7.7)=1$\\
                &$q({\bf 5}.3.8)=3$&$q({\bf 5}.4.8)=5$&                  &$q({\bf 5}.6.8)=2$&$q({\bf 5}.7.8)=4$\\
                &$q({\bf 5}.3.9)=5$&$q({\bf 5}.4.9)=3$&                  &$q({\bf 5}.6.9)=4$&$q({\bf 5}.7.9)=3$\\
                &$q({\bf 5}.3.10)=4$&$q({\bf 5}.4.10)=4$&                 &$q({\bf 5}.6.10)=5$&$q({\bf 5}.7.10)=2$\\
                &                  &                &                   &$q({\bf 5}.6.11)=4$&                  \\
                &                  &                &                   &$q({\bf 5}.6.12)=3$&                  \\
                &                  &                &                   &$q({\bf 5}.6.13)=2$&                  \\
                &                  &                &                   &$q({\bf 5}.6.14)=4$&                  \\
                &                  &                &                   &$q({\bf 5}.6.15)=3$&                  \\
\noalign{\smallskip}\hline\noalign{\smallskip}
\end{tabular}
\vspace{9mm}
\end{table}

\begin{table}
\centering
\caption{$T-t=0.0035.$ Stochastic integral 
$I_{(00000)T,t}^{(i_1i_2i_3i_4i_5)}.$ The condition (\ref{tred3}).}
\label{tab:9}      
\begin{tabular}{p{1.9cm}p{2.2cm}p{2.2cm}p{2.2cm}p{2.2cm}p{2.2cm}}
\hline\noalign{\smallskip}
$q({\bf 5}.1)=\fbox{6}$&$q({\bf 5}.3.1)=6$&$q({\bf 5}.4.1)=3$&$q({\bf 5}.5.1)=1$&$q({\bf 5}.6.1)=3$&$q({\bf 5}.7.1)=3$\\
                &$q({\bf 5}.3.2)=6$&$q({\bf 5}.4.2)=4$&$q({\bf 5}.5.2)=3$&$q({\bf 5}.6.2)=6$&$q({\bf 5}.7.2)=4$\\
                &$q({\bf 5}.3.3)=6$&$q({\bf 5}.4.3)=4$&$q({\bf 5}.5.3)=3$&$q({\bf 5}.6.3)=4$&$q({\bf 5}.7.3)=3$\\
                &$q({\bf 5}.3.4)=6$&$q({\bf 5}.4.4)=2$&$q({\bf 5}.5.4)=3$&$q({\bf 5}.6.4)=4$&$q({\bf 5}.7.4)=2$\\
                &$q({\bf 5}.3.5)=4$&$q({\bf 5}.4.5)=4$&$q({\bf 5}.5.5)=1$&$q({\bf 5}.6.5)=4$&$q({\bf 5}.7.5)=4$\\
                &$q({\bf 5}.3.6)=6$&$q({\bf 5}.4.6)=4$&                  &$q({\bf 5}.6.6)=6$&$q({\bf 5}.7.6)=4$\\
                &$q({\bf 5}.3.7)=6$&$q({\bf 5}.4.7)=3$&                  &$q({\bf 5}.6.7)=6$&$q({\bf 5}.7.7)=1$\\
                &$q({\bf 5}.3.8)=4$&$q({\bf 5}.4.8)=6$&                  &$q({\bf 5}.6.8)=3$&$q({\bf 5}.7.8)=6$\\
                &$q({\bf 5}.3.9)=6$&$q({\bf 5}.4.9)=4$&                  &$q({\bf 5}.6.9)=6$&$q({\bf 5}.7.9)=4$\\
                &$q({\bf 5}.3.10)=4$&$q({\bf 5}.4.10)=4$&                 &$q({\bf 5}.6.10)=6$&$q({\bf 5}.7.10)=3$\\
                &                  &                &                   &$q({\bf 5}.6.11)=4$&                  \\
                &                  &                &                   &$q({\bf 5}.6.12)=4$&                  \\
                &                  &                &                   &$q({\bf 5}.6.13)=3$&                  \\
                &                  &                &                   &$q({\bf 5}.6.14)=6$&                  \\
                &                  &                &                   &$q({\bf 5}.6.15)=4$&                  \\
\noalign{\smallskip}\hline\noalign{\smallskip}
\end{tabular}
\vspace{9mm}
\end{table}

\begin{table}
\centering
\caption{Milstein scheme. Stochastic integral $I_{(00)T,t}^{(i_1i_2)}$.
The condition (\ref{tred1}).}
\label{tab:10}      
\begin{tabular}{p{1.7cm}p{1.7cm}p{1.7cm}p{1.7cm}p{1.7cm}}
\hline\noalign{\smallskip}
$T-t$&$2^{-1}$&$2^{-4}$&$2^{-8}$&$2^{-12}$\\
\noalign{\smallskip}\hline\noalign{\smallskip}
$q({\bf 2}.1.a)$&$1$&$2$&$32$&$512$\\
\noalign{\smallskip}\hline\noalign{\smallskip}
\end{tabular}
\vspace{9mm}
\end{table}

\begin{table}
\centering
\caption{Scheme with strong order $1.5$. Stochastic integrals $I_{(00)T,t}^{(i_1i_2)},$
$I_{(000)T,t}^{(i_1i_2i_3)}$. The condition (\ref{tred1}).}
\label{tab:11}      
\begin{tabular}{p{1.7cm}p{1.7cm}p{1.7cm}p{1.7cm}p{1.7cm}}
\hline\noalign{\smallskip}
$T-t$&$2^{-1}$&$2^{-3}$&$2^{-5}$&$2^{-8}$\\
\noalign{\smallskip}\hline\noalign{\smallskip}
$q({\bf 2}.1.a)$&$1$&$8$&$128$&$8192$\\
\hline\noalign{\smallskip}
$q({\bf 3}.1.a)$&$\fbox{0}$&$\fbox{1}$&$\fbox{4}$&$\fbox{32}$\\
$q({\bf 3}.3.1.a)$&$0$&$0$&$2$&$16$\\
$q({\bf 3}.3.2.a)$&$0$&$0$&$2$&$16$\\
$q({\bf 3}.3.3.a)$&$0$&$0$&$4$&$\fbox{33}$\\
\noalign{\smallskip}\hline\noalign{\smallskip}
\end{tabular}
\vspace{9mm}
\end{table}

\begin{table}
\centering
\caption{Scheme with strong order $2.0$. 
Stochastic integrals $I_{(00)T,t}^{(i_1i_2)},$
$I_{(000)T,t}^{(i_1i_2i_3)},$ 
$I_{(01)T,t}^{(i_1i_2)},$ 
$I_{(10)T,t}^{(i_1i_2)},$ $I_{(0000)T,t}^{(i_1i_2i_3i_4)}.$
The condition (\ref{tred2}).}
\label{tab:12}      
\begin{tabular}{p{1.7cm}p{1.7cm}p{1.7cm}p{1.7cm}p{1.7cm}}
\hline\noalign{\smallskip}
$T-t$&$2^{-1}$&$2^{-2}$&$2^{-3}$&$2^{-4}$\\
\noalign{\smallskip}\hline\noalign{\smallskip}
$q({\bf 2}.1.a)$&$1$&$8$&$64$&$512$\\
\hline\noalign{\smallskip}
$q({\bf 3}.1.a)$&$\fbox{0}$&$\fbox{2}$&$\fbox{8}$&$\fbox{32}$\\
$q({\bf 3}.3.1.a)$&$0$&$1$&$4$&$16$\\
$q({\bf 3}.3.2.a)$&$0$&$1$&$4$&$16$\\
$q({\bf 3}.3.3.a)$&$0$&$2$&$8$&$\fbox{33}$\\
\hline\noalign{\smallskip}
$q({\bf 2}.1.b)$&$\fbox{0}$&$\fbox{0}$&$\fbox{1}$&$\fbox{1}$\\
$q({\bf 2}.2.b)$&$0$&$0$&$0$&$0$\\
\noalign{\smallskip}\hline\noalign{\smallskip}
$q({\bf 2}.1.c)$&$\fbox{0}$&$\fbox{0}$&$\fbox{0}$&$\fbox{0}$\\
$q({\bf 2}.2.c)$&$0$&$0$&$0$&$0$\\
\hline\noalign{\smallskip}
$q({\bf 4}.1)$&$\fbox{0}$&$\fbox{0}$&$\fbox{0}$&$\fbox{0}$\\
$q({\bf 4}.3.1)$&$0$&$0$&$0$&$0$\\
$\cdot\cdot\cdot$&$\cdot\cdot\cdot$&$\cdot\cdot\cdot$&$\cdot\cdot\cdot$&$\cdot\cdot\cdot$\\
$q({\bf 4}.3.6)$&$0$&$0$&$0$&$0$\\
$q({\bf 4}.4.1)$&$0$&$0$&$0$&$0$\\
$\cdot\cdot\cdot$&$\cdot\cdot\cdot$&$\cdot\cdot\cdot$&$\cdot\cdot\cdot$&$\cdot\cdot\cdot$\\
$q({\bf 4}.4.4)$&$0$&$0$&$0$&$0$\\
$q({\bf 4}.5.1)$&$0$&$0$&$0$&$0$\\
$q({\bf 4}.5.2)$&$0$&$0$&$0$&$0$\\
$q({\bf 4}.5.3)$&$0$&$0$&$0$&$0$\\
\noalign{\smallskip}\hline\noalign{\smallskip}
\end{tabular}
\vspace{3mm}
\end{table}

\begin{table}
\centering
\caption{Scheme with strong order $2.5$. 
Stochastic integrals $I_{(00)T,t}^{(i_1i_2)},$
$I_{(000)T,t}^{(i_1i_2i_3)},$
$I_{(01)T,t}^{(i_1i_2)},$ 
$I_{(10)T,t}^{(i_1i_2)},$ $I_{(0000)T,t}^{(i_1i_2i_4)},$
$I_{(001)T,t}^{(i_1i_2i_3)},$
$I_{(010)T,t}^{(i_1i_2i_3)},$
$I_{(100)T,t}^{(i_1i_2i_3)},$
$I_{(00000)T,t}^{(i_1i_2i_3i_4i_5)}.$ The condition (\ref{tred3}).}
\label{tab:13}      
\vspace{5mm}
\begin{tabular}{p{1.7cm}p{1.7cm}p{1.7cm}p{1.7cm}p{1.7cm}}
\hline\noalign{\smallskip}
$T-t$&$2^{-1}$&$2^{-3/2}$&$2^{-2}$&$2^{-5/2}$\\
\noalign{\smallskip}\hline\noalign{\smallskip}
$q({\bf 2}.1.a)$&$2$&$8$&$32$&$128$\\
\hline\noalign{\smallskip}
$q({\bf 3}.1.a)$&$\fbox{1}$&$\fbox{3}$&$\fbox{8}$&$\fbox{23}$\\
$q({\bf 3}.3.1.a)$&$0$&$1$&$4$&$11$\\
$q({\bf 3}.3.2.a)$&$0$&$1$&$4$&$11$\\
$q({\bf 3}.3.3.a)$&$0$&$3$&$8$&$23$\\
\hline\noalign{\smallskip}
$q({\bf 2}.1.b)$&$\fbox{0}$&$\fbox{1}$&$\fbox{1}$&$\fbox{2}$\\
$q({\bf 2}.2.b)$&$0$&$0$&$0$&$0$\\
\hline\noalign{\smallskip}
$q({\bf 2}.1.c)$&$\fbox{0}$&$\fbox{0}$&$\fbox{0}$&$\fbox{2}$\\
$q({\bf 2}.2.c)$&$0$&$0$&$0$&$0$\\
\hline\noalign{\smallskip}
$q({\bf 4}.1)$&$\fbox{0}$&$\fbox{0}$&$\fbox{0}$&$\fbox{2}$\\
$q({\bf 4}.3.1)$&$0$&$0$&$0$&$1$\\
$q({\bf 4}.3.2)$&$0$&$0$&$0$&$1$\\
$q({\bf 4}.3.3)$&$0$&$0$&$0$&$2$\\
$q({\bf 4}.3.4)$&$0$&$0$&$0$&$1$\\
$q({\bf 4}.3.5)$&$0$&$0$&$0$&$1$\\
$q({\bf 4}.3.6)$&$0$&$0$&$0$&$1$\\
$q({\bf 4}.4.1)$&$0$&$0$&$0$&$0$\\
$q({\bf 4}.4.2)$&$0$&$0$&$0$&$0$\\
$q({\bf 4}.4.3)$&$0$&$0$&$0$&$0$\\
$q({\bf 4}.4.4)$&$0$&$0$&$0$&$0$\\
$q({\bf 4}.5.1)$&$0$&$0$&$0$&$1$\\
$q({\bf 4}.5.2)$&$0$&$0$&$0$&$1$\\
$q({\bf 4}.5.3)$&$0$&$0$&$0$&$1$\\
\noalign{\smallskip}\hline\noalign{\smallskip}
$q({\bf 3}.1.b)$&$\fbox{0}$&$\fbox{0}$&$\fbox{0}$&$\fbox{0}$\\
$q({\bf 3}.2.b)$&$0$&$0$&$0$&$0$\\
$q({\bf 3}.3.1.b)$&$0$&$0$&$0$&$0$\\
$q({\bf 3}.3.2.b)$&$0$&$0$&$0$&$0$\\
$q({\bf 3}.3.3.b)$&$0$&$0$&$0$&$0$\\
\noalign{\smallskip}\hline\noalign{\smallskip}
$q({\bf 3}.1.c)$&$\fbox{0}$&$\fbox{0}$&$\fbox{0}$&$\fbox{0}$\\
$q({\bf 3}.2.c)$&$0$&$0$&$0$&$0$\\
$q({\bf 3}.3.1.c)$&$0$&$0$&$0$&$0$\\
$q({\bf 3}.3.2.c)$&$0$&$0$&$0$&$0$\\
$q({\bf 3}.3.3.c)$&$0$&$0$&$0$&$0$\\
\noalign{\smallskip}\hline\noalign{\smallskip}
$q({\bf 3}.1.d)$&$\fbox{0}$&$\fbox{0}$&$\fbox{0}$&$\fbox{0}$\\
$q({\bf 3}.2.d)$&$0$&$0$&$0$&$0$\\
$q({\bf 3}.3.1.d)$&$0$&$0$&$0$&$0$\\
$q({\bf 3}.3.2.d)$&$0$&$0$&$0$&$0$\\
$q({\bf 3}.3.3.d)$&$0$&$0$&$0$&$0$\\
\noalign{\smallskip}\hline\noalign{\smallskip}
\begin{tabular}{p{10.2cm}}
\vspace{6mm}
~~~~~~~~~All numbers $q(\alpha)$ 
for $I_{(00000)T,t}^{(i_1i_2i_3i_4i_5)}$
are equal to zero\\
\end{tabular}
\end{tabular}
\end{table}

In Tables 1--13, we can see the results of numerical 
experiments. These results confirm (in most situations) the inequalities
(\ref{som1})-(\ref{som8})).

Let us show by numerical experiments that
we can choose the minimal natural numbers $p$  
satisfying the inequalities
(\ref{tred1})-(\ref{tred3}) only for the values
$E^p_2$, $E^p_3$, $E^p_4$, $E^p_5$ 
defined by the relations 
(\ref{prod1}),
(\ref{prod3}),
(\ref{prod5}),
(\ref{prod7}),
(\ref{prod11}),
(\ref{prod16}),
(\ref{prod21}),
(\ref{prod26}),
(\ref{prod40}).
At that, we can suppose $i_1,\ldots,i_5=1,\ldots,m$ in these
relations and use the above numbers $p$ for all remaining cases. This means that we can
ignore all the formulas
(\ref{prod4}),
(\ref{prod6}),
(\ref{prod8})-(\ref{prod10}),
(\ref{prod12})-(\ref{prod15}),
(\ref{prod17})-(\ref{prod20}),
(\ref{prod22})-(\ref{prod25}),
(\ref{prod27})-(\ref{prod39}),
(\ref{prod41})-(\ref{prod90}).
As a result, we will not get a noticeable loss
of the mean-square approximation accuracy for iterated 
Ito stochastic integrals.
The detailed numerical confirmation of the above hypothesis
can be found in Tables 14--22.

Taking into account the results of this article,
we can recommend the following 
conditions for choosing the minimal
natural numbers 
$q, q_1, q_2,\ldots, q_{7}$ (see Sect.~5) for the numerical schemes
(\ref{al1})--(\ref{al4}) (constant $C$ (see below)
has the same meaning as in the 
condition (\ref{uslov})).

\vspace{10mm}

\centerline{\bf Milstein scheme (\ref{al1})}

\vspace{3mm}

$$
\frac{(T-t)^2}{2}\left(\frac{1}{2}-\sum_{i=1}^q
\frac{1}{4i^2-1}\right)\le C (T-t)^3.
$$

\vspace{10mm}

\centerline{\bf Strong Taylor--Ito scheme with convergence order 1.5
(\ref{al2})}

\vspace{4mm}

$$
\frac{(T-t)^2}{2}\left(\frac{1}{2}-\sum_{i=1}^q
\frac{1}{4i^2-1}\right)\le C (T-t)^4,
$$

\vspace{-3mm}

$$
(T-t)^3\left(\frac{1}{6}-\frac{1}{64}\sum_{j_1,j_2,j_3=0}^{q_1}
(2j_1+1)(2j_2+1)(2j_3+1)
\left(\bar C_{j_3j_2j_1}^{000}\right)^2\right)\le 
$$

$$
\le C(T-t)^4.
$$

\vspace{13mm}

\centerline{\bf Strong Taylor--Ito scheme with convergence order 2.0
(\ref{al3})}

\vspace{4mm}

$$
\frac{(T-t)^2}{2}\left(\frac{1}{2}-\sum_{i=1}^q
\frac{1}{4i^2-1}\right)\le C (T-t)^5,
$$

\vspace{4mm}

$$
(T-t)^3\left(\frac{1}{6}-\frac{1}{64}\sum_{j_1,j_2,j_3=0}^{q_1}
(2j_1+1)(2j_2+1)(2j_3+1)
\left(\bar C_{j_3j_2j_1}^{000}\right)^2\right)\le C(T-t)^5,
$$

\vspace{6mm}

$$
(T-t)^4\left(\frac{1}{4}-\frac{1}{64}\sum_{j_1,j_2=0}^{q_2}
(2j_1+1)(2j_2+1)
\left(\bar C_{j_2j_1}^{01}\right)^2\right)\le C(T-t)^5,
$$

\vspace{6mm}

$$
(T-t)^4\left(\frac{1}{12}-\frac{1}{64}\sum_{j_1,j_2=0}^{q_2}
(2j_1+1)(2j_2+1)
\left(\bar C_{j_2j_1}^{10}\right)^2\right)\le C(T-t)^5,
$$

\vspace{6mm}

$$
(T-t)^4\left(\frac{1}{24}-\frac{1}{256}\sum_{j_1,\ldots,j_4=0}^{q_3}
(2j_1+1)(2j_2+1)(2j_3+1)(2j_4+1)
\left(\bar C_{j_4\ldots j_1}^{0000}\right)^2\right)\le 
$$

$$
\le C(T-t)^5.
$$

\vspace{15mm}

\begin{table}
\centering
\caption{Stochastic integral $I_{(000)T,t}^{(i_1i_2i_3)}$. The values $E_3^p/(T-t)^3
\stackrel{\sf def}{=}E$.}
\label{tab:14}      
\begin{tabular}{p{1.7cm}p{1.7cm}p{1.7cm}p{1.7cm}p{1.7cm}p{1.7cm}p{1.7cm}}
\hline\noalign{\smallskip}
$T-t$&$0.011$&$0.008$&$0.0045$&$0.0035$&$0.0027$&$0.0025$\\
\noalign{\smallskip}\hline\noalign{\smallskip}
$q({\bf 3}.1.a)$&$12$&$16$&$28$&$36$&$47$&$50$\\
$E$&0.010154&0.007681&0.004433&0.003456&0.002652&0.002494\\
$q({\bf 3}.3.1.a)$&$12$&$16$&$28$&$36$&$47$&$50$\\
$E$&0.005077&0.003841&0.002216&0.001728&0.001326&0.001247\\
$q({\bf 3}.3.2.a)$&$12$&$16$&$28$&$36$&$47$&$50$\\
$E$&0.005077&0.003841&0.002216&0.001728&0.001326&0.001247\\
$q({\bf 3}.3.3.a)$&$12$&$16$&$28$&$36$&$47$&$50$\\
$E$&0.010308&0.007787&0.004480&0.003488&0.002673&0.002513\\
\noalign{\smallskip}\hline\noalign{\smallskip}
\end{tabular}
\vspace{6mm}
\end{table}

\begin{table}
\centering
\caption{Stochastic integral $I_{(0000)T,t}^{(i_1i_2i_3i_4)}$. The values $E_4^p/(T-t)^4
\stackrel{\sf def}{=}E$.}
\label{tab:15}      
\begin{tabular}{p{1.7cm}p{1.7cm}p{1.7cm}p{1.7cm}p{1.7cm}}
\hline\noalign{\smallskip}
$T-t$&$0.011$&$0.008$&$0.0045$&$0.0042$\\
\noalign{\smallskip}\hline\noalign{\smallskip}
$q({\bf 4}.1)$&6&8&14&15\\
$E$&0.009636&0.007425&0.004378&0.004096\\
$q({\bf 4}.3.1)$&6&8&14&15\\
$E$&0.006771&0.005191&0.003041&0.002843\\
$q({\bf 4}.3.2)$&6&8&14&15\\
$E$&0.009722&0.007502&0.004424&0.004139\\
$q({\bf 4}.3.3)$&6&8&14&15\\
$E$&0.009641&0.007427&0.004379&0.004097\\
$q({\bf 4}.3.4)$&6&8&14&15\\
$E$&0.005997&0.004614&0.002720&0.002545\\
$q({\bf 4}.3.5)$&6&8&14&15\\
$E$&0.009722&0.007502&0.004424&0.004139\\
$q({\bf 4}.3.6)$&6&8&14&15\\
$E$&0.006771&0.005191&0.003041&0.002843\\
$q({\bf 4}.4.1)$&6&8&14&15\\
$E$&0.003095&0.002364&0.001379&0.001290\\
$q({\bf 4}.4.2)$&6&8&14&15\\
$E$&0.003095&0.002364&0.001379&0.001290\\
$q({\bf 4}.4.3)$&6&8&14&15\\
$E$&0.006885&0.005282&0.003090&0.002889\\
$q({\bf 4}.4.4)$&6&8&14&15\\
$E$&0.006885&0.005282&0.003090&0.002889\\
$q({\bf 4}.5.1)$&6&8&14&15\\
$E$&0.003690&0.002834&0.001663&0.001555\\
$q({\bf 4}.5.2)$&6&8&14&15\\
$E$&0.009756&0.007545&0.004457&0.004170\\
$q({\bf 4}.5.3)$&6&8&14&15\\
$E$&0.006010&0.004621&0.002722&0.002547\\
\noalign{\smallskip}\hline\noalign{\smallskip}
\end{tabular}
\vspace{15mm}
\end{table}

\begin{table}
\centering
\caption{Stochastic integrals $I_{(01)T,t}^{(i_1i_2)},$
$I_{(10)T,t}^{(i_1i_2)}$. The values $E_2^p/(T-t)^4\stackrel{\sf def}{=}E$.}
\label{tab:16}      
\begin{tabular}{p{1.7cm}p{1.7cm}p{1.7cm}p{1.7cm}}
\hline\noalign{\smallskip}
$T-t$&$0.010$&$0.005$&$0.0025$\\
\noalign{\smallskip}\hline\noalign{\smallskip}
$q({\bf 2}.1.b)$&4&8&16\\
$E$&0.008950&0.004660&0.002383\\
$q({\bf 2}.2.b)$&4&8&16\\
$E$&0.000042&0.000006&0.000001\\
$q({\bf 2}.1.c)$&4&8&16\\
$E$&0.008950&0.004660&0.002383\\
$q({\bf 2}.2.c)$&4&8&16\\
$E$&0.000042&0.000006&0.000001\\
\noalign{\smallskip}\hline\noalign{\smallskip}
\end{tabular}
\vspace{6mm}
\end{table}

\begin{table}
\centering
\caption{$T-t=0.011.$ Stochastic integral 
$I_{(00000)T,t}^{(i_1i_2i_3i_4i_5)}.$ The values $E_5^p/(T-t)^5\stackrel{\sf def}{=}E$.}
\label{tab:17}      
\begin{tabular}{p{1.9cm}p{2.2cm}p{2.2cm}p{2.2cm}p{2.2cm}p{2.2cm}}
\hline\noalign{\smallskip}
\\
$q({\bf 5}.1)=0$&$q({\bf 5}.3.1)=0$&$q({\bf 5}.4.1)=0$&$q({\bf 5}.5.1)=0$&$q({\bf 5}.6.1)=0$&$q({\bf 5}.7.1)=0$\\
$E$=0.008264    &$E$=0.008195      &$E$=0.007917      &$E$=0.006667      &$E$=0.008056      &$E$=0.007500       \\
\\                
                &$q({\bf 5}.3.2)=0$&$q({\bf 5}.4.2)=0$&$q({\bf 5}.5.2)=0$&$q({\bf 5}.6.2)=0$&$q({\bf 5}.7.2)=0$\\
                &$E$=0.008195      &$E$=0.007917      &$E$=0.006667      &$E$=0.008056      &$E$=0.007500       \\
\\
                &$q({\bf 5}.3.3)=0$&$q({\bf 5}.4.3)=0$&$q({\bf 5}.5.3)=0$&$q({\bf 5}.6.3)=0$&$q({\bf 5}.7.3)=0$\\
                &$E$=0.008195      &$E$=0.007917      &$E$=0.006667      &$E$=0.008056      &$E$=0.007500       \\
\\
                &$q({\bf 5}.3.4)=0$&$q({\bf 5}.4.4)=0$&$q({\bf 5}.5.4)=0$&$q({\bf 5}.6.4)=0$&$q({\bf 5}.7.4)=0$\\
                &$E$=0.008195      &$E$=0.007917      &$E$=0.006667      &$E$=0.008056      &$E$=0.007500       \\
\\
                &$q({\bf 5}.3.5)=0$&$q({\bf 5}.4.5)=0$&$q({\bf 5}.5.5)=0$&$q({\bf 5}.6.5)=0$&$q({\bf 5}.7.5)=0$\\
                &$E$=0.008195      &$E$=0.007917      &$E$=0.006667      &$E$=0.008056      &$E$=0.007500       \\
\\
                &$q({\bf 5}.3.6)=0$&$q({\bf 5}.4.6)=0$&                  &$q({\bf 5}.6.6)=0$&$q({\bf 5}.7.6)=0$\\
                &$E$=0.008195      &$E$=0.007917      &                  &$E$=0.008056      &$E$=0.007500       \\
\\
                &$q({\bf 5}.3.7)=0$&$q({\bf 5}.4.7)=0$&                  &$q({\bf 5}.6.7)=0$&$q({\bf 5}.7.7)=0$\\
                &$E$=0.008195      &$E$=0.007917      &                  &$E$=0.008056      &$E$=0.007500       \\
\\
                &$q({\bf 5}.3.8)=0$&$q({\bf 5}.4.8)=0$&                  &$q({\bf 5}.6.8)=0$&$q({\bf 5}.7.8)=0$\\
                &$E$=0.008195      &$E$=0.007917      &                  &$E$=0.008056      &$E$=0.007500       \\
\\
                &$q({\bf 5}.3.9)=0$&$q({\bf 5}.4.9)=0$&                  &$q({\bf 5}.6.9)=0$&$q({\bf 5}.7.9)=0$\\
                &$E$=0.008195      &$E$=0.007917      &                  &$E$=0.008056      &$E$=0.007500        \\
\\
                &$q({\bf 5}.3.10)=0$&$q({\bf 5}.4.10)=0$&                 &$q({\bf 5}.6.10)=0$&$q({\bf 5}.7.10)=0$\\
                &$E$=0.008195       &$E$=0.007917       &                 &$E$=0.008056       &$E$=0.007500       \\
\\
                &                   &                   &                 &$q({\bf 5}.6.11)=0$&                  \\
                &                   &                   &                 &$E$=0.008056       &                  \\
\\
                &                   &                   &                 &$q({\bf 5}.6.12)=0$&                  \\
                &                   &                   &                 &$E$=0.008056       &                  \\
\\
                &                   &                   &                 &$q({\bf 5}.6.13)=0$&                  \\
                &                   &                   &                 &$E$=0.008056       &                  \\
\\
                &                   &                   &                 &$q({\bf 5}.6.14)=0$&                  \\
                &                   &                   &                 &$E$=0.008056       &                  \\
\\
                &                   &                   &                 &$q({\bf 5}.6.15)=0$&                  \\
                &                   &                   &                 &$E$=0.008056       &                   \\
\\
\noalign{\smallskip}\hline\noalign{\smallskip}                           
\end{tabular}
\vspace{9mm}
\end{table}

\begin{table}
\centering
\caption{$T-t=0.008.$ Stochastic integral 
$I_{(00000)T,t}^{(i_1i_2i_3i_4i_5)}.$ The values $E_5^p/(T-t)^5\stackrel{\sf def}{=}E$.}
\label{tab:18}      
\begin{tabular}{p{1.9cm}p{2.2cm}p{2.2cm}p{2.2cm}p{2.2cm}p{2.2cm}}
\hline\noalign{\smallskip}
\\
$q({\bf 5}.1)=1$&$q({\bf 5}.3.1)=1$&$q({\bf 5}.4.1)=1$&$q({\bf 5}.5.1)=1$&$q({\bf 5}.6.1)=1$&$q({\bf 5}.7.1)=1$\\
$E$=0.007590    &$E$=0.007570      &$E$=0.005488      &$E$=0.003272      &$E$=0.006052      &$E$=0.004175       \\
\\                
                &$q({\bf 5}.3.2)=1$&$q({\bf 5}.4.2)=1$&$q({\bf 5}.5.2)=1$&$q({\bf 5}.6.2)=1$&$q({\bf 5}.7.2)=1$\\
                &$E$=0.007300      &$E$=0.006701      &$E$=0.005292      &$E$=0.007058      &$E$=0.006105      \\
\\
                &$q({\bf 5}.3.3)=1$&$q({\bf 5}.4.3)=1$&$q({\bf 5}.5.3)=1$&$q({\bf 5}.6.3)=1$&$q({\bf 5}.7.3)=1$\\
                &$E$=0.007558      &$E$=0.006976      &$E$=0.005774      &$E$=0.007014      &$E$=0.006072      \\
\\
                &$q({\bf 5}.3.4)=1$&$q({\bf 5}.4.4)=1$&$q({\bf 5}.5.4)=1$&$q({\bf 5}.6.4)=1$&$q({\bf 5}.7.4)=1$\\
                &$E$=0.007570      &$E$=0.005995      &$E$=0.005292      &$E$=0.006467      &$E$=0.005955       \\
\\
                &$q({\bf 5}.3.5)=1$&$q({\bf 5}.4.5)=1$&$q({\bf 5}.5.5)=1$&$q({\bf 5}.6.5)=1$&$q({\bf 5}.7.5)=1$\\
                &$E$=0.007084      &$E$=0.006679      &$E$=0.003272      &$E$=0.007054      &$E$=0.006576       \\
\\
                &$q({\bf 5}.3.6)=1$&$q({\bf 5}.4.6)=1$&                  &$q({\bf 5}.6.6)=1$&$q({\bf 5}.7.6)=1$\\
                &$E$=0.007432      &$E$=0.006701      &                  &$E$=0.007260      &$E$=0.006105       \\
\\
                &$q({\bf 5}.3.7)=1$&$q({\bf 5}.4.7)=1$&                  &$q({\bf 5}.6.7)=1$&$q({\bf 5}.7.7)=1$\\
                &$E$=0.007558      &$E$=0.005488      &                  &$E$=0.007521      &$E$=0.003236       \\
\\
                &$q({\bf 5}.3.8)=1$&$q({\bf 5}.4.8)=1$&                  &$q({\bf 5}.6.8)=1$&$q({\bf 5}.7.8)=1$\\
                &$E$=0.007084      &$E$=0.007134      &                  &$E$=0.005819      &$E$=0.006797       \\
\\
                &$q({\bf 5}.3.9)=1$&$q({\bf 5}.4.9)=1$&                  &$q({\bf 5}.6.9)=1$&$q({\bf 5}.7.9)=1$\\
                &$E$=0.007300      &$E$=0.006679      &                  &$E$=0.007412      &$E$=0.006576        \\
\\
                &$q({\bf 5}.3.10)=1$&$q({\bf 5}.4.10)=1$&                 &$q({\bf 5}.6.10)=1$&$q({\bf 5}.7.10)=1$\\
                &$E$=0.006962       &$E$=0.006976       &                 &$E$=0.007260       &$E$=0.006072       \\
\\
                &                   &                   &                 &$q({\bf 5}.6.11)=1$&                  \\
                &                   &                   &                 &$E$=0.006467        &                  \\
\\
                &                   &                   &                 &$q({\bf 5}.6.12)=1$&                  \\
                &                   &                   &                 &$E$=0.007054        &                  \\
\\
                &                   &                   &                 &$q({\bf 5}.6.13)=1$&                  \\
                &                   &                   &                 &$E$=0.006052        &                  \\
\\
                &                   &                   &                 &$q({\bf 5}.6.14)=1$&                  \\
                &                   &                   &                 &$E$=0.007058        &                  \\
\\
                &                   &                   &                 &$q({\bf 5}.6.15)=1$&                  \\
                &                   &                   &                 &$E$=0.007014        &                   \\
\\
\noalign{\smallskip}\hline\noalign{\smallskip}                           
\end{tabular}
\vspace{9mm}
\end{table}

\begin{table}
\centering
\caption{$T-t=0.0045.$ Stochastic integral 
$I_{(00000)T,t}^{(i_1i_2i_3i_4i_5)}.$ The values $E_5^p/(T-t)^5\stackrel{\sf def}{=}E$.}
\label{tab:19}      
\begin{tabular}{p{1.9cm}p{2.2cm}p{2.2cm}p{2.2cm}p{2.2cm}p{2.2cm}}
\hline\noalign{\smallskip}
\\
$q({\bf 5}.1)=4$&$q({\bf 5}.3.1)=4$&$q({\bf 5}.4.1)=4$&$q({\bf 5}.5.1)=4$&$q({\bf 5}.6.1)=4$&$q({\bf 5}.7.1)=4$\\
$E$=0.004209    &$E$=0.004208      &$E$=0.002351      &$E$=0.001055      &$E$=0.002247      &$E$=0.002065       \\
\\                
                &$q({\bf 5}.3.2)=4$&$q({\bf 5}.4.2)=4$&$q({\bf 5}.5.2)=4$&$q({\bf 5}.6.2)=4$&$q({\bf 5}.7.2)=4$\\
                &$E$=0.004204      &$E$=0.003461      &$E$=0.002379      &$E$=0.004149      &$E$=0.003428      \\
\\
                &$q({\bf 5}.3.3)=4$&$q({\bf 5}.4.3)=4$&$q({\bf 5}.5.3)=4$&$q({\bf 5}.6.3)=4$&$q({\bf 5}.7.3)=4$\\
                &$E$=0.004212      &$E$=0.003460      &$E$=0.002624      &$E$=0.003168      &$E$=0.002256      \\
\\
                &$q({\bf 5}.3.4)=4$&$q({\bf 5}.4.4)=4$&$q({\bf 5}.5.4)=4$&$q({\bf 5}.6.4)=4$&$q({\bf 5}.7.4)=4$\\
                &$E$=0.004208      &$E$=0.001982      &$E$=0.002379      &$E$=0.003451      &$E$=0.001982        \\
\\
                &$q({\bf 5}.3.5)=4$&$q({\bf 5}.4.5)=4$&$q({\bf 5}.5.5)=4$&$q({\bf 5}.6.5)=4$&$q({\bf 5}.7.5)=4$\\
                &$E$=0.003161      &$E$=0.003189      &$E$=0.001055      &$E$=0.003160      &$E$=0.003191       \\
\\
                &$q({\bf 5}.3.6)=4$&$q({\bf 5}.4.6)=4$&                  &$q({\bf 5}.6.6)=4$&$q({\bf 5}.7.6)=4$\\
                &$E$=0.004180      &$E$=0.003461      &                  &$E$=0.004206      &$E$=0.003428        \\
\\
                &$q({\bf 5}.3.7)=4$&$q({\bf 5}.4.7)=4$&                  &$q({\bf 5}.6.7)=4$&$q({\bf 5}.7.7)=4$\\
                &$E$=0.004212      &$E$=0.002351      &                  &$E$=0.004214      &$E$=0.001318       \\
\\
                &$q({\bf 5}.3.8)=4$&$q({\bf 5}.4.8)=4$&                  &$q({\bf 5}.6.8)=4$&$q({\bf 5}.7.8)=4$\\
                &$E$=0.003161      &$E$=0.004201      &                  &$E$=0.002590      &$E$=0.004124       \\
\\
                &$q({\bf 5}.3.9)=4$&$q({\bf 5}.4.9)=4$&                  &$q({\bf 5}.6.9)=4$&$q({\bf 5}.7.9)=4$\\
                &$E$=0.004204      &$E$=0.003189      &                  &$E$=0.004180      &$E$=0.003191        \\
\\
                &$q({\bf 5}.3.10)=4$&$q({\bf 5}.4.10)=4$&                 &$q({\bf 5}.6.10)=4$&$q({\bf 5}.7.10)=4$\\
                &$E$=0.003456       &$E$=0.003460       &                 &$E$=0.004206       &$E$=0.002256       \\
\\
                &                   &                   &                 &$q({\bf 5}.6.11)=4$&                  \\
                &                   &                   &                 &$E$=0.003451        &                  \\
\\
                &                   &                   &                 &$q({\bf 5}.6.12)=4$&                  \\
                &                   &                   &                 &$E$=0.003160        &                  \\
\\
                &                   &                   &                 &$q({\bf 5}.6.13)=4$&                  \\
                &                   &                   &                 &$E$=0.002247        &                  \\
\\
                &                   &                   &                 &$q({\bf 5}.6.14)=4$&                  \\
                &                   &                   &                 &$E$=0.00414         &                  \\
\\
                &                   &                   &                 &$q({\bf 5}.6.15)=4$&                  \\
                &                   &                   &                 &$E$=0.003168        &                   \\
\\
\noalign{\smallskip}\hline\noalign{\smallskip}                           
\end{tabular}
\vspace{9mm}
\end{table}

\begin{table}
\centering
\caption{$T-t=0.0042.$ Stochastic integral 
$I_{(00000)T,t}^{(i_1i_2i_3i_4i_5)}.$ The values $E_5^p/(T-t)^5\stackrel{\sf def}{=}E$.}
\label{tab:20}      
\begin{tabular}{p{1.9cm}p{2.2cm}p{2.2cm}p{2.2cm}p{2.2cm}p{2.2cm}}
\hline\noalign{\smallskip}
\\
$q({\bf 5}.1)=5$&$q({\bf 5}.3.1)=5$&$q({\bf 5}.4.1)=5$&$q({\bf 5}.5.1)=5$&$q({\bf 5}.6.1)=5$&$q({\bf 5}.7.1)=5$\\
$E$=0.003557    &$E$=0.003556      &$E$=0.001940      &$E$=0.000863      &$E$=0.001863      &$E$=0.001728       \\
\\                
                &$q({\bf 5}.3.2)=5$&$q({\bf 5}.4.2)=5$&$q({\bf 5}.5.2)=5$&$q({\bf 5}.6.2)=5$&$q({\bf 5}.7.2)=5$\\
                &$E$=0.003564      &$E$=0.002910      &$E$=0.001969      &$E$=0.003539      &$E$=0.002897      \\
\\
                &$q({\bf 5}.3.3)=5$&$q({\bf 5}.4.3)=5$&$q({\bf 5}.5.3)=5$&$q({\bf 5}.6.3)=5$&$q({\bf 5}.7.3)=5$\\
                &$E$=0.003559      &$E$=0.002897      &$E$=0.002188      &$E$=0.002639      &$E$=0.001869      \\
\\
                &$q({\bf 5}.3.4)=5$&$q({\bf 5}.4.4)=5$&$q({\bf 5}.5.4)=5$&$q({\bf 5}.6.4)=5$&$q({\bf 5}.7.4)=5$\\
                &$E$=0.003556      &$E$=0.001642      &$E$=0.001969      &$E$=0.002903      &$E$=0.001641       \\
\\
                &$q({\bf 5}.3.5)=5$&$q({\bf 5}.4.5)=5$&$q({\bf 5}.5.5)=5$&$q({\bf 5}.6.5)=5$&$q({\bf 5}.7.5)=5$\\
                &$E$=0.002634      &$E$=0.002661      &$E$=0.000863      &$E$=0.002634      &$E$=0.002664       \\
\\
                &$q({\bf 5}.3.6)=5$&$q({\bf 5}.4.6)=5$&                  &$q({\bf 5}.6.6)=5$&$q({\bf 5}.7.6)=5$\\
                &$E$=0.003552      &$E$=0.002910       &                 &$E$=0.003566      &$E$=0.002897       \\
\\
                &$q({\bf 5}.3.7)=5$&$q({\bf 5}.4.7)=5$&                  &$q({\bf 5}.6.7)=5$&$q({\bf 5}.7.7)=5$\\
                &$E$=0.003559      &$E$=0.001940       &                 &$E$=0.003561      &$E$=0.001090       \\
\\
                &$q({\bf 5}.3.8)=5$&$q({\bf 5}.4.8)=5$&                  &$q({\bf 5}.6.8)=5$&$q({\bf 5}.7.8)=5$\\
                &$E$=0.002634      &$E$=0.003572       &                 &$E$=0.002155      &$E$=0.003531        \\
\\
                &$q({\bf 5}.3.9)=5$&$q({\bf 5}.4.9)=5$&                  &$q({\bf 5}.6.9)=5$&$q({\bf 5}.7.9)=5$\\
                &$E$=0.003564      &$E$=0.002661       &                 &$E$=0.003552      &$E$=0.002664        \\
\\
                &$q({\bf 5}.3.10)=5$&$q({\bf 5}.4.10)=5$&                 &$q({\bf 5}.6.10)=5$&$q({\bf 5}.7.10)=5$\\
                &$E$=0.002894       &$E$=0.002897        &                &$E$=0.003566       &$E$=0.001869        \\
\\
                &                   &                   &                 &$q({\bf 5}.6.11)=5$&                  \\
                &                   &                   &                 &$E$=0.002903        &                  \\
\\
                &                   &                   &                 &$q({\bf 5}.6.12)=5$&                  \\
                &                   &                   &                 &$E$=0.002634        &                  \\
\\
                &                   &                   &                 &$q({\bf 5}.6.13)=5$&                  \\
                &                   &                   &                 &$E$=0.001863        &                  \\
\\
                &                   &                   &                 &$q({\bf 5}.6.14)=5$&                  \\
                &                   &                   &                 &$E$=0.003539        &                  \\
\\
                &                   &                   &                 &$q({\bf 5}.6.15)=5$&                  \\
                &                   &                   &                 &$E$=0.002639        &                   \\
\\
\noalign{\smallskip}\hline\noalign{\smallskip}                           
\end{tabular}
\vspace{9mm}
\end{table}

\begin{table}
\centering
\caption{$T-t=0.0035.$ Stochastic integral 
$I_{(00000)T,t}^{(i_1i_2i_3i_4i_5)}.$ The values $E_5^p/(T-t)^5\stackrel{\sf def}{=}E$.}
\label{tab:21}      
\begin{tabular}{p{1.9cm}p{2.2cm}p{2.2cm}p{2.2cm}p{2.2cm}p{2.2cm}}
\hline\noalign{\smallskip}
\\
$q({\bf 5}.1)=6$&$q({\bf 5}.3.1)=6$&$q({\bf 5}.4.1)=6$&$q({\bf 5}.5.1)=6$&$q({\bf 5}.6.1)=6$&$q({\bf 5}.7.1)=6$\\
$E$=0.003071    &$E$=0.003071      &$E$=0.001650      &$E$=0.000729      &$E$=0.001591      &$E$=0.001591       \\
\\                
                &$q({\bf 5}.3.2)=6$&$q({\bf 5}.4.2)=6$&$q({\bf 5}.5.2)=6$&$q({\bf 5}.6.2)=6$&$q({\bf 5}.7.2)=6$\\
                &$E$=0.003083      &$E$=0.002503      &$E$=0.001676      &$E$=0.003074      &$E$=0.002500      \\
\\
                &$q({\bf 5}.3.3)=6$&$q({\bf 5}.4.3)=6$&$q({\bf 5}.5.3)=6$&$q({\bf 5}.6.3)=6$&$q({\bf 5}.7.3)=6$\\
                &$E$=0.003073      &$E$=0.002486      &$E$=0.001872      &$E$=0.002260      &$E$=0.001596      \\
\\
                &$q({\bf 5}.3.4)=6$&$q({\bf 5}.4.4)=6$&$q({\bf 5}.5.4)=6$&$q({\bf 5}.6.4)=6$&$q({\bf 5}.7.4)=6$\\
                &$E$=0.003071      &$E$=0.001399      &$E$=0.001676      &$E$=0.002497      &$E$=0.001399        \\
\\
                &$q({\bf 5}.3.5)=6$&$q({\bf 5}.4.5)=6$&$q({\bf 5}.5.5)=6$&$q({\bf 5}.6.5)=6$&$q({\bf 5}.7.5)=6$\\
                &$E$=0.002256      &$E$=0.002281      &$E$=0.000729      &$E$=0.002256      &$E$=0.002284       \\
\\
                &$q({\bf 5}.3.6)=6$&$q({\bf 5}.4.6)=6$&                  &$q({\bf 5}.6.6)=6$&$q({\bf 5}.7.6)=6$\\
                &$E$=0.003077      &$E$=0.002503       &                 &$E$=0.003085      &$E$=0.002500       \\
\\
                &$q({\bf 5}.3.7)=6$&$q({\bf 5}.4.7)=6$&                  &$q({\bf 5}.6.7)=6$&$q({\bf 5}.7.7)=6$\\
                &$E$=0.003073      &$E$=0.001650       &                 &$E$=0.003074      &$E$=0.000928        \\
\\
                &$q({\bf 5}.3.8)=6$&$q({\bf 5}.4.8)=6$&                  &$q({\bf 5}.6.8)=6$&$q({\bf 5}.7.8)=6$\\
                &$E$=0.002256      &$E$=0.003096       &                 &$E$=0.001841      &$E$=0.003074        \\
\\
                &$q({\bf 5}.3.9)=6$&$q({\bf 5}.4.9)=6$&                  &$q({\bf 5}.6.9)=6$&$q({\bf 5}.7.9)=6$\\
                &$E$=0.003083      &$E$=0.002281       &                 &$E$=0.003077      &$E$=0.002284         \\
\\
                &$q({\bf 5}.3.10)=6$&$q({\bf 5}.4.10)=6$&                 &$q({\bf 5}.6.10)=6$&$q({\bf 5}.7.10)=6$\\
                &$E$=0.002484       &$E$=0.002486        &                &$E$=0.003085       &$E$=0.001596       \\
\\
                &                   &                   &                 &$q({\bf 5}.6.11)=6$&                  \\
                &                   &                   &                 &$E$=0.002497        &                  \\
\\
                &                   &                   &                 &$q({\bf 5}.6.12)=6$&                  \\
                &                   &                   &                 &$E$=0.002256        &                  \\
\\
                &                   &                   &                 &$q({\bf 5}.6.13)=6$&                  \\
                &                   &                   &                 &$E$=0.001591        &                  \\
\\
                &                   &                   &                 &$q({\bf 5}.6.14)=6$&                  \\
                &                   &                   &                 &$E$=0.003074        &                  \\
\\
                &                   &                   &                 &$q({\bf 5}.6.15)=6$&                  \\
                &                   &                   &                 &$E$=0.002260        &                   \\
\\
\noalign{\smallskip}\hline\noalign{\smallskip}                           
\end{tabular}
\vspace{9mm}
\end{table}

\begin{table}
\centering
\caption{$T-t=0.01.$ The values $E_3^p/(T-t)^5\stackrel{\sf def}{=}E$.}
\label{tab:22}      
\begin{tabular}{p{3.0cm}p{3.0cm}p{2.2cm}}
\hline\noalign{\smallskip}
$I_{(001)T,t}^{(i_1i_2i_3)}$&$I_{(010)T,t}^{(i_1i_2i_3)}$&$I_{(100)T,t}^{(i_1i_2i_3)}$\\
\noalign{\smallskip}\hline\noalign{\smallskip}
\\
$q({\bf 3}.1.b)=6$&$q({\bf 3}.1.c)=4$&$q({\bf 3}.1.d)=2$\\
$E$=0.009425      &$E$=0.009051      &$E$=0.008154  \\
\\
$q({\bf 3}.2.b)=0$&$q({\bf 3}.2.c)=4$&$q({\bf 3}.2.d)=2$\\
$E$=0.000007       &$E$=0.000049     &$E$=0.000147       \\
\\
$q({\bf 3}.3.1.b)=6$&$q({\bf 3}.3.1.c)=4$&$q({\bf 3}.3.1.d)=2$\\
$E$=0.004361        &$E$=0.006366        &$E$=0.004142   \\
\\
$q({\bf 3}.3.2.b)=6$&$q({\bf 3}.3.2.c)=4$&$q({\bf 3}.3.2.d)=2$\\
$E$=0.005044        &$E$=0.002731        &$E$=0.004778   \\
\\
$q({\bf 3}.3.3.b)=6$&$q({\bf 3}.3.3.c)=4$&$q({\bf 3}.3.3.d)=2$\\
$E$=0.009557        &$E$=0.009152        &$E$=0.007963   \\
\\
\noalign{\smallskip}\hline\noalign{\smallskip}
\end{tabular}
\vspace{7mm}
\end{table}

\centerline{\bf Strong Taylor--Ito scheme with convergence order 2.5
(\ref{al4})}

\vspace{4mm}

$$
\frac{(T-t)^2}{2}\left(\frac{1}{2}-\sum_{i=1}^q
\frac{1}{4i^2-1}\right)\le C (T-t)^6,
$$

\vspace{4mm}

$$
(T-t)^3\left(\frac{1}{6}-\frac{1}{64}\sum_{j_1,j_2,j_3=0}^{q_1}
(2j_1+1)(2j_2+1)(2j_3+1)
\left(\bar C_{j_3j_2j_1}^{000}\right)^2\right)\le C(T-t)^6,
$$

\vspace{4mm}

$$
(T-t)^4\left(\frac{1}{4}-\frac{1}{64}\sum_{j_1,j_2=0}^{q_2}
(2j_1+1)(2j_2+1)
\left(\bar C_{j_2j_1}^{01}\right)^2\right)\le C(T-t)^6,
$$

\vspace{4mm}

$$
(T-t)^4\left(\frac{1}{12}-\frac{1}{64}\sum_{j_1,j_2=0}^{q_2}
(2j_1+1)(2j_2+1)
\left(\bar C_{j_2j_1}^{10}\right)^2\right)\le C(T-t)^6,
$$

\vspace{4mm}

$$
(T-t)^4\left(\frac{1}{24}-\frac{1}{256}\sum_{j_1,\ldots,j_4=0}^{q_3}
(2j_1+1)(2j_2+1)(2j_3+1)(2j_4+1)
\left(\bar C_{j_4\ldots j_1}^{0000}\right)^2\right)\le C(T-t)^6,
$$

\vspace{4mm}

$$
(T-t)^4\left(\frac{1}{120}-\frac{1}{32^2}\sum_{j_1,\ldots,j_5=0}^{q_{4}}
(2j_1+1)(2j_2+1)(2j_3+1)(2j_4+1)(2j_5+1)
\left(\bar C_{j_5\ldots j_1}^{00000}\right)^2\right)\le C(T-t)^6,
$$

\vspace{4mm}

$$
(T-t)^5\left(\frac{1}{10}-\frac{1}{256}\sum_{j_1,j_2,j_3=0}^{q_5}
(2j_1+1)(2j_2+1)(2j_3+1)
\left(\bar C_{j_3j_2j_1}^{001}\right)^2\right)\le C(T-t)^6,
$$

\vspace{4mm}

$$
(T-t)^5\left(\frac{1}{20}-\frac{1}{256}\sum_{j_1,j_2,j_3=0}^{q_6}
(2j_1+1)(2j_2+1)(2j_3+1)
\left(\bar C_{j_3j_2j_1}^{010}\right)^2\right)\le C(T-t)^6,
$$

\vspace{4mm}

$$
(T-t)^5\left(\frac{1}{60}-\frac{1}{256}\sum_{j_1,j_2,j_3=0}^{q_{7}}
(2j_1+1)(2j_2+1)(2j_3+1)
\left(\bar C_{j_3j_2j_1}^{100}\right)^2\right)\le C(T-t)^6.
$$

\vspace{8mm}

Note that in this paper we use
the database with 270,000 exactly calculated
Fourier--Legendre coefficients
(see \cite{Kuz-Kuz}, \cite{Mikh-1}) for detail).

\begin{table}
\centering
\caption{Comparison of the conditions (\ref{eeee1}), (\ref{eeee2}).}
\label{tab:23}      
\begin{tabular}{p{1.5cm}p{1.5cm}p{1.5cm}p{1.5cm}p{1.5cm}p{1.5cm}p{1.5cm}}
\hline\noalign{\smallskip}
$T-t$&$2^{-1}$&$2^{-2}$&$2^{-3}$&$2^{-4}$&$2^{-5}$&$2^{-6}$\\
\noalign{\smallskip}\hline\noalign{\smallskip}
$p$&$0$&$0$&$1$&$2$&$4$&$8$\\
$(p+1)^3$&1&1&8&27&125&729\\
$p'$&$1$&$3$&$6$&$12$&$24$&$48$\\
$(p'+1)^3$&8&64&343&2197&15625&117649\\
\noalign{\smallskip}\hline\noalign{\smallskip}
\end{tabular}
\vspace{2mm}
\end{table}

\begin{table}
\centering
\caption{Comparison of the conditions (\ref{eeee3}), (\ref{eeee4}).}
\label{tab:24}      
\begin{tabular}{p{1.5cm}p{1.5cm}p{1.5cm}p{1.5cm}p{1.5cm}p{1.5cm}p{1.5cm}}
\hline\noalign{\smallskip}
$T-t$&$2^{-1}$&$2^{-3/2}$&$2^{-2}$&$2^{-5/2}$&$2^{-3}$&$2^{-7/2}$\\
\noalign{\smallskip}\hline\noalign{\smallskip}
$q$&$0$&$0$&$0$&$0$&$0$&$0$\\
$(q+1)^4$&1&1&1&1&1&1\\
$q'$&$3$&$4$&$6$&$9$&$12$&$17$\\
$(q'+1)^4$&256&625&2401&10000&28561&104976\\
\noalign{\smallskip}\hline\noalign{\smallskip}
\end{tabular}
\vspace{8mm}
\end{table}

\begin{table}
\centering
\caption{Comparison of the conditions (\ref{eeee5}), (\ref{eeee6}).}
\label{tab:25}      
\begin{tabular}{p{1.5cm}p{1.5cm}p{1.5cm}p{1.5cm}p{1.5cm}p{1.5cm}p{1.5cm}}
\hline\noalign{\smallskip}
$T-t$&$2^{-1/8}$&$2^{-1/4}$&$2^{-1/2}$&$2^{-3/4}$&$2^{-1}$\\
\noalign{\smallskip}\hline\noalign{\smallskip}
$r$&$0$&$0$&$0$&$0$&$0$\\
$(r+1)^5$&$1$&$1$&$1$&$1$&$1$\\
$r'$&$1$&$2$&$3$&$4$&$5$\\
$(r'+1)^5$&$32$&$243$&$1024$&$3125$&$7776$\\
\noalign{\smallskip}\hline\noalign{\smallskip}
\end{tabular}
\vspace{8mm}
\end{table}

Let us consider the minimal natural numbers $p, p', q, q', r, r'$
satisfying the conditions

\vspace{1mm}
\begin{equation}
\label{eeee1}
(T-t)^3\left(\frac{1}{6}-\frac{1}{64}\sum_{j_1,j_2,j_3=0}^{p}
(2j_1+1)(2j_2+1)(2j_3+1)
\left(\bar C_{j_3j_2j_1}^{000}\right)^2\right)\le (T-t)^4,
\end{equation}

\vspace{4mm}

\begin{equation}
\label{eeee2}
6(T-t)^3\left(\frac{1}{6}-\frac{1}{64}\sum_{j_1,j_2,j_3=0}^{p'}
(2j_1+1)(2j_2+1)(2j_3+1)
\left(\bar C_{j_3j_2j_1}^{000}\right)^2\right)\le (T-t)^4,
\end{equation}

\vspace{4mm}

\begin{equation}
\label{eeee3}
(T-t)^4\left(\frac{1}{24}-\frac{1}{256}\sum_{j_1,\ldots,j_4=0}^{q}
(2j_1+1)(2j_2+1)(2j_3+1)(2j_4+1)
\left(\bar C_{j_4\ldots j_1}^{0000}\right)^2\right)\le (T-t)^5,
\end{equation}

\vspace{4mm}

\begin{equation}
\label{eeee4}
24(T-t)^4\left(\frac{1}{24}-\frac{1}{256}\sum_{j_1,\ldots,j_4=0}^{q'}
(2j_1+1)(2j_2+1)(2j_3+1)(2j_4+1)
\left(\bar C_{j_4\ldots j_1}^{0000}\right)^2\right)\le (T-t)^5,
\end{equation}

\vspace{4mm}

$$
(T-t)^4\left(\frac{1}{120}-\frac{1}{32^2}\sum_{j_1,\ldots,j_5=0}^{r}
(2j_1+1)(2j_2+1)(2j_3+1)(2j_4+1)(2j_5+1)
\left(\bar C_{j_5\ldots j_1}^{00000}\right)^2\right)\le 
$$

\begin{equation}
\label{eeee5}
\le (T-t)^6,
\end{equation}

\vspace{4mm}

$$
120(T-t)^4\left(\frac{1}{120}-\frac{1}{32^2}\sum_{j_1,\ldots,j_5=0}^{r'}
(2j_1+1)(2j_2+1)(2j_3+1)(2j_4+1)(2j_5+1)
\left(\bar C_{j_5\ldots j_1}^{00000}\right)^2\right)\le 
$$

\begin{equation}
\label{eeee6}
\le (T-t)^6,
\end{equation}

\vspace{6mm}
\noindent
where the inequalities (\ref{eeee2}), (\ref{eeee4}), (\ref{eeee6}) 
are particular cases
of the formula (\ref{sme}) for $r=3, 4,$ and $5$.

In Tables 23--25, we can see the numerical comparison
of the conditions (\ref{eeee1}), (\ref{eeee3}), (\ref{eeee5})
with the conditions (\ref{eeee2}), (\ref{eeee4}), (\ref{eeee6}),
respectively.
Obviously, the conditions (\ref{eeee1}), (\ref{eeee3}), (\ref{eeee5})
(i.e. conditions without the multiplier factors $3!,$ $4!,$ and $5!$)
essentially reduce the calculation costs for
the mean-square approximations of iterated Ito
stochastic integrals 

\vspace{2mm}
$$
I_{(000)T,t}^{(i_{1}i_{2}i_{3})},\ \ \ 
I_{(0000)T,t}^{(i_{1}i_{2}i_{3}i_4)},\ \ \ 
I_{(00000)T,t}^{(i_{1}i_{2}i_{3}i_4i_5)}\ \ \ 
(i_1,i_2,i_3,i_4,i_5=1,\ldots,m).
$$

\end{document}